\documentclass{tran-l}

\newtheorem{theorem}{Theorem}[section]

\usepackage{epsfig,psfrag}
\theoremstyle{definition}

\theoremstyle{remark}
\newtheorem{remark}[theorem]{Remark}

\numberwithin{equation}{section}

\swapnumbers \theoremstyle{plain}

\newtheorem{thm}{Theorem}[section]

\newtheorem{lem}[thm]{Lemma}
\newtheorem{cor}[thm]{Corollary}

\newtheorem{prop}[thm]{Proposition}

\theoremstyle{definition}

\theoremstyle{remark}

\newcommand{\T}{\mathcal T}
\newcommand{\C}{\mathcal C}

\renewcommand{\P}{{\mathcal P}}
\renewcommand{\S}{{\mathcal S}}
\newcommand{\bdy}{\partial}
\newcommand{\bbb}{\mathbb}
\newcommand{\rp}{\mathbb{R}P^3}
\newcommand{\rpp}{\mathbb{R}P^2}

\newcommand{\td}{\tilde}
\newcommand{\open}[1]{\stackrel{\circ}{#1}}
\newcommand{\ve}{\varepsilon}

\newcommand{\abs}[1]{\lvert#1\rvert}

\begin{document}

\title{$\mathbf{0}$--efficient triangulations of
$3$--manifolds}

\author{William Jaco}
\address{Department of Mathematics, Oklahoma State University,
Stillwater, OK 74078}

\email{jaco@math.okstate.edu}
\thanks{The first author was partially supported by NSF Grants
 DMS9704833 and DMS9971719,
The Grayce B. Kerr Foundation, The American Institute of Mathematics (AIM), and The Visiting
Research Scholar Program at University of Melbourne (Australia)}

\author{J.~Hyam Rubinstein}
\address{Department of Mathematics and Statistics,
University of Melbourne, Parkville, VIC 3052, Australia}
\email{rubin@maths.unimelb.edu.au}
\thanks{The second author was partially supported by The Australian Research Council, The
Grayce B. Kerr Foundation, Stanford University and The American Institute of Mathematics
(AIM)}


\date{January 11, 2001}


\begin{abstract}
$0$--efficient triangulations of $3$--manifolds are defined and studied.  It is shown
that any triangulation of a closed, orientable, irreducible $3$--manifold $M$ can be modified
to a $0$--efficient triangulation or $M$ can be shown to be one of the manifolds $S^3,
\rp$ or $L(3,1)$. Similarly, any triangulation of a compact, orientable, irreducible,
$\bdy$--irreducible $3$--manifold can be modified to a $0$--efficient triangulation. The
notion of a $0$--efficient ideal triangulation is defined. It is shown if $M$ is a compact,
orientable, irreducible, $\bdy$--irreducible $3$--manifold having no essential annuli and
distinct from the $3$--cell, then
$\stackrel{\circ}{M}$ admits an ideal triangulation; furthermore, it is
shown that any ideal triangulation of such a $3$--manifold can be modified to a
$0$--efficient ideal triangulation. A $0$--efficient triangulation of a
closed manifold has only one vertex or the manifold is $S^3$ and the triangulation has
precisely two vertices. $0$--efficient triangulations of $3$--manifolds with boundary,
 and distinct from the
$3$--cell, have all their vertices in the boundary and then just
one vertex in each boundary component. As tools, we introduce the
concepts of barrier surface and shrinking, as well as the notion
of crushing a triangulation along a normal surface. A number of
applications are given, including an algorithm to construct  an
irreducible decomposition of a closed, orientable $3$--manifold,
an algorithm to construct a maximal collection of pairwise
disjoint, normal $2$--spheres in a closed $3$--manifold, an
alternate algorithm for the $3$--sphere recognition problem,
results on edges of low valence in minimal triangulations of
$3$--manifolds, and a construction of irreducible knots in closed
$3$--manifolds.
\end{abstract}

\maketitle

\section{introduction}\label{sect-intro}
This is the first in a series of papers concerning triangulations
of 3-manifolds. Our aim is to give a theory of {\it efficient}
triangulations of $3$--manifolds, which has applications to
finiteness theorems, knot theory, Dehn fillings, decision
problems, algorithms, computational complexity and Heegaard
splittings. A motivation for these investigations is to achieve
properties of triangulations which have features similar to those
of geometric structures on manifolds. For example, with
triangulations, the analogue of a stable minimal surface is a
normal surface and the analogue of an unstable surface of index
one is an almost normal surface. Our work very strongly uses
normal and almost normal surfaces and is guided by the
similarities to minimal surfaces. An appropriately efficient
triangulation has features similar to those of a complete
hyperbolic structure of finite volume on a $3$-manifold. In the
latter, any essential surface can be homotoped to a (stable)
minimal surface; hence, by the Gauss-Bonnet formula, there can not
be any essential surfaces of nonnegative Euler characteristic. In
a manifold with a $1$--efficient triangulation, we obtain similar
results, limiting embedded normal surfaces with nonnegative Euler
characteristic.

In this paper, we define and study $0$--efficient triangulations of $3$--manifolds. A
triangulation of a closed, orientable $3$--manifold is {\it $0$--efficient} if and only if the
only embedded, normal $2$--spheres are vertex-linking. For a compact, orientable
$3$--manifold with nonempty boundary, a triangulation is {\it $0$--efficient} if and only if
the only properly embedded, normal disks are vertex-linking. We show that a closed
$3$--manifold with a $0$--efficient triangulation is irreducible, not
$\mathbb{R}P^3$ and the triangulation has precisely one-vertex or the manifold is
$S^3$ and the triangulation has precisely two vertices. In the case of a compact
$3$--manifold with boundary, no component of which is a $2$--sphere, a $0$--efficient
triangulation implies that the manifold is irreducible and $\partial$--irreducible and
the triangulation has all vertices in the boundary and has precisely one vertex in each
boundary component.

One of the main theorems in Section \ref{sect-0-eff}, is that
any triangulation of a closed,  orientable, irreducible
$3$--manifold can be modified to a
$0$--efficient triangulation or it can be shown the $3$--manifold is one of $S^3$,
$\mathbb{R}P^3$, or
$L(3,1)$. Similarly, it is shown that any triangulation of a compact, orientable
irreducible, $\bdy$--irreducible  $3$--manifold with nonempty boundary can be
modified to a $0$--efficient triangulation.

 We show that a minimal triangulation of
a closed, orientable, irreducible $3$--manifold is $0$--efficient,
except for $\rp$ and $L(3,1)$, and therefore has just one vertex
or is $S^3$. The $3$--sphere has two distinct, one-tetrahedron
triangulations, both are $0$--efficient but one has two vertices.
For the other exceptions, there are two distinct, two-tetrahedra
(minimal) triangulations of $\rp$, neither are $0$--efficient, and
four distinct, two-tetrahedra (minimal) triangulations of
$L(3,1)$, only one of which is $0$--efficient. Unlike the case of
$2$--manifolds, where one immediately has from Euler
characteristic that a minimal triangulation of a closed
$2$--manifold with non positive Euler characteristic has just one
vertex, there does not seem to be an immediate way to see that a
minimal triangulation of an irreducible $3$--manifold, other than
one of $S^3, \rp$ and $L(3,1)$, has just one vertex.  One-vertex
triangulations are  very interesting and have a number of
applications. In particular, one-vertex triangulations of
$3$-manifolds appear to be very well-suited for normal surface
theory. See \cite{jac-let-rub1, jac-let-rub2, jac-let-rub-sed,
jac-rub4, jac-rub-sed, jac-sed1} for numerous applications of
one-vertex and efficient triangulations.

We also study $0$--efficient (no normal $2$--spheres), ideal triangulations and show that the
interior of any compact, orientable, irreducible, $\bdy$--irreducible $3$--manifold that has
no essential annuli admits an ideal triangulation; furthermore, in such a $3$--manifold,
any ideal triangulation can be modified to a $0$--efficient, ideal triangulation. It follows
that a minimal, ideal triangulation of such a bounded $3$--manifold is $0$--efficient.

Our main technique is to identify a
 normal
$2$-sphere, which is not vertex-linking and which bounds a $3$--cell, and crush  such a $2$--sphere and
the
$3$-cell it bounds to a point. This gives us back our manifold but wrecks havoc with the triangulation.
Most work then is to show that we can recover a ``nicer" triangulation; we must do this with great care
not to add tetrahedra (or new normal $2$--spheres). To achieve this we introduce in Section
\ref{sect-crush} the notion of crushing a triangulation along a normal surface. This has turned out to be
a useful tool throughout our work.

In the sequel to this paper, we show that further modifications of
the triangulation can be performed for closed, orientable
$3$--manifolds which are not only irreducible but also assumed
atoroidal. In such cases, we show that a triangulation can be
obtained so that any embedded, normal torus is of a very special
form or it can be shown the $3$--manifold is $S^3$, a lens space
or a small Seifert fiber space. This property along with
$0$--efficiency is called $1$-efficiency. The reader can
immediately see our motivation if they recall that in Haken's
theory of normal surfaces, there is a finite constructible
collection of fundamental normal surfaces and every normal surface
is a sum of these fundamental surfaces. Hence, if we have control
of  all surfaces with nonnegative Euler characteristic, then we
can expect a bounded number of such sums giving all normal
surfaces of bounded genus, up to isotopy. A similar result works
for almost normal surfaces and this completes the solution of the
Waldhausen Conjecture that a closed, orientable $3$--manifold has
only a finite number (up to homeomorphism and isotopy) of Heegaard
splittings of bounded genus \cite{jac-rub5}. This was shown in the
case of Haken manifolds by K. Johannson \cite{joh}.

One of the main applications in this first paper is an algorithm
to determine the prime decomposition of a $3$-manifold. Given a
closed, orientable $3$--manifold, $M$, we give an algorithm to
write $M$ as a connected sum decomposition of $3$--manifolds where
it is known that each factor is either, $S^3$, $S^2\times S^1$,
$\mathbb{R}P^3$, $L(3,1)$, or has a $0$--efficient triangulation.
Now, a $0$--efficient triangulation is precisely the environment
for implementation of the $3$--sphere recognition algorithm,
\cite{rub}. Hence, we arrive at the prime decomposition of the
given $3$--manifold $M$. Ben Burton and David Letscher have
written a program to implement this algorithm. We observe certain
methods which have implications toward the complexity of these
algorithms. In particular, we obtain results similar to those
announced by A. Casson who observed that the algorithm for
achiving a prime decomposition can be implemented in time
essentially of order $p(t)\mathcal{O}( 3^t)$, where $t$ is the
number of tetrahedra in a given triangulation for $M$ and $p(t)$
is a low degree polynomial in $t$.

We provide the background  from normal and almost normal surface
theory needed in this paper in Section \ref{sect-triang}, along
with our generalized idea of triangulations. However, we assume
the reader has some familiarity with normal and almost normal
surface theory. In Section \ref{sect-barriers}, we define and
provide a study of barrier surfaces and shrinking within the
context of normal and almost normal surface theory. These ideas
come from geometric surface theory and are consistent with the
usage there. More complete developments of barriers and shrinking
can be found in \cite{jac-rub1, jac-let-rub1}.

We would like to thank Ben Burton, Dave Letscher and Eric Sedgwick
for many helpful conversations throughout the development of this
theory. Also, we want to acknowledge Andrew Casson for helpful
comments and ideas on this project. In particular,  Andrew Casson
announced a beautiful program to prove geometrization for
irreducible, atoroidal $3$-manifolds with tori boundaries at a
meeting in Montreal in 1995. His method involved a version of
efficient {\it ideal} triangulations and became our model for
efficient triangulations of closed and bounded $3$--manifolds.
Finally, we wish to thank Thomas Faulkenberry, who has been so
helpful in drafting many of the figures which appear in this
paper.

\section{triangulations, cell-decompositions and normal surfaces}\label{sect-triang}
In this section we give some basic definitions and results needed
for the work presented in this paper. Our approach to
triangulations is more general than the notion one typically
finds. We develope a broad study of triangulations from this point
of view in \cite{jac-let-rub1}. At various places in our proofs
and constructions, we use cell-decompositions having more general
cells than triangulations; however, a major feature of our methods
is that these various cell-decompositions are, themselves, quite
nice. We shall assume the reader has a basic understanding of
normal surface theory.  The references \cite{jac-ree, jac-rub4,
jac-tol} are sources to review normal surface theory from our
point of view. While we use very little normal surface theory
here, it is an indispensable part of our program and serves to
organize what might otherwise be quite messy.

\subsection{Triangulations and cell-decompositions.} We assume most readers are familiar with
cell-decompositions; however, since we are using them in a slightly broader sense than usual, we
include a brief discussion with some definitions.

Let
$\boldsymbol{\Delta} = \{\td{\Delta}_1,\ldots,\td{\Delta}_t\}$ be a pairwise-disjoint collection
of oriented,
 compact, convex, linear  cells. Suppose
$\Phi$ is a family of affine isomorphisms pairing faces of the
cells in $\boldsymbol{\Delta}$ so that if $\phi\in\Phi$, then
$\phi$ is an orientation-reversing affine isomorphism from a face
$\sigma_i\in\td{\Delta}_i$ to a face $\sigma_j\in\td{\Delta}_j$,
possibly $i = j$.   We use $\boldsymbol{\Delta}/\Phi$ to denote
the space obtained from the disjoint union of the $\td{\Delta}_i$
by setting $x\in\td\sigma_i$ equal to $\phi(x)\in\td\sigma_j$,
with the identification topology. Then $\boldsymbol{\Delta}/\Phi$
is a $3$--manifold, except possibly at the images of the vertices
of the $\td{\Delta}_i$. (In a completely general setting, the
identification space $\boldsymbol{\Delta}/\Phi$ may not be a
$3$--manifold at the image of the centers of some edges, as well
as the images of the vertices; however, we have avoided this
problem by orienting the $\td{\Delta}_i$ and choosing the affine
isomorphisms $\phi\in\Phi$ orientation-reversing.) We collect all
this information into a single symbol $\T$ and call $\T$ a {\it
cell-decomposition} of $\boldsymbol{\Delta}/\Phi$; in this case,
we also use
 just $\abs{\T}$ to denote the space $\boldsymbol{\Delta}/\Phi$. A {\it cell (tetrahedron), face,
edge}, or {\it vertex} in this cell decomposition is,
respectively, the image of a cell (tetrahedron), face, edge, or
vertex from the collection $\boldsymbol{\Delta} =
\{\td{\Delta}_1,\ldots,\td{\Delta}_t\}$. We will denote the image
of the faces by $\T^{(2)}$, the image of the edges by $\T^{(1)}$
and the image of the vertices by $\T^{(0)}$. We call $\T^{(i)}$
the {\it $i$--skeleton of $\T$}; but, generally, we just refer to
these as the faces, edges or vertices of $\T$. We will denote the
image of $\td{\Delta}_i$ by $\Delta_i$ and call $\td{\Delta}_i$
the {\it lift of $\Delta_i$}. A cell is the quotient of a unique
cell and a face is the quotient of one or two faces; edges and
vertices may be the quotient of a number of edges or vertices,
respectively. While the cells are not necessarily embedded, the
interior of each cell is embedded. We define the {\it order} or
{\it valence of an edge} $e$ of $\T$ to be the number of edges in
the collection $\boldsymbol{\Delta} = \{\td{\Delta}_1,\ldots$
$\ldots,\td{\Delta}_t\}$, which are identified to $e$.  If the
link (here we mean the boundary of a small regular neighborhood
and not the combinatorial link) of each vertex  is either a
$2$--sphere or a $2$--cell, then the underlying point set is an
oriented $3$--manifold, possibly with boundary, and we say $\T$ is
a {\it cell-decomposition} of the $3$--manifold  $M = \abs{\T}$.
If each cell in $\boldsymbol{\Delta}$ is a tetrahedron, then we
say $\T$ is a {\it triangulation} of the $3$--manifold $M =
\abs{\T}$. In the literature, one is likely to find this notion of
a triangulation referred to as a {\it pseudo}-triangulation and
the term triangulation reserved to mean that tetrahedra are
embedded and if two simplices meet at all, then they meet in a
face of each. Our cells (tetrahedra) are not necessarily embedded;
however, they are embedded on the interiors of each cell
(tetrahedron), face and edge  and the intersection of two cells is
a union of sub-cells of each. If the link of some vertex is a
closed surface, distinct from the $2$--sphere, we say $\T$ is an
{\it ideal cell-decomposition} (or {\it ideal triangulation}) of
the $3$--manifold $M = \abs{\T}\setminus\abs{\T^{(0)}}$. In this
case the vertices of $\T$ are called {\it ideal vertices } and the
{\it index of an ideal vertex} is the genus of its linking
surface. In some cases where we have a mix of genus zero and
higher genus ideal vertices, we include the genus zero vertices
into our manifold. Generally, however, we don't have genus zero
vertices in an ideal triangulation.

 Similarly, we use cell-decompositions of surfaces. If
$\boldsymbol{\sigma} = \{\td{\sigma}_1,\ldots,\td{\sigma}_n\}$ is
a pairwise disjoint collection of compact, convex, planar polygons
and $\Psi$ is a family of linear isomorphisms pairing edges of the
polygons in $\boldsymbol{\sigma}$ so that $\psi\in\Psi$, then
$\psi$ is a linear isomorphism of  an edge $e_i$ of
$\td{\sigma}_i$ to an edge $e_j$ of $\td{\sigma}_j$, possibly $i =
j$. We use $\boldsymbol{\sigma}/\Psi$ to denote the space obtained
from the disjoint union of the $\td{\sigma}_i$ by setting
$x\in\td{e}_i$ equal to $\psi(x)\in\td{e}_j$, with the
identification topology. We have that $\boldsymbol{\sigma}/\Psi$
is always a $2$--manifold, possibly nonorientable. In this
situation, we say we have a cell-decomposition of the
$2$--manifold $\boldsymbol{\sigma}/\Psi$. If each $\td{\sigma}_i$
is a triangle, we say we have a {\it triangulation} of the
$2$--manifold $\boldsymbol{\sigma}/\Psi$. Similar to the case for
$3$--manifolds, in a cell-decomposition of a $2$--manifold our
cells are not embedded; however, the open cells are embedded.

\begin{figure}[htbp]
            \psfrag{a}{$v_0$}
            \psfrag{b}{$v_1$}
            \psfrag{c}{$v_2$}
            \psfrag{1}{($1$)}
            \psfrag{2}{($2$)}
            \psfrag{3}{($3$)}
            \psfrag{4}{($4$)}
            \psfrag{5}{($5$)}
            \psfrag{6}{($6$)}
            \psfrag{7}{($7$)}
            \psfrag{8}{($8$)}
        \vspace{0 in}
        \begin{center}
        \epsfxsize=4 in
        \epsfbox{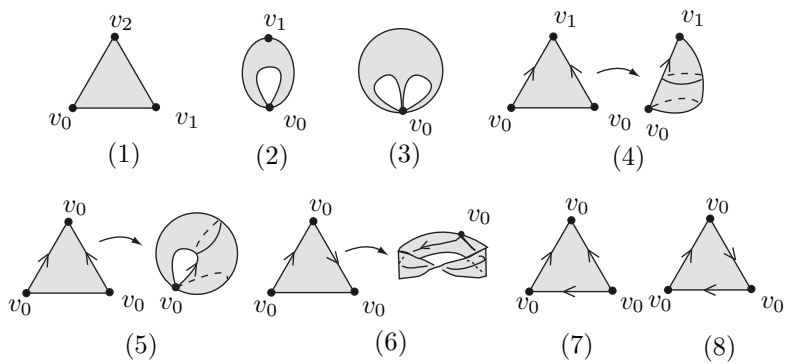}
        \caption{Possible configurations for faces in a
        triangulation.}
        \label{f-faces}
        \end{center}
        \end{figure}

\begin{figure}[htbp]
            \psfrag{a}{$v_0$}
            \psfrag{b}{$v_1$}
            \psfrag{c}{$v_2$}
            \psfrag{d}{$v_3$}
            \psfrag{1}{($1$) $\mathbb{B}^3$}
            \psfrag{2}{($2$) $\mathbb{B}^3$}
            \psfrag{3}{($3$) $\mathbb{B}^2\times S^1$}
            \psfrag{4}{($4$) $S^3$}
            \psfrag{5}{($5$) $S^3$}
            \psfrag{6}{($6$) $L(4,1)$}
            \psfrag{7}{($7$) $L(5,2)$}
            \psfrag{e}{$e$}
        \vspace{0 in}
        \begin{center}
        \epsfxsize=4 in
        \epsfbox{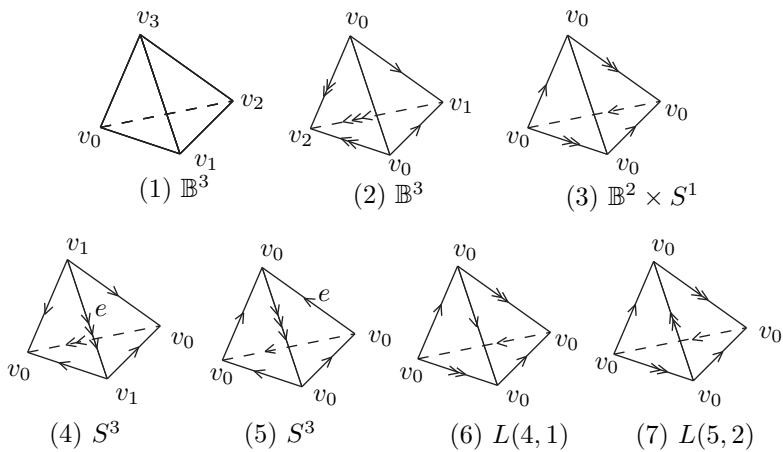}
        \caption{One-tetrahedron triangulations of orientable
$3$--manifolds.}
        \label{f-tetra}
        \end{center}
\end{figure}

In our cell-decompositions, an edge can be a simple closed curve,
an edge in a cell with end points (vertices) identified. In fact,
we will be working toward having triangulations with just one
vertex. In this case, every edge is a simple closed curve. Faces
can take on some interesting configurations. In Figure
\ref{f-faces} we give the possibilities. In Figure \ref{f-faces},
parts (4) and (5), we have two edges identified to give faces
which are {\it cones} (the latter is a {\it pinched cone}); in
(6), we have a face which is a M\"obius band; and in (7) and (8),
we have all three edges identified, giving in (7) the classical
{\it dunce hat} (see for example Figure \ref{f-tetra} (5)) and in
(8) a spine for $L(3,1)$ (see for example the classical
presentation for $L(3,1)$ in Figure \ref{f-not-0-eff}). In the
case of a triangulation of a $2$--manifold, we only have those
possibilities (1), (2), (3), (4) and (6) in Figure \ref{f-faces}
for the triangles. As for tetrahedra, we have in Figure
\ref{f-tetra} the seven distinct identifications of a single
tetrahedron to give an orientable $3$--manifold.  In Figure
\ref{f-tetra}, the edge $e$ in part (4) has order $1$ and in part
(5) $e$ has order $5$.

\subsection {Normal surfaces.}
  A cell-decomposition (or ideal cell-decomposition) of a $3$--manifold distinguishes certain
surfaces, the normal surfaces relative to that cell-decomposition.
They are analogous to stable (minimal) surfaces in geometric
analysis and are determined by how they meet the cells of the
cell-decomposition. However, since our cells are not embedded, we
provide the definition of a normal surface in this setting.
Furthermore, in this paper, we, generally, only use normal surface
theory in triangulations. Even when we use more general cell
decompositions (truncated tetrahedra, truncated prisms, regions
between parallel quads and triangles), we restrict our
considerations to very special types of normal surfaces.

 If
$\td{\Delta}$ is a compact, convex, linear cell and $\td\sigma$ is
a face of $\td{\Delta}$, we say a spanning arc in $\td\sigma$ is a
{\it normal arc} if its end points are in distinct edges of
$\td\sigma$. A {\it normal curve} in the boundary of a compact,
convex, linear cell is a curve which meets each face in a
collection of normal arcs. The elementary components of normal
surface theory are the normal disks in the cells of the
cell-decomposition. We call a properly embedded disk in a compact,
convex, linear cell a {\it normal disk} if its boundary is a
normal curve and it meets no edge more than once. If $\T$ is a
cell-decomposition of the manifold $M$, then an isotopy of $M$ is
called a {\it normal isotopy} if it is invariant on the cells,
faces, edges and vertices of $\T$. Up to normal isotopy there are
only finitely many equivalence classes of normal disks in a
compact, convex, linear cell; these are called {\it normal disk
types}.

 Suppose
$\T$ is a cell-decomposition or ideal cell-decomposition of the $3$--manifold $M$ and $S$ is a
properly embedded surface transverse to the $2$--skeleton of $\T$. Suppose
$c$ is a component of $S$ in the cell $\Delta_i$. Then $c$ is the image of a properly
embedded surface,
$\td{c}$ in
$\td{\Delta}_i$. We will call $\td{c}$ the {\it lift of $c$}.

\begin{figure}[htbp]
            \psfrag{a}{\begin{tabular}{c}
            normal\\
            triangles\\
            \end{tabular}}
            \psfrag{b}{\begin{tabular}{c}
            normal\\
           quad\\
            \end{tabular}}
            \psfrag{c}{\begin{tabular}{c}
            normal\\
            quad\\
            \end{tabular}}
            \psfrag{d}{\begin{tabular}{c}
            normal triangle\\
           and hexagon\\
            \end{tabular}}
            \psfrag{e}{\begin{tabular}{c}
            normal\\
            octagon\\
            \end{tabular}}
            \psfrag{f}{\begin{tabular}{c}
            normal triangle\\
           and a normal quad\\
            \end{tabular}}

        \vspace{0 in}
        \begin{center}
        \epsfxsize=4 in
        \epsfbox{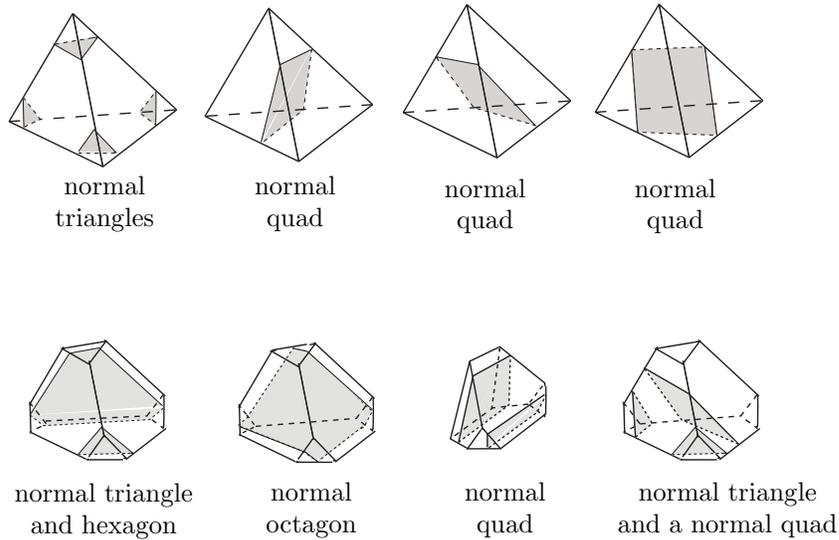}
        \caption{Examples of normal disk types (normal triangles and normal quads for a triangulation).}
        \label{f-normal}
        \end{center}
\end{figure}

\begin{figure}[htbp]
            \psfrag{a}{a}
            \psfrag{b}{b}
            \psfrag{c}{c}
            \psfrag{d}{d}
            \psfrag{e}{e}
            \psfrag{f}{f}
            \psfrag{1}{($1$) vertex-linking $S^2$ in $S^3$}
            \psfrag{2}{($2$) Heegaard torus in $S^3$}
            \psfrag{3}{($3$) M\"obius band in $\mathbb{B}^2\times S^1$}
            \psfrag{4}{ $S^2$}

        \vspace{0 in}
        \begin{center}
        \epsfxsize=4 in
        \epsfbox{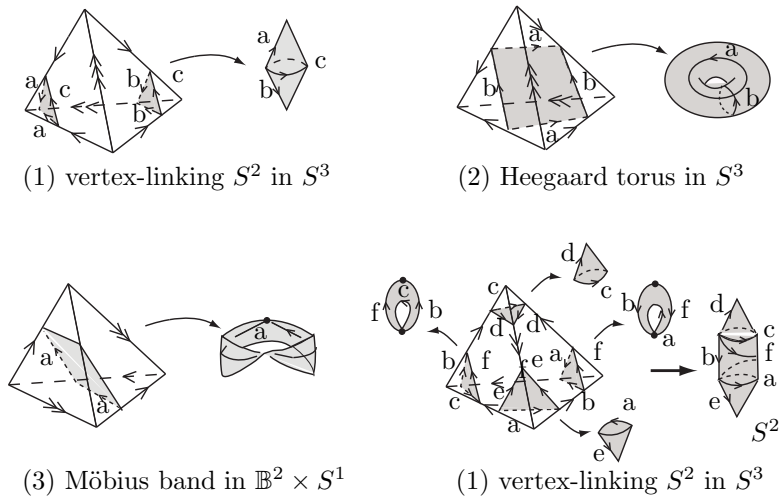}
        \caption{Examples of normal surfaces in triangulated
$3$--manifolds.}
        \label{f-normal-quotient}
        \end{center}
\end{figure}

Now,  if $\T$ is a cell-decomposition of  the $3$--manifold $M$, we say a surface $F$ is a {\it
normal surface} in $M$ (with respect to
$\T$) if $F$ meets each cell of $\T$ in the images of a collection of normal
disks in the cells of $\boldsymbol{\Delta} =
\{\td{\Delta}_1,\ldots,\td{\Delta}_n\}$. That is, the surface $F$ is a normal surface if and
only if the lift
 of every component of $F$ in a cell of $\boldsymbol{\Delta}$
 is a normal disk. The elementary components of normal surface theory for triangulations are the
{\it normal triangles} and {\it normal quadrilaterals} ({\it
normal quads}) in a tetrahedron. The normal triangles and  normal
quads are shown in Figure \ref{f-normal}. There are four types of
normal triangles and three types of normal quads in each
tetrahedron (no identification). Also, in Figure \ref{f-normal},
we give some examples of normal disks in cells which are truncated
tetrahedra and truncated-prisms. In a truncated tetrahedra,
besides normal triangles and normal quads, one may have other
normal disks. For example, we show a normal hexagon and a normal
octagon in a truncated tetrahedra. We also show some normal quads
in a truncated prism and normal triangles and normal quads in a
truncated tetrahedron, which are also normal in the tetrahedron
before truncation. In Figure \ref{f-normal-quotient}, we show a
$2$--sphere made up of two normal triangles (Figure
\ref{f-normal-quotient}(1)), a normal torus made up of a normal
quad (Figure \ref{f-normal-quotient}(2)), and a M\"obius band made
up of a single normal quad (Figure \ref{f-normal-quotient}(3)). In
a tetrahedron in $M$, after identifications, a normal triangle can
take on any of the possible orientable identifications of a
triangle shown in Figure \ref{f-faces}, Parts (1)--(6). (A normal
triangle can identify to a M\"obius band but then the resulting
identification space is not a manifold at the associated vertex.
It can be made into an ideal triangulation of a manifold; but the
manifold is nonorientable.) Similarly, a normal quadrilateral,
after identification, can take on numerous forms, including some
which are nonorientable. The normal disk types give a normal
surface $F$ a cell-decomposition made up of normal quads and
normal triangles; we call this the {\it cell-decomposition induced
on $F$}(or the induced cell-decomposition). Finally, we note that
an embedded normal surface must be properly embedded.

If $\T$ is a cell-decomposition, as we noticed above, there are
only finitely many normal isotopy classes of normal disks. If $n$
is the number of normal isotopy classes of normal disks in $\T$
and we select an ordering of these normal disk types, say
$\mathit{d}_1,\ldots,\mathit{d}_n$ , then a normal isotopy class
of a normal surface has a parameterization as an $n$--tuple of
non-negative integers $(x_1,\dots,x_n)$ in $\mathbb R^n$, where
$x_i$ is the number of elementary disks of type $\mathit{d}_i,
1\leq i\leq n$. For a triangulation there are four normal triangle
types and three normal quad types in each tetrahedra; so, $n = 7t$
for a triangulation $\T$ with $t$ tetrahedra.   Associated with
the cell-decomposition $\T$ is a system of linear equations.
Non-negative integer solutions to this system give a
parameterization of the normal isotopy classes of normal surfaces.
We add the condition $x_i \geq 0, \forall i$, and obtain a cone in
the positive-orthant of $\mathbb R^n$. We will refer to this cone
as the {\it solution cone} and write $\mathcal{S}(M,\T)$. There
are a finite number of nonnegative integer lattice points
(solutions) in the solution cone so that any integer lattice point
in the cone can be written as  a finite linear combination of
these using only nonnegative integer coefficients. These are a
Hilbert basis for the nonnegative integer lattice points in the
solution cone and are called {\it fundamental solutions} after
Haken \cite{haken1}. If, in addition to the above equations, we
add the equation $\sum_1^{n} x_i = 1$, we determine a compact,
convex linear cell.  The rational points in this cell correspond
to projective classes of normal isotopy classes of normal surfaces
in $\mathcal{S}(M,\T)$. We denote this compact, convex linear cell
by $\P(M,\T)$ and call it the {\it projective solution space} (of
$\mathcal{S}(M,\T)$). If $F$ is a normal surface in $M$, we do not
distinguish and let $F$ denote not only the surface $F$ but its
normal isotopy class, and its representation as an $n$--tuple in
$\mathbb R^n$. We denote the projective class of $F$ by $\bar{F}
\in \P(M,\T)$. The {\it carrier} of a normal surface $F$ is the
unique minimal face of $\P(M,\T)$ which contains $\bar{F}$ and is
denoted $\mathcal{C}(F)$.  If the projective representation of the
normal isotopy class of the normal surface $F$ is at a vertex of
the projective solution space, we say $F$ is a {\it vertex
solution}. A vertex solution is a fundamental solution but the
converse is not necessarily true. We have  algebraic
characterizations of fundamental and vertex solutions. A solution
$F$ is a fundamental solution if and only if it can {\it not} be
written as a nontrivial sum $F = X + Y$; a solution $F$ is a
vertex solution if and only $kF$ can {\it not} be written as a
nontrivial sum $kF = nX + mY$ for some positive integers $k,n$ and
$m$.

In a triangulation $\T$, a normal surface is embedded if and only
if it does not meet a tetrahedron in more than one quad type. This
is sometimes call the {\it quadrilateral condition} or the {\it
compatibility condition}. Any normal surface determines a unique
nonnegative, integer $n$--tuple in $\mathbb R^n$; but, this
correspondence is not one-one. However, if an integer lattice
point represents an embedded normal surface, then there is a
unique such embedded, normal surface associated with that lattice
point. One can realize a face (or a cone) of embedded normal
surfaces by adding the conditions that two of the quadrilateral
types in each tetrahedron are zero. There are a maximum $3^t$ such
faces. These contain the representations of all embedded normal
surfaces.

In a cell-decomposition or an ideal cell-decomposition of the $3$--manifold $M$, the boundary of
a small regular neighborhood of a vertex is normally isotopic to a normal surface,   each
component of which is called a {\it vertex-linking surface}. If $\T$ is a cell-decomposition of
the
$3$--manifold $M$, a vertex-linking surface is either a disk (the vertex is in $\bdy M$) or a
$2$--sphere (the vertex is in
$\stackrel{\circ}{M}$, the interior of $M$). If $\T$ is an ideal cell-decomposition and $v$ is an
ideal vertex, a vertex-linking surface about $v$ is a closed, orientable $2$--manifold
possibly having genus $g\geq 1$, where $g$ is the index of the vertex $v$. In this case we
sometimes refer to the vertex-linking surface as a {\it surface-at-infinity}. If $\T$ is a
triangulation or an ideal triangulation, the entire collection of elementary triangles form an
embedded, normal surface, each component of which is a vertex-linking surface. The elementary
disk types in a vertex-linking surface of a more general cell-decomposition do not allow such a
simple combinatorial description. The vertex-linking surfaces give examples of normal surfaces.

If $S$ is a properly embedded surface in a $3$--manifold $M$ and $N(S)$ is a small regular
neighborhood of
$S$, the manifold $M' = M\setminus\stackrel{\circ}{N}(S)$ is  said to be obtained from $M$
by {\it splitting along $S$}. If $S$ is one-sided in $M$, then $S$ is nonorientable and there is
a  copy, say $S'$, of the orientable double cover of $S$ in $\bdy M'$. If $S$ is two-sided, then
there are two homeomorphic copies,
$S'$ and
$S''$, of
$S$ in
$\bdy M'$.  They are in
the same component of
$M'$ if and only if
$S$ does not separate
$M$. If
$S$ is a normal surface, then we choose $N(S)$ in such a way that the components of its
frontier, $S'$ and $S''$, are normally isotopic to $S$ (or if $S$ is one-sided, then $S' = 2S$).

 If $S$ is a normal surface and $M'$
is obtained by splitting $M$ along $S$, there is a natural
cell-decomposition on $M'$. A cell in this decomposition is just a
component of a $\td\Delta_i$ split along the normal disks in the
lifts of $S$ and the face identifications are just the face
identifications for $\T$ restricted to the faces of our new cells.
We call this the {\it induced cell-decomposition} on $M'$, the
manifold obtained by splitting $M$ along $S$. Let $\C$ denote the
induced cell-decomposition on $M'$. A normal disk in a cell in
$\C$, which  misses $S'\cup S''$, also is a normal disk in $\T$;
however, there are, possibly, numerous  distinct disk types of
this sort in $\C$, which all represent the same disk type in $\T$.
Let
$\mathit{d}_{1,1},\ldots,\mathit{d}_{1,n_1},\mathit{d}_{2,1},\ldots,\mathit{d}_{n,1},
\ldots,\mathit{d}_{n,n_n}$ denote the normal disk types in $\C$,
which miss $S'\cup S''$, where
$\mathit{d}_{i,1},\ldots,\mathit{d}_{i,n_i}$ are all normally
isotopic to the normal disk type $\mathit{d}_i$ in $\T$. If $F$ is
a normal surface in $X$, which misses $S'\cup S''$, and
$(x_{1,1},\ldots,x_{1,n_1},\ldots,x_{i,1},\ldots,x_{i,n_i},\ldots,x_{n,n_n}$)
is a parametrization of $F$ as a normal surface in $\C$, then
$(x_1,\ldots,x_n)$, where $\sum_{j=1}^{j=n_i} x_{i,j} = x_i $, is
the parametrization of $F$ as a normal surface in $\T$. We can
think of the lattice point
$(x_{1,1},\ldots,x_{1,n_1},\ldots,x_{i,1},\ldots,x_{i,n_i},\ldots,x_{n,n_n}$)
as a point in the sub-cone of the solution cone
$\mathcal{S}(X,\mathcal{C})$, where we set all normal disk types
in $\C$ that meet $S'\cup S''$ equal to zero. We call the
parametrization
$(x_{1,1},\ldots,x_{1,n_1},\ldots,x_{i,1},\ldots,x_{i,n_i},\ldots,x_{n,n_n}$)
of $F$, a {\it re-writing} of the parametrization
$(x_1,\ldots,x_n)$ of $F$.

Finally, if $S$ is a two-sided normal surface, $M'$ is the
manifold obtained by splitting $M$ along $S$ and $S'$ and $S''$
the copies of $S$ in $\bdy M'$, then there is a natural
identification of $S'$ and $S''$ to recover the $3$--manifold $M$
with $S'$ and $S''$ being identified to $S$ in $M$. We will refer
to this as {\it re-attaching along $S'$ and $S''$}. In addition,
if $S^*$ is a normal surface in the cell-decomposition induced on
$M'$  ($M$ split along $S$) and possibly now $S^*$ meets $S'\cup
S''$, then when we re-attach along $S'$ and $S''$, we get a
subcomplex, denoted $S\cup S^*$, which is the image of $S'\cup
S''\cup S^*$ in $M$. We will call this subcomplex the {\it
piecewise linear normal surface} obtained from $S$ and $S^*$.

We list below some existence results for normal surfaces. We say a
$2$--sphere $S$ embedded in a $3$--manifold $M$ is {\it
inessential in $M$} if $S$  bounds a $3$--cell in $M$; otherwise
$S$ is {\it essential}.  A properly embedded disk $D$ in a
$3$--manifold $M$ is {\it inessential in $M$} if $\bdy D$  bounds
a disk $D'$ in $\bdy M$; otherwise, $D$ is {\it essential in $M$}.

\begin{thm}\label{normaldisk}\cite{haken1}
Let $M$ be a 3--manifold.  If there is an essential, properly
embedded disk in $M$, then for any cell-decomposition $\T$ of $M$
there is an essential, normal disk embedded in $M$.
\end{thm}

\begin{thm}\label{normalsphere}\cite{kne, sch}
Let $M$ be a 3--manifold.  If there is an essential, embedded $2$--sphere in $M$,
then for any cell-decomposition $\T$ of $M$ there is an essential, embedded,  normal $2$--sphere
in
$M$.
\end{thm}

\begin{thm}\label{kneser}$[Kneser's Finiteness Theorem]$\cite{kne} Suppose $\T$ is a
triangulation of the compact
$3$--manifold $M$. There is a nonnegative integer $N_0$ so that whenever $F_1,\ldots,F_n$
is a pairwise disjoint collection of normal surfaces in $M$ and $n\ge N_0$, then for some
$i\ne j, F_i = F_j$.
\end{thm}

The equality, $F_i = F_j$, can be interpreted as $F_i$ is normally isotopic to $F_j$ or,
equivalently, they have the same parameterization. A similar result is true for
cell-decompositions; however, Theorem \ref{kneser} is quite familiar and is the form in which we
use this result. If $M$ is closed, then one has, for example, $N_0\leq 5t$, where $t$ is
the number of tetrahedra in $\T$.

There are numerous other existence results for normal surfaces, not needed for this work. See, for
example,
\cite{haken1, jac-ree, jac-rub4, jac-tol, jac-oer}.

\begin{figure}[htbp]
            \psfrag{A}{(A) normal octagon}
            \psfrag{B}{(B) normal tube}
        \vspace{0 in}
        \begin{center}
        \epsfxsize = 3 in
        \epsfbox{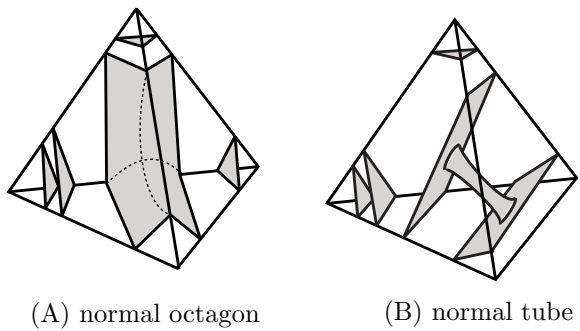}
        \caption{A normal octagon and a normal tube
between two normal quads.}
        \label{f-almostnormal}
        \end{center}
\end{figure}

In this work we also use almost normal surfaces in triangulations.
The elementary components of almost normal surface theory include
the normal triangles and normal quads of normal surface theory but
allow more general elementary components, normal octagons and
normal tubes. A {\it normal octagon} is a properly embedded disk
in a tetrahedron having boundary consisting of eight normal arcs
in the boundary of the tetrahedron; whereas, a {\it normal tube}
is a properly embedded annulus in a tetrahedron formed from two
disjoint normal triangles, two disjoint normal quads or a normal
triangle and a disjoint normal quad by joining them via a tube
parallel to an edge of the tetrahedron.  A normal octagon and a
normal tube are shown in Figure \ref{f-almostnormal}. There are
three types of normal octagons and twenty-five types of normal
tubes in each tetrahedron.

If $\T$ is a triangulation or ideal triangulation of the
3--manifold $M$, we say a surface $F$ is {\it almost normal} (with
respect to $\T$)
 if and only if the lift
of every component of $F$ in a tetrahedron is either a normal
triangle or a normal quadrilateral, except for at most one
tetrahedron, where we allow the lift of $F$ in the exceptional
tetrahedron to be precisely one of:
 \begin{itemize} \item  a collection of
normal triangles and one normal octagon, or \item  a collection of
normal triangles and normal quads with one normal tube.
\end{itemize}

In the case of a normal tube, we do not allow that the normal tube
is along an edge between two copies of the same normal surface.
So, an almost normal surface never contains both a normal octagon
and a normal tube but may contain one of them. An almost normal
surface with a normal octagon is called an {\it octagonal almost
normal surface} and one with a normal tube is called a {\it tubed
almost normal surface}. A compression of the tube gives a normal
surface and does not give two copies of the same normal surface.

Finally, we have from \cite{jac-ree, jac-tol} that if $M$ has an
essential $2$ sphere, then for any triangulation of $M$ there is a
vertex solution which is an essential, embedded, normal
$2$--sphere. Similarly, if $M$ has an essential, properly embedded
disk, then for any triangulation of $M$ there is a vertex solution
which is an essential, normal disk. Under the appropriate
hypothesis (which is one of the motivations of this work), one can
also conclude that if there are almost normal $2$--spheres, then
there is one that is a vertex-solution. We will need additional
results analogous to these for this work. See Propositions
\ref{0-decide}, \ref{anoctagonal} and \ref{0-b-decide}.

\section{Barriers and shrinking.}\label{sect-barriers}

In this section we collect a number of facts, which we need later
and which are, in their own right, very useful to many
applications.  The basic notions are that of a ``barrier surface"
and the operation of ``shrinking a surface" to a (possibly
disconnected) normal surface. Similar results may be found in
\cite{jac-rub1} and a complete study of barrier surfaces and
shrinking is in \cite{jac-let-rub1}. The ideas comes from
geometric analysis and minimal surface theory and were introduced
into normal surface theory in \cite{rub}. Given a surface which is
not normal, one wants to try to normalize the surface preserving
some control. The control comes from the barrier surface and the
shrinking is defined in terms of a sequence of permissable moves
and a complexity. Traditionally,  this complexity has been the
weight (area) of the surface; the generality in which we use this
notion requires a new complexity. Each permissable move reduces
the complexity; if there are no permissable moves, then the
surface is {\it stable}. It will follow that a stable surface has
components that are either normal surfaces or are $0$--weight
$2$--spheres or $0$-weight disks, each properly embedded in a cell
in our cell-decomposition.

Suppose $\T$ is a cell-decomposition (or ideal cell-decomposition) of the
$3$--manifold $M$ and $S$ is a properly
embedded surface in
$M$, which is transverse to the $2$--skeleton of $\T$, $\T^{(2)}$. The \emph{weight} of
$S$ is the cardinality of $S\cap \T^{(1)}$, $wt(S) = \abs{S\cap \T^{(1)}}$.  Now, recall that the
surface
$S$ is a normal surface if and only if the lift
 of every  component of the intersection of $S$ with a cell
is a normal disk (and in the case of a triangulation, a normal
triangle or a normal quad). So, for normal surface theory, we
would like the lifts of the components of $S$ in the various cells
to be disks, in general, and  normal disks, in particular. If the
lifts of all the components of $S$ in the cells are not disks, we
need a measure of this variance. We define the {\it local Euler
number of $S$}, written $\lambda_{\chi}(S)$, to be the sum
$\lambda_{\chi}(S) = \Sigma_{c\ne S^2} (1 - \chi(\tilde{c}))$,
where $c$ runs over all non $2$--sphere components of $S$ in the
cells of $\T$. Notice that $\lambda_{\chi}(S) = 0$ if and only if
each non spherical component of $S$ in a cell of the decomposition
$\T$ lifts to a disk. Finally, it is convenient to clean up some
of the $0$--weight intersections with the faces of $\T$. A
$0$--weight curve of intersection of $S$ with $\T^{(2)}$ is a
simple closed curve lying entirely in the interior of a face of
$\T$.
 Let $\sigma (S)$ denote the number of $0$--weight curves of the intersection of $S$ with
faces of $\T$, which are also in $\open{M}$, the interior of the
$3$--manifold $M$. We define the {\it complexity of $S$} to be
$\mathit{C}(S) = (wt(S),\sigma (S),\lambda_{\chi}(S))$, where we
consider the set of triples under lexicographical order from the
left.

\subsection{Shrinking surfaces: normalization}In this work we use four basic moves in shrinking
(normalizing)
 a properly embedded surface
$S$: a compression, an
 isotopy, a
$\bdy$--compression and finally a ``cleaning up" move, which is
not really necessary but brings a nice order to the notion of
stable surface. Of course, a $\bdy$--compression is relevant only
when $S$ (and hence, the manifold $M$) has nonempty boundary. We
show these moves in Figures
\ref{f-normalmoves1-new}--\ref{f-normalmoves2-new}.

We begin with a properly embedded surface $S$ meeting the
$2$--skeleton of the cell-decomposition transversely.  To keep
notation simple, we refer to the surface at each step of the
shrinking as $S$, understanding that it may have changed
considerably from the original surface $S$. The target of
shrinking is to arrive at a surface (a stable surface) having
components which  are normal surfaces or are properly embedded,
$0$--weight $2$--spheres and $0$--weight disks, each of which is
contained entirely in  some  cell of our cell-decomposition.
Hence, the lifts of the components of $S$
 in a cell will be  normal disks and properly embedded,
$0$--weight $2$--spheres and $0$--weight disks in the cell. Recall
that normal disks are characterized by their boundary curves in
the cells, which are made up of a finite number of normal spanning
arcs in the faces of the cells and which do not meet an edge in
the cell more than once. Finally, we use the terms ``compression"
and ``$\bdy$--compression" in a more general context than usual;
we do not require the boundary of a compressing disk to be an
essential curve in our surface nor do we require the arc in which
a $\bdy$--compressing disk meets our surface to be an essential
arc in the surface. Hence, a compression or a $\bdy$--compression
may split off only a trivial bit of the surface.

\vspace{.1 in}
The {\it normal moves} are:
\begin{enumerate}
\item {\it A compression in the interior of a cell. This move reduces the
local Euler number and does not change weight or the number of $0$--weight curves of intersection of $S$
with the faces of $\T$.} (See Figure
\ref{f-normalmoves1-new}.)\\

 A compression reducing the local Euler number can be made in the interior of some cell
whenever the local Euler number is not zero, $\lambda_{\chi}(S)\ne
0$. In this case, there is a component $c$ of the intersection of
$S$ with some cell, say $\Delta_i$ and for $\tilde{c}$ the lift of
$c$, we have $1 - \chi(\tilde{c}) > 0$; hence, there is a
compression of $\tilde{c}$ along an essential curve in $\td{c}$ in
$\td{\Delta}_i$. It follows, there is a disk $\td{D}$ embedded in
$\td{\Delta}_i$ so that $\td{D}\cap\td{c} = \bdy \td{D}$ and $\bdy
\td{D}$ is not trivial in $\td{c}$. Of course, it is possible that
$\td{D}$ meets other lifts of the components of $S$ in $\Delta_i$.
However, if this is the case, then we may assume the intersection
of the lifts of the components of $S$ in $\Delta_i$ meet $\td{D}$
in simple closed curves in the interior of $\td{D}$. Either we can
change our choice of $\td{D}$ to eliminate such intersections or
there is a lift $\td{c}'$ of a component $c'$ of $S$ meeting
$\Delta_i$ and a disk $\td{D}'\subset\td{D}$ so that $\bdy\td{D}'$
is an essential curve in $\td{c}'$ (in particular, $1 -
\chi(\tilde{c}') > 0$) and $\td{D}'$ does not meet any other lifts
of the components of $S$ in $\Delta_i$. Let's assume that $D$ is
such an innermost disk so we don't have to drag the prime notation
along. We let $D$ denote the image of $\td{D}$ in $\Delta_i$ and
compress $c$ along $D$ (which induces a compression of $\td{c}$
along $\td{D}$).  Notice such a compression
 does not affect the weight but since $\bdy\td{D}$ is essential in $\td{c}$, the
compression decreases the local Euler number; furthermore, this move does not affect the
intersection of $S$ with the interior of the faces or the edges of $\T$. Hence, it reduces the complexity
of the surface
$S$.

Note that this move may be an essential compression of the surface
$S$ and thereby a change of its topological type. If $S$ were
incompressible, we could argue that we still have a component,
after this move, homeomorphic with $S$; and if $M$ were also
irreducible, we could argue that we have a surface isotopic with
$S$. In our case we do not care. We wait until we arrive at a
stable surface and then make an analysis.

\begin{figure}[htbp]
            \psfrag{c}{compression}
 \vspace{0 in}
        \begin{center}
        \epsfxsize = 3 in
        \epsfbox{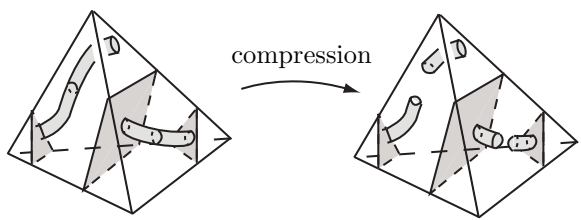}
        \caption{A compression in the  interior of a cell. This move
reduces the local Euler number and does not change weight or the
$0$-weight curves in the interior of faces of $\T$.}
        \label{f-normalmoves1-new}
        \end{center}
\end{figure}

\vspace{.25 in}

\item {\it An isotopy reducing the number of times the boundary of a lift of a
component meets an edge of a cell, where the edge is in the interior of $M$. This move reduces
$wt(S)$.} (See  Figure \ref{f-normalmoves4-new}.)\\

 At this stage we may assume all the lifts of the components of $S$ in a cell of $\T$ are either
properly embedded  disks or  $2$--spheres. Now, suppose
$\tilde{c}$ is a lift of a component of $S$ in a cell $\Delta_i$;
$\td{c}$ is a properly embedded disk; and $\tilde{c}$ meets an
edge $e$ of $\td{\Delta}_i$, the lift of $\Delta_i$, more than
once. Then if we consider the curve $\bdy\td{c}$, it divides the
edge $e$ into a number of subarcs and there is at least one, say
$\td{\beta}$ which has both its end points in $\bdy\td{c}$. Hence,
there is a disk $\td{D}$ embedded in $\td{\Delta}_i$ so that
$\td{D}\cap\td{c}$ is an arc $\td{\alpha}\subset\bdy \td{D}$,
$\td{\alpha}\cup\td{\beta} = \bdy \td{D}$ and
$\td{\alpha}\cap\td{\beta} = \bdy\td{\alpha} =\bdy\td{\beta}$. See
Figure \ref{f-normalmoves4-new}. However, it is possible that
$\td{D}$ meets other lifts in $\td{\Delta}_i$ of the components of
$S$ in $\Delta_i$. If this is the case, then $\td{D}$ meets lifts
other than $\td{c}$ in simple closed curves in the interior of
$\td{D}$ and spanning arcs in $\td{D}$ having both their end
points in $\td{\beta}$. Standard techniques allow us to choose
$\td{D}$ so that there are no such simple closed curve components.
Hence, if there are spanning arcs remaining, we can choose one
$\td{\alpha}'$, which is ``outermost" in the sense that   there is
a disk $\td{D}'\subset \td{D}$ and a lift $\td{c}'$ of a component
$c'$ of $S$ in $\Delta_i$ so that $\td{D}'\cap\td{c}' =
\td{\alpha}'$, $\bdy \td{D}' = \td{\alpha}'\cup\td{\beta}'$, where
$\td{\beta}'\subset\td{\beta}$ and $\td{\alpha}'\cap\td{\beta}' =
\bdy\td{\alpha}' = \bdy\td{\beta}'$. Furthermore, $\td{D}'$ does
not meet any other lifts of components of $S$ in $\Delta_i$. As
above, having demonstrated that we can find such an outermost
disk, we assume the original disk $\td{D}$ has this property so we
do not need to drag along the prime notation.

We consider the image $D$ of $\td{D}$ in $\Delta_i$. Then $D$ is
an embedded disk in ${\Delta_i}$; $D$ only meets $S$ in $c$;
$D\cap c = \alpha\subset\bdy D$ is a spanning arc of $c$; $D$
meets the boundary of $\Delta_i$ in the arc $\beta$ in the edge
$e$ of $\Delta_i$ ($e$ is also used for the image of the edge $e$
in $\td{\Delta}_i$); and $\beta\subset \bdy D$, where $\alpha\cap
\beta = \bdy\alpha = \bdy\beta$ and $\alpha\cup\beta = \bdy D$.
See Figure \ref{f-normalmoves4-new}. Since we have assumed the
edge $e$ containing $\beta$ is in the interior of the
$3$--manifold, then there is an isotopy of $S$, splitting $c$ into
two disks and reducing  $wt(S)$. Of course, this move may increase
the value $\sigma (S)$ and the local Euler number; however, it
reduces the complexity of $S$.

\begin{figure}[htbp]
            \psfrag{a}{$\alpha$}
            \psfrag{b}{$\beta$}
            \psfrag{D}{D}
            \psfrag{e}{e}
            \psfrag{c}{\begin{tabular}{c}
            \small{$\bdy$--compression}\\
            \small{$e\subset\bdy M$} or\\
            \end{tabular}}
            \psfrag{i}{\begin{tabular}{c}
            \small{isotopy}\\
           \small{$e\subset\open{M}$}\\
            \end{tabular}}
        \vspace{0 in}
        \begin{center}
        \epsfxsize=4 in
        \epsfbox{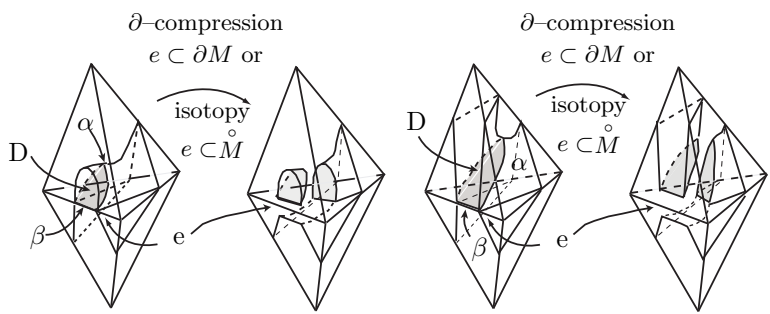}
        \caption{An isotopy or $\bdy$--compression reducing the  number of
times a lift of a component meets an edge of a cell. This move
reduces $wt(S)$.}
        \label{f-normalmoves4-new}
        \end{center}

\end{figure}

\vspace{.25 in}

\item {\it A $\bdy$--compression reducing the number of times a lift of a component meets an
edge of a cell, where the edge is in the boundary of $M$. This moves  reduces  $wt(S)$.}
(See Figure
\ref{f-normalmoves4-new}.)\\

Instead of, as above, where the edge $e$  is in the interior of
$M$, we now have the edge $e$ in $\bdy M$. In this case, the move
must be accomplished by a $\bdy$--compression rather than an
isotopy.  This move reduces $wt(S)$; however,  as above in the
isotopy move, it may increase the value $\sigma (S)$ and the local
Euler number. In any case, it reduces the complexity of $S$.
Again, see Figure \ref{f-normalmoves4-new}.

Note that as above when we had a compression, a $\bdy$--compression can change the topological
type of $S$. However, if $S$ were $\bdy$--incompressible, then we would still have a component
topological equivalent to the one before the $\bdy$--compression and if $M$ were
$\bdy$--irreducible and irreducible, we would have a component isotopic to the one before the
$\bdy$--compression. But just as above, in our case, this does not matter and we wait until we
have a stable situation to make an analysis.

\vspace{.125 in}

\item {\it A compression eliminating $0$--weight simple closed curve components
from the intersection of
$S$ with the faces of $\T$ in the interior of $M$.} (See Figure
\ref{f-normalmoves2-new}.)\\

 We may assume there are no essential compressions in the interior of any
cell; actually, at this stage, we may assume the lift of any
component of $S$ in a cell is either a normal disk or a properly
embedded, $0$--weight $2$--sphere or $0$--weight disk in the cell.
If there is a $0$--weight simple closed curve common to $S$ and a
face $\sigma$ of a cell, say $\Delta_i$,  then there is an
innermost one. Such a $0$--weight curve bounds a  component $c$
 of
$S$ in $\Delta_i$, which is a $0$--weight disk,
 properly embedded in $\Delta_i$, and having its boundary entirely in the interior of the face
$\sigma$. Furthermore, its boundary bounds a disk $D$ in the
interior of $\sigma$. If $\sigma$ is in the interior of the
manifold, then there is a similarly embedded disk $c'$ on the
other side of this face, $\bdy D = \bdy c'$ and $c\cup c'$ is a
small $2$--sphere, which can be isotoped entirely into the
interior of one of the cells or we can perform a compression along
the
 disk $D$ in $\sigma$ and create two $0$--weight $2$--spheres, one in each cell
sharing the face, $\sigma$. We choose to do the latter and
therefore get two $2$--spheres, each embedded entirely in the
interior of a cell; this is one of our stable situations. Thus we
eliminate all $0$--weight simple closed curves in the interior of
faces of $\T$, which are also in the interior of $M$.
 These moves reduce the value
$\sigma (S)$ and do not affect the weight or the local Euler number.

Note that we could just throw away all of the $0$--weight pieces,
which in practice is essentially what we do; but what we have done
here reduces our work later when we need to analyze what we have
after we shrink a surface.

\vspace{.25 in}

\begin{figure}[htbp]
            \psfrag{a}{after}
            \psfrag{b}{before}
            \psfrag{c}{$c$}
            \psfrag{d}{$c'$}
            \psfrag{t}{$\sigma$}
            \psfrag{o}{compression}
            \psfrag{S}{$S^2$}
            \psfrag{D}{$\mathbb{D}^2$}
            \psfrag{f}{\begin{tabular}{c}
           {faces}\\
           {identified}\\
            \end{tabular}}
            \psfrag{M}{$\sigma\subset\bdy M$}
            \psfrag{n}{$\sigma\subset\open{M}$}
        \vspace{0 in}
        \begin{center}
        \epsfxsize=4.25 in
        \epsfbox{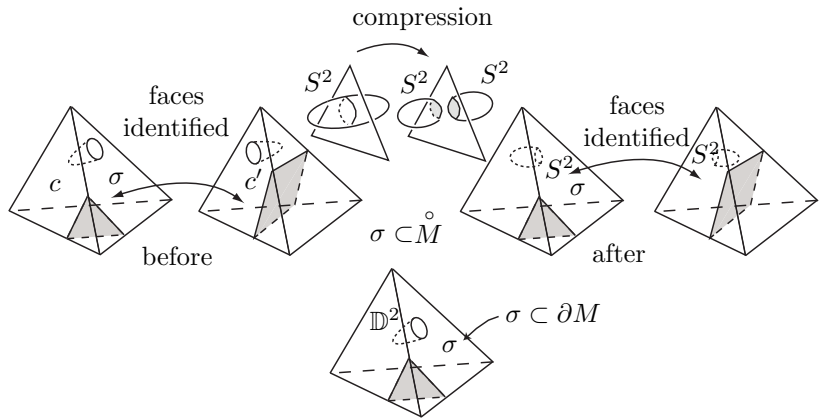}
        \caption{A
compression removing a  simple closed curve component from the
intersection of $S$ with a face in the interior of the manifold.}
        \label{f-normalmoves2-new}
        \end{center}

\end{figure}

\end{enumerate}

As mentioned above, compressions and $\bdy$--compressions can
alter the properties of the surface $S$. We can place conditions
on the surface and the $3$--manifold so we can make these
modifications and still maintain certain properties of the
surface; for example, as indicated above, the surface is
incompressible and $\bdy$--incompressible and the manifold is
irreducible and $\bdy$--irreducible. However, these normal moves
can be made on any surface. The moves never increase weight. We
call a sequence of these normal moves a \emph{shrinking} of the
surface $S$ and we allow the surface $S$ to be compressed and
$\bdy$--compressed, possibly resulting in a number of distinct
components. After a finite number of steps, the components will
either be normal or will be properly embedded, $0$--weight
$2$--spheres and $0$--weight disks contained entirely in various
cells of our cell-decomposition. We say $S$ has been \emph
{shrunk} to the resulting surfaces or sometimes we throw away the
$0$--weight components and say $S$ has been shrunk to the
resulting normal surface(s). Of course, it may be that $S$ has
shrunk until it disappears (has only $0$--weight components).

\subsection{Barrier surfaces} Suppose $B$ be a properly embedded surface
in a $3$--manifold
$M$ and let
$N$ be a component of the complement of $B$, $M\setminus B$. We
say $B$ is a \emph{barrier surface for $N$}, or simply a \emph{barrier}, if
any properly embedded, compact
surface $F$ in $N$ can be shrunk in $N$.

\begin{figure}[htbp]
            \psfrag{A}{(A)}
            \psfrag{B}{(B)}
            \psfrag{C}{(C)}
            \psfrag{b}{$\td{b}$}
            \psfrag{c}{$\td{C}$}
            \psfrag{d}{$\td{D}$}
            \psfrag{s}{$\td{\alpha}$}
            \psfrag{t}{$\td{\beta}$}
            \psfrag{r}{$\sigma$}
            \psfrag{e}{$\td{E}$}
            \psfrag{f}{$\td{E}'$}
            \psfrag{x}{\begin{tabular}{c}
           $\bdy$--compression\\
           disk $\td{D}$\\
            \end{tabular}}
        \vspace{0 in}
        \begin{center}
        \epsfxsize=3.5 in
        \epsfbox{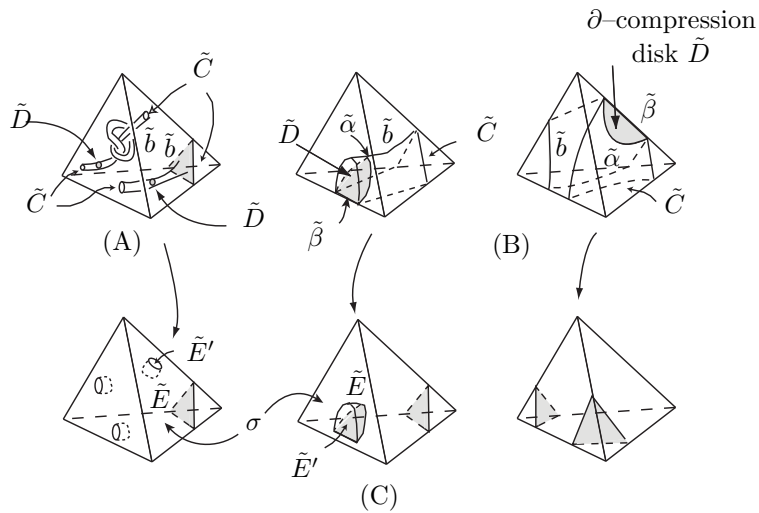}
        \caption{Complete system of  compressing disks for a barrier
surface.}
        \label{f-compressionbarrier-new}
        \end{center}

\end{figure}

We give criteria
for a  surface
$B$, properly embedded in $M$, to be a barrier surface for a
component $N$ of its complement  in $M$.
Suppose $\Delta$ is a cell of $\T$ and $C$ is the closure of
a component of $\Delta\cap
N$, $\td{\Delta}$ the lift of $\Delta$ and  $\td{C}$ the lift of $C$. Let $b = C\cap B$ and
$\td{b}$ denote the lift of $b$. A collection of pairwise disjoint disks in
$\td{C}$ is said to be a \emph{complete system of compressing disks for $B$ in $C$} if
\begin{itemize}
\item a disk, $\td{D}$, in the collection meets $\td{b}$ only in
its boundary, which is an essential curve in $\td{b}$; i.e.,
$\td{D}$ is an essential compressing disk for $\td{b}$ in $\td{C}$
(see Figure \ref{f-compressionbarrier-new}(A)),

\item a disk, $\td{D}$, in the collection meets $\td{b}$ in a
properly embedded arc $\td\alpha$ and meets boundary of
$\td{\Delta}$ in an arc $\td\beta$, which is entirely in the
interior of an edge of $\td{\Delta}$, $\td\alpha\cup\td\beta =
\bdy \td{D}$ and $\td\alpha\cap\td\beta = \bdy\td\alpha =
\bdy\td\beta$; i.e., $\td{D}$ is a $\bdy$--compression (not
necessarily essential) of $\td{b}$ in $\td{C}$ (see Figure
\ref{f-compressionbarrier-new}(B)), and \item each component
remaining after $\td{b}$ has been compressed and
$\bdy$--compressed along the collection of disks is either a
normal disk for $\T$ or a properly embedded $0$--weight disk
$\td{E}$ having its boundary entirely in the interior of a face
$\sigma$ of $\td{\Delta}$, $\bdy\td{E}$ bounds a disk
$\td{E}'\subset \sigma$ and the $2$--sphere $\td{E}\cup\td{E}'$
bounds a $3$--cell in $\td{C}$ (see Figure
\ref{f-compressionbarrier-new}(C)).
\end{itemize}

\begin{lem} Suppose $\T$ is a cell-decomposition (or ideal cell-decomposition) of
the $3$--manifold $M$ and $B$ is a properly embedded surface in
$M$. The surface $B$ is a barrier surface for the component $N$ of
$M\setminus B$ if for each cell $\Delta$ of $\T$ and the closure
of each component $C$ of $\Delta\cap N$, there is a complete
system of compressing disks for $B$ in $C$.
\end{lem}

\begin{proof} Suppose $F$ is a properly embedded surface in $N$. The surface
$F$ misses $B$ ($N$ is a component of $M\setminus B$) and we may assume that $F$ is transverse to
$\T^{(2)}$, the $2$--skeleton of $\T$. We choose a properly embedded surface
$S$ that has  minimal complexity among  all surfaces that can be obtained by shrinking
$F$ missing $B$; i.e., sequences of compressions and $\bdy$--compressions,
missing
$B$, and
 isotopies, which are identity on $B$. We claim each component of $S$ is
either normal or is a $0$--weight $2$--sphere or disk properly
embedded in a cell of $\T$.

To see this, suppose $\Delta$ is a cell of $\T$. Let $C$ be a
component of $N\cap\Delta$. Let $\td{C}$ be the lift of $C$; let
$b = B\cap C$ and let $\td{b}$ be the lift of $b$. By hypothesis,
there is a complete system of disks for $B$ in $C$. Suppose $c =
S\cap C$ and $\td{c}$ is the lift of $c$. If $\td{c}$ meets the
complete system of compressing disks for $B$ in $C$, it does so in
a very nice way. Specifically, if $\td{D}$ is a compressing disk
in $\td{C}$ for $\td{b}$, any component of $\td{c}\cap \td{D}$ is
a simple closed curve in the interior of $\td{D}$. If $\td{D}$ is
a $\bdy$--compressing disk in $\td{C}$ for $\td{b}$, any component
of $\td{c}\cap \td{D}$ is either a simple closed curve in the
interior of $\td{D}$ or a spanning arc; furthermore, the spanning
arcs have both end points in the same edge of $\td{\Delta}$. It
follows that none of these  intersections of $\td{c}$ with the
complete system of compressing disks are essential in $\td{c}$ or
we could make a sequence of moves, as defined above, to reduce the
complexity of $S$. Hence, we can assume that $\td{c}$ misses a
complete system of compressing disks for $B$ in $C$. This is true
in every cell of $\T$. Notice {\it a priori} it may seem we have
to go back to a cell in which we have made these moves before. But
this could only happen if we reduce the weight of $S$, which would
contradict our choice of $S$.

Now, for any cell $\Delta$ and any component $C$ of $N\cap\Delta$, if $c$ is
a component of $S$ meeting $C$ and $\td{c}$ is the lift of $c$, we can make all
compressions and
$\bdy$--compressions along a complete system of disks for $B$ in $C$, missing $\td{c}$.  Hence,
$\td{c}$ lies in a component of $\td{\Delta}$ determined by normal disks and properly embedded
 $0$--weight disks parallel into the interior of a face of $\td\Delta$ through a $3$--cell in
$\td C$.  It follows, there is no obstruction to making further moves on $\td{c}$ to reduce
complexity unless each component of $\td{c}$ is a
 normal disk or  a properly embedded, $0$--weight $2$--sphere or $0$--weight disk
 lying entirely in a cell of $\T$.
\end{proof}

We now list several examples of barrier surfaces. Additional
examples are given in \cite{jac-rub1} and \cite{jac-let-rub1}. As
above, if $S$ is a two-sided, normal surface in $M$ and $M'$ is
the manifold obtained from $M$ by splitting along $S$, then we let
$S'$ and $S''$ denote the copies of $S$ in $\bdy M'$. Furthermore,
if $S^*$ is a normal surface in the induced cell structure on
$M'$, we let $S\cup S^*$ denote the piecewise normal surface
obtained from $S$ and $S^*$. Similarly, if $\mathcal{K}$ is a
subcomplex of the cell-decomposition induced on $M'$ by $\T$, then
when we re-attach along $S'$ and $S''$, we get a subpolyhedron,
denoted $S\cup \abs{\mathcal{K}}$, which is the image of $S'\cup
S''\cup \abs{\mathcal{K}}$ in $M$.

\begin{thm}\label{barrier} Suppose $\T$ is a cell-decomposition (or ideal
cell-decomposition) of the $3$--manifold $M$.
\begin{enumerate}
\item If $S$ is a normal surface or an almost normal surface in
$M$ and $B$ is the boundary of a small regular neighborhood of $S$
in $M$, then $B$ is a barrier surface for each component of its
complement not meeting $S$. Often in this case, we just say the
normal surface or almost normal surface $S$ is a barrier surface
for each component of its complement.

\item If $S$ is a two-sided, normal surface in $M$ and $S^*$ is a normal surface in the
induced cell structure on
$M$ split along $S$ and $B$ is
 the boundary of a small regular neighborhood (in $M$) of the piecewise normal surface $S\cup S^*$,
then $B$ is a barrier surface
for each component of its complement  not meeting $S\cup S^*$.

\item If $X$ is a finite union of normal surfaces in $M$, which meet
transversely, and $B$ is the boundary of a small regular neighborhood of $X$, then $B$
is a barrier surface for each
 component of its complement  not meeting $X$.

\item If $\mathcal{K}$ is a subcomplex of $\T$ and $B$ is  the boundary of a small regular
neighborhood of the underlying point set of $\mathcal{K}$, $\abs{\mathcal{K}}$, then $B$ is a
barrier surface for any component of its complement  not meeting $\abs{\mathcal{K}}$.

\item If $S$ is a normal surface in $M$ and $\mathcal{K}$ is a subcomplex of the
cell structure induced by $\T$ on $M$ split along $S$ and $B$ is the boundary  of a
small regular neighborhood of $S\cup \abs{\mathcal{K}}$, then $B$ is a barrier surface
for any component of its complement  not meeting $S\cup\abs{\mathcal{K}}$.

\item If $S$ is a two-sided, normal surface in $M$ and $\mathcal{K}$ is a
subcomplex of the cell structure induced by $\T$ on $M'$, the manifold obtained by  spliting $M$ along
$S$, and
$F$ is  the frontier of a small regular neighborhood of  $\abs{\mathcal{K}}$, then $F$ is a barrier
surface in $M'$ for any component of its complement  not meeting $\abs{\mathcal{K}}$.
\end{enumerate}
\end{thm}

\begin{proof}
The proof is straight forward. In situations 3, 4, 5 and 6 there are numerous
cases to consider; but listing the cases is the only task. \end{proof}

Barrier surfaces are used as a tool throughout this work  and appear in the middle of certain
arguments; however, we give here some general notions, which exhibit the use of barrier surfaces.

Suppose $B_1,\ldots,B_n$ is a pairwise disjoint collection of $3$--cells in
$S^3$. We call the $3$--manifold $M = S^3 \setminus\bigcup_{i = 1}^{n}\stackrel{\circ}{B}_i$ a
{\it punctured $3$--sphere}. In particular, the collection may be empty; so, we allow that the
$3$--sphere, itself, is a punctured $3$--sphere, of course, without any punctures. In the case we
know the boundary is not empty, we may also say we have a {\it punctured $3$--cell}.

In our
definition of a barrier surface $B$, we have that whenever $B$ is a barrier surface for a component of its
complement, then in all cases $B$ is two-sided. So, if $B$ is a
barrier surface for the component
$N$ of its complement, then by taking a small regular neighborhood of $B$, we have a copy of
$B$ in
$N$. We can then shrink this copy. We use the phrase ``shrink $B$ in $N$" in this
situation. Also, if the barrier surface $B$ is a normal surface or has components which are normal
surfaces, then we may say shrink $B$, understanding that each component of $B$ which is normal is
stable and there are no normalization moves on such components. Such components survive (are
never touched) in the shrinking.

Our first result uses a barrier surface to engulf the vertices of a triangulation.

\begin{prop}\label{engulf-sphere} Suppose $\T$ is a triangulation of the  closed,
orientable, irreducible
$3$--manifold $M$ and suppose $S$ is a normal $2$-sphere in $M$ which bounds a
$3$--cell $E$ in
$M$. Then there is a  normal $2$--sphere $S'$ bounding a
$3$--cell $E'$ in $M$, which contains $E$ and all the vertices of $\T$ or $M = S^3$.
\end{prop}
\begin{proof} If all the vertices of $\T$ are in $E$, then there is nothing to prove. Otherwise,
split $M$ along $S$ and let $M'$ denote the component not meeting $E$. We will continue to use
$S$ to denote the copy of $S$
in the boundary of
$M'$. There is a subcomplex
$\Lambda$ of the $1$--skeleton of the induced cell structure on $M'$ so that each component of
$\Lambda$ is a tree and meets $S$ in precisely one point and
$\Lambda$ contains all vertices of
$\T$ not in
$E$. By the above theorem, the frontier $B$ of a small regular neighborhood $N = N(S\cup\Lambda)$
of
$S\cup\Lambda$ is a barrier surface in the component of its complement in  $M'$ not meeting
$S\cup\Lambda$. Note that $N$
is a punctured $3$--cell  in
$M'$ (actually, $B$ is isotopic to $S$ and in this
case $N = S^2\times I$), which contains  all the vertices of $\T$ not in $E$.

If $B$ is normal, then it is itself stable and does not shrink. In
this case, we let $E' = E\cup N$ and let $S' = B$. So, we may
assume $B$ is not normal. We can shrink $B$ in $M'$; furthermore,
the point of $B$ being a barrier is that this shrinking will not
meet $S\cup\Lambda$.   A shrinking (a finite sequence of normal
moves) involves either a compression or an isotopy move ($M$ is
closed). So, assume we are at a stage in our shrinking where we
have a finite number of pairwise disjoint $2$--spheres
$S_1,\ldots,S_n$, and a punctured $3$--cell $P_k$ with
$S\cup\Lambda\subset P_k$ and $\bdy P_k$ includes $S$ and the
spheres  $S_1,\ldots,S_n$. If the collection $S_1,\ldots,S_n$ is
not stable, then there is either an isotopy normal move or a
compression normal move on one of these $2$--spheres.

An isotopy move is across an edge in $N$ missing $S\cup\Lambda$;
hence, we have an isotopy move of some $2$--sphere, say $S_i$, in
$\bdy P_k$. We get a new collection of $2$--spheres
$S_1,\ldots,\S_i',\ldots,S_n$ where $S_i'$ replaces $S_i$ and
$wt(S_i') < wt(S_i)$. We let $P_{k+1}$ denote the image of $P_k$
under this isotopy. Then $P_{k+1}$ is a punctured $3$--cell
containing $S\cup\Lambda$.

If there is a compression on one of the $2$--spheres, say $S_i$,
then let $D$ denote the compressing disk. Not only does $D$ not
meet any $2$--sphere in the collection except for $S_i$, which it
meets in its boundary, $D$ does not meet $S\cup\Lambda$. If
$D\subset P_k$, then $D$ splits $P_k$ into two punctured
$3$--cells, one, say $P_{k+1}$ containing $S\cup\Lambda$, and
$\bdy D$ splits $S_i$ into two $2$--spheres, $S_i'$ and $S_i''$,
with, say $S_i'\subset P_{k+1}$. We have a new collection of
$2$--spheres, $S_{j_1}\ldots,S_i',\ldots, S_{j_m}$, which along
with $S$ make up the boundary of our new punctured $3$--cell
$P_{k+1}$. If $D$ is not in  $P_k$, then again $\bdy D$ splits
$S_i$  into two $2$--spheres, $S_i'$ and $S_i''$; however, in this
case, a compression is adding a $2$--handle to $S_i$ and we get a
new punctured $3$--cell, $P_{k+1}$, containing $P_k$ and having
both $S_i'$ and $S_i''$ in its boundary. Also,
$S\cup\Lambda\subset P_{k+1}$.

It follows that in shrinking $B$ and in the stable situation we have a punctured $3$--cell $P$,
$S\cup\Lambda\subset P$ and the boundary of $P$ consists of $S$ along with possibly
some other normal $2$--spheres and possibly some $0$--weight $2$--spheres entirely in the
interior of cells in the induced cell structure on $M'$. Each $0$--weight $2$--sphere bounds a
$3$--cell whose interior misses $P$. We fill in these $2$--spheres with these
$3$--cells and continue to call our punctured $3$--cell $P$. Now, since
$M$ is irreducible, each normal
$2$--sphere in the boundary of
$P$ bounds a $3$--cell in $M$. If such a boundary component, other than $S$, bounds a
$3$--cell whose interior misses $P$ we add that $3$--cell to $P$. We will continue to call the
punctured
$3$--cell $P$. So, we now have that $S\cup\Lambda\subset P$ and any $2$--sphere in boundary
of $P$ other than $S$ does not bound a $3$--cell whose interior misses $P$. If $S$ is
the only boundary component of $P$, then $M$ is $S^3$. If $S'$ is a component of the
boundary of $P$ distinct from $S$, then by $M$ irreducible, $S'$ bounds a $3$--cell,
say $E'$, in $M$. But then we have $E\cup P\subset E'$. So, such an $S'$ and $E'$ satisfy
the conclusions of our proposition.\end{proof}

There is a useful variation to the previous proposition when we do
not assume the $3$--manifold $M$ is irreducible; namely, we have
either there is a collection $S_1\ldots, S_n$ of normal
$2$--spheres bounding a punctured $3$--cell $P$ in $M$, where $P$
contains $E$ and all the vertices of $\T$, or $M$ is $S^3$. We
also can use Proposition \ref{engulf-sphere} by, say, choosing $S$
to be a vertex-linking normal $2$--sphere and $E$ the $3$--cell it
bounds to conclude that for any triangulation $\T$ of a closed,
orientable, irreducible $3$--manifold $M$, there is a normal
$2$--sphere bounding a $3$--cell containing all the vertices of
$\T$ or it follows that $M$ is $S^3$. There are triangulations of
$S^3$ for which there is no normal $2$--sphere bounding a
$3$--cell containing all the vertices of the triangulation. For
example, the triangulation given in Figure \ref{f-tetra}(4) is
such a triangulation. There also are more interesting ones.

Suppose $F$ is a closed, orientable surface. Let $F\times [0,1]$ be the product of $F$
with the unit interval and let $\gamma_1,\ldots,\gamma_n$  be a finite, pairwise disjoint
collection of simple closed curves in
$F\times 0$; it is not necessary that the $\gamma_i$ be essential. Choose small
regular neighborhoods
$N(\gamma_1),\ldots,N(\gamma_n)$ of the
$\gamma_i, 1\leq i\leq n,$ in $F\times 0$ so that $N(\gamma_i)\cap N(\gamma_j) = \emptyset, i\neq j$. Let
$D_1\times[0,1],\ldots,D_n\times [0,1]$ be a collection of $2$--handles, where $D_i, 1\leq i\leq
n$, is a $2$--cell. A $3$--manifold is obtained by attaching the $2$-handles,
$D_i\times [0,1]$ along the
$\gamma_i$; i.e., identifying  the annulus $\bdy D_i\times [0,1]$ with the annulus $N(\gamma_i)$
 for $1\leq i\leq n$. $F\times 1$ is a component of the boundary of this $3$--manifold. There may
be some number of $2$--sphere components in the boundary as well.
We may or may not fill in some of the $2$--sphere boundary
components with $3$--cells ($3$--handles). We call the resulting
$3$--manifold, say $H$, a {\it compression body} and denote the
boundary component $F\times 1$ by $\bdy_+H$ and denote the
remaining boundary, which may not be connected, by $\bdy_{-}H$. A
component of $\bdy_- H$, which is not a $2$--sphere, is
incompressible in $H$. If $\bdy_-H = \emptyset$, then $H$ is a
handlebody and if each component of $\bdy_-H$ is a $2$--sphere,
then $H$ is a punctured handlebody. Finally, $F\times [0,1]$ is
itself a compression body, as is a punctured $F\times [0,1]$.

We have the following result which is analogous to Proposition \ref{engulf-sphere}.

\begin{prop}\label{engulf-surface} Suppose $\T$ is a triangulation  of the closed,
 orientable
$3$--manifold $M$ and suppose $F$ is a normal, two-sided surface, embedded in $M$. Then there are
compression bodies $H'$ and $H''$ embedded in
$M$ so that $H'\cap H'' = F = \bdy_+H' = \bdy_+H''$, each component of $\bdy_-H'$ and $\bdy_-H''$
is normal    and
$H'\cup H''$ contains all vertices of
$\T$.
\end{prop}
\begin{proof} The proof is very similar to the proof of Proposition \ref{engulf-sphere}. Let $M'$
denote the manifold we get by splitting $M$ along $F$; let $F'$ and $F''$ denote the copies of
$F$ in $\bdy M'$.  There are  disjoint subcomplexes
$\Lambda'$ and $\Lambda''$ of the $1$--skeleton of the induced cell structure on $M'$ so that each
component of
$\Lambda'$ and $\Lambda''$ is a tree and meets  $F'$ and $F''$, respectively, in precisely one
point,
$\Lambda'\cap F'' = \emptyset = \Lambda''\cap F'$ and
$\Lambda'\cup\Lambda''$ contains all vertices of
$\T$. Let $B'$ and $B''$ be the boundaries of  small regular neighborhoods of $F'\cup\Lambda'$
 and $F''\cup\Lambda''$, respectively. Then $B'\cup B''$ is a barrier surface for the components
of their complements not meeting $F'\cup\Lambda'\cup
F''\cup\Lambda''$; furthermore, $B'$ and $F'$ are the boundaries
of a compression body as well as $B''$ and $F''$; actually, in
these cases the compression bodies are products. We shrink $B'\cup
B''$. In shrinking $B'\cup B''$ we obtain two compression bodies
$G'$ and $G''$ so that $\bdy_+G' = F'$, $F'\cup \Lambda'\subset
G'$, $\bdy_+G'' = F''$, $F''\cup \Lambda''\subset G''$ and each
component of $\bdy_-G'$ and of $\bdy_-G''$ is either a normal
surface or a $0$--weight $2$--sphere contained entirely in the
interior of a cell in the induced cell structure on $M'$. Any such
$0$--weight $2$--sphere bounds a $3$--cell missing
$F'\cup\Lambda'$ and $F''\cup\Lambda''$. We fill in these
$0$--weight $2$--spheres with such $3$--cells.

Now, when we reattach $F'$ and $F''$ to get $M$ and set $H'$ equal to the image of $G'$ and set
$H''$ equal to the image of $G''$, we get the desired compression bodies.\end{proof}

The following is essentially a direct generalization of Proposition \ref{engulf-sphere} for
engulfing the vertices of a triangulation by a handlebody; it uses the previous
proposition and has numerous useful variants.

\begin{prop}\label{engulf-handlebody} Suppose $\T$ is a triangulation  of the closed
 orientable, irreducible
$3$--manifold $M$ and suppose $F$ is an embedded normal surface,
which bounds a handlebody $H$ in $M$. If $F$ is incompressible in
$M\setminus \stackrel{\circ}{H}$ and is not contained in a
$3$--cell in $M$, then there is a normal surface $F'$ embedded in
$M$, $F'$ is parallel to $F$ and bounds a handlebody $H'$ in $M$
so that $H\subset H'$ and $H'$ contains all the vertices of
$\T$.\end{prop}

\begin{proof} We can repeat the argument used in the proof of Proposition \ref{engulf-surface}.
In this case we split $M$ along $F$ and consider only the
component that does not meet the handlebody $H$, call it $M'$. We
denote the copy of $F$ in $\bdy M'$ by $F'$. We have the
subcomplex $\Lambda'$ as above and we let $B'$ denote the boundary
of a small regular neighborhood of $F'\cup\Lambda'$. $B'$ is a
barrier surface in the component of its complement not meeting
$F'\cup\Lambda'$. Furthermore, $F'$ and $B'$ bound a compression
body which is homeomorphic to $F'\times I$. We shrink $B'$.
However, since $F'$ is incompressible in $M'$ (hence, $B'$ is
incompressible in $M'$), each normal move which is a compression
is an inessential compression. It follows that in the stable
situation we have a surface which is a copy of $F'$ and every
other component is either a normal $2$--sphere or a $0$--weight
$2$--sphere contained entirely in the interior of a cell in the
induced cell decomposition of $M'$. Since $M$ is irreducible and
$F$ is not contained in a $3$--cell, we can fill in each
$2$--sphere boundary component with a $3$--cell whose interior
does not meet the compression body.
 It follows that we have a compression body
$G'$ with $\bdy_+G' = F'$, $\bdy_-G'$ a normal surface isotopic to $F'$ and $F'\cup\Lambda'\subset
G'$. When we reattach $M'$ to $H$ to get $M$, the image of $G'$ along with $H$ gives us the
desired handlebody $H'$.\end{proof}

The next application is referred to as a ``double barrier" argument.

\begin{prop}\label{doublebarrier} Suppose $\T$ is a triangulation  of the compact,
 orientable
$3$--manifold $M$ and suppose $K$ and $L$ are disjoint subcomplexes in $\T$. Then there is a
 normal surface $F$ in $M$ separating $K$ and $L$.\end{prop}

\begin{proof} Let $B_K$ and $B_L$ denote the frontiers of small
regular neighborhoods of $K$ and $L$, respectively, chosen so that
$B_K\cap B_L = \emptyset$. Then $B_K$  and $B_L$ are barrier
surfaces in the component of the complement of $B_K\cup B_L$ not
meeting $K\cup L$. Furthermore, $B_K$ separates $K$ and $L$. We
shrink $B_K$. In shrinking $B_K$, we have compressions,
$\bdy$--compressions and isotopy moves. An isotopy move or
$\bdy$--compression occurs through the interior of an edge that
meets $B_K$ and so is away from $K$ or $L$. A compression is
entirely in the interior of a cell or the face of a cell and so
does not run through $K$ or $L$. So, in our stable situation we
have components which are either normal surfaces or $0$--weight
$2$--spheres and disks which are properly embedded in the cells of
our induced cell decomposition; furthermore, the union of these
components separate $K$ from $L$. We wish to eliminate the
$0$--weight $2$--spheres and disks. None of the $0$--weight
$2$--spheres and disks separate any components of $K$ from $L$;
so, we can discard these components. Since we must have $K$
separated from $L$, we have the desired normal surface.\end{proof}

Notice that in shrinking a surface, we do not increase the genus of the surface, even in the
bounded case. Hence, in the previous proposition, the  separating normal surface may
be found so that its genus is no more than the minimal genus of the surfaces $B_K$ and $B_L$. In
particular, if one of
$K$ or
$L$ is simply connected, then we can separate $K$ and $L$ by normal $2$--spheres.

Finally, we have

\begin{prop} Suppose $\T$ is a triangulation of the compact, orientable $3$--manifold $M$ and $S$ is a
closed, two-sided normal surface in $M$. Let $M'$ be the manifold
obtained by splitting $M$ along $S$. Suppose $D_1,\ldots,D_n$ is a
collection of pairwise disjoint, properly embedded disks in $M'$,
which are normal in the induced cell-decomposition on $M'$.
Furthermore, suppose the $D_i$ are all on the same side of $S$
(only meet $S'$, say, and so not meet $S''$ in $M'$). Then there
is a compression body $H$ embedded in $M$, $\bdy_+H = S$, each
component of $\bdy_-H$ is a normal surface in $M$ and
$D_1\cup\ldots\cup D_n\subset H$.\end{prop}

\begin{proof} The proof of this proposition follows along the very same lines as \ref{engulf-surface},
except we replace the graph $\Lambda$ in that argument with the subcomplex $D_1\cup\ldots\cup
D_n$.\end{proof}

More discussions of barriers and shrinking appear in
\cite{jac-let-rub1}.

\section{crushing triangulations}\label{sect-crush}

In the Introduction we pointed out that our techniques evolve from
the idea of finding a non vertex-linking, normal $2$--sphere,
which bounds a $3$--cell in our manifold, and then ``crushing" the
$2$--sphere (and $3$--cell) to a point. Generally speaking, this
is what we do; however,  while crushing a cell to a point gives us
back our manifold, in general, it wrecks havoc with the
triangulation.  Of course, crushing a $3$--cell, which is bounded
by a normal $2$--sphere, to a point induces a cell-decomposition
on the resulting manifold and there are many straight forward
techniques to construct a triangulation from such a
cell-decomposition. We emphasize that this is not what we do. Such
methods generally result in having a large number of tetrahedra
and our notion of a simpler triangulation is one with fewer
tetrahedra.  Our main task then, in this situation,  is to
reconstruct a triangulation of our manifold, which is simpler than
the original triangulation. In this section we organize such a
construction, which we think of as crushing a triangulation along
a normal surface. This has turned out to be a very useful concept.

Suppose $\T$ is a triangulation of a closed, orientable
$3$--manifold or an ideal triangulation of the interior of a
compact, orientable $3$--manifold $M$. Suppose $S$ is a normal
surface embedded in $M$ and $X$ is the closure of a component of
the complement of $S$; furthermore, suppose $X$ does not contain
any vertices of $\T$. In this situation, we give sufficient
conditions for constructing a particularly nice {\it ideal}
triangulation of $\open{X}$. We will organize our  construction
into a theorem at the end of this section.

\begin{figure}[htbp]
            \psfrag{X}{$X$}
            \psfrag{s}{tetrahedron}
            \psfrag{f}{face}
            \psfrag{e}{edge}
            \psfrag{t}{\begin{tabular}{c}
          truncated-tetrahedron\\
            Type I\\
            \end{tabular}}
            \psfrag{p}{\begin{tabular}{c}
            truncated-prism\\
         Type II\\
            \end{tabular}}
            \psfrag{q}{\begin{tabular}{c}
         parallel\\
          triangular\\
            Type III\\
            \end{tabular}}
            \psfrag{r}{\begin{tabular}{c}
         parallel\\
          quadrilateral\\
            Type IV\\
            \end{tabular}}

        \vspace{0 in}
        \begin{center}
        \epsfxsize=3 in
        \epsfbox{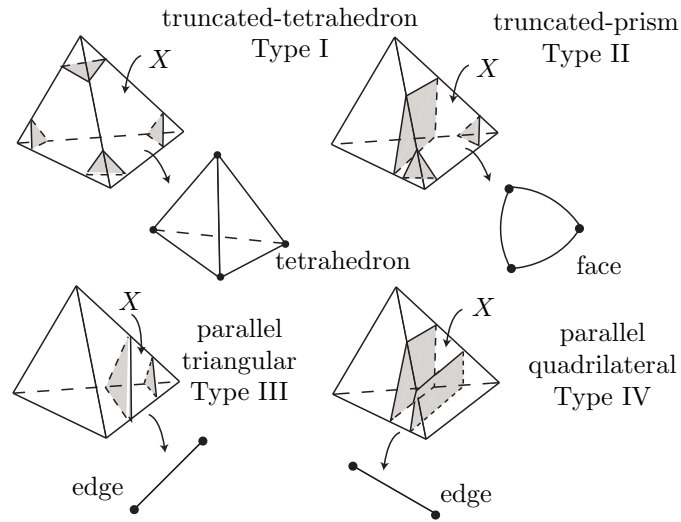}
        \caption{Cells in induced cell-decomposition of $X$.}
        \label{f-cell-decomp}
        \end{center}

\end{figure}

Notice that the manifold $X$ has a nicely described
cell-decomposition, say $\mathcal{C}$, induced from the
triangulation $\T$ and the way in which the  normal surface $S$
sits in $\T$ (no vertices are in $X$). The induced cells are of
four types:  truncated  tetrahedra, truncated prisms,  triangular
parallel regions, and quadrilateral parallel regions. We say the
truncated tetrahedra in $\mathcal{C}$ are {\it cells of type $I$},
the truncated prisms are {\it cells of type $II$}, the parallel
triangular regions are {\it cells of type $III$} and the parallel
quadrilateral regions are {\it cells of type $IV$}. See Figure
\ref{f-cell-decomp}.

 We show
that under the right conditions, one can get an ideal
triangulation of $\open{X}$ by essentially  replacing the
truncated tetrahedra in $\C$ by tetrahedra; replacing the
truncated prisms in $\C$ by faces and replacing the parallel
product cells by edges (again, see Figure \ref{f-cell-decomp}).
However, as one might expect, things are in general a bit more
complicated than this; hence, to see that something like this
works, we need to analyze the total structure of the collection of
parallel triangular and quadrilateral regions and truncated prisms
in $X$.

\begin{figure}[htbp]

            \psfrag{t}{trapezoid}
            \psfrag{p}{\begin{tabular}{c}
            truncated-prism\\
         Type II\\
            \end{tabular}}
            \psfrag{q}{\begin{tabular}{c}
         parallel\\
          triangular\\
            Type III\\
            \end{tabular}}
            \psfrag{r}{\begin{tabular}{c}
         parallel\\
          quadrilateral\\
            Type IV\\
            \end{tabular}}

        \vspace{0 in}
        \begin{center}
        \epsfxsize=3 in
        \epsfbox{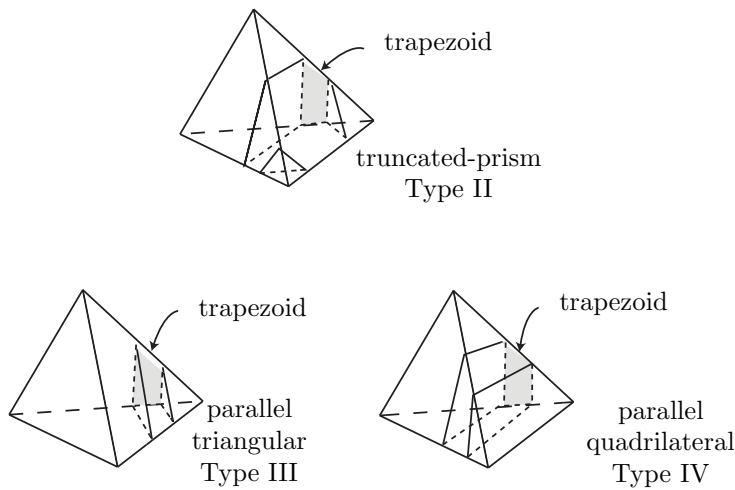}
        \caption{Trapezoidal faces in  the $2$--skeleton of the
cell-decomposition of $X$.}
        \label{f-trapezoid}
        \end{center}

\end{figure}

First, we define what we will call a product region for $X$, determined from the cell-decomposition
$\mathcal{C}$. Notice that in
 cells of type III and type IV  in $\C$  there are
quadrilateral faces in the $2$--skeleton of $\C$ (and in faces of
the tetrahedra of $\T$), which are complementary to the faces in
$S$. See Figure \ref{f-trapezoid}. Also, in the cells of type II,
the truncated prisms, there are quadrilateral faces in the
$2$--skeleton of $\C$ (and, again, in the faces of tetrahedra of
$\T$), which are complementary to the faces in $S$. We will call
these quadrilaterals {\it trapezoids} to distinguish them from the
normal quadrilaterals in $S$, which also are in the faces of cells
of types II and IV in $\C$. We let $\bbb{P}(\C)$ denote the
following union. $\bbb{P}(\C) = \{$edges of $\C$ not in $S\}\cup$
$\{$cells of type III and type IV in $\C\} \cup$ $\{$all
trapezoidal faces of $\C\}$. Each component of $\bbb{P}(\C)$ is an
$I$-bundle. In earlier drafts, we only considered cells of type
III and type IV in $\C$ in defining $\bbb{P}(\C)$; the method we
are using here was suggested by Nathan Dunfield, Marc Culler and
Peter Shalen and supported by a number of other colleagues.

Suppose $\bbb{P}(\C) \ne X$ and all of the components of
$\bbb{P}(\C)$ are product $I$-bundles. In applications, we will
need to establish that these conditions are satisfied. Under these
assumptions, each component of $\bbb{P}(\C)$, say $\bbb{P}_i$, for
some $i = 1,\ldots,k$, where $k$ is the number of components of
$\bbb{P}(\C)$, has the structure $\bbb{P}_i = K_i\times[0,1]$,
where $K_i$ is isomorphic to a subcomplex of the induced normal
cell structure on $S$.  Let $K_{i}^\ve = K_i\times \ve, \ve =
0,1$. Then $K_{i}^0$ and $K_{i}^1$ are disjoint, isomorphic
subcomplexes of the induced normal cell structure on $S$. Some of
these product components may be just edges in the $1$--skeleton of
$\C$ or made up entirely of trapezoids, in which case  $K_{i}^\ve,
\ve = 0, 1$ are points or graphs, respectively. The ideal
situation would be that each $K_{i}^\ve$ is simply connected and
therefore a cellular planar complex. (We have assumed
$\bbb{P}(\C)\ne X$.) While, in general, we can not have each
$K_{i}^\ve$ simply connected, in practice there is a construction
which alters the products $K_i\times [0,1]$ to a new collection of
products $D_j\times [0,1]$, where $D_j\times \ve, \ve = 0, 1$ is a
subcomplex of $S$ and the inclusion homomorphism $\pi_1(D_j\times
\ve)$ into $\pi_1(S)$ is injective for $\ve =0,1$.

\begin{figure}[htbp]

            \psfrag{S}{\large{$S$}}
            \psfrag{a}{\small{$D_4^0$}}
            \psfrag{b}{\small{$D_4^1$}}
            \psfrag{c}{\small{$D_3^0$}}
            \psfrag{d}{\small{$D_3^1$}}
            \psfrag{1}{\small{$K_1\times [0,1] = D_1\times [0,1]$}}
            \psfrag{2}{\small{$K_2\times [0,1]\subset D_2\times [0,1]$}}
            \psfrag{3}{\small{$K_3\times [0,1]\subset D_3\times
            [0,1]$}}
            \psfrag{4}{\small{$K_4\times [0,1]\subset D_4\times
            [0,1]$}}
            \psfrag{5}{\begin{tabular}{c}
            \small{$K_5\times [0,1] =$}\\
        \small{$D_5\times [0,1]$}\\
            \end{tabular}}
        \vspace{0 in}
        \begin{center}
        \epsfxsize=3 in
        \epsfbox{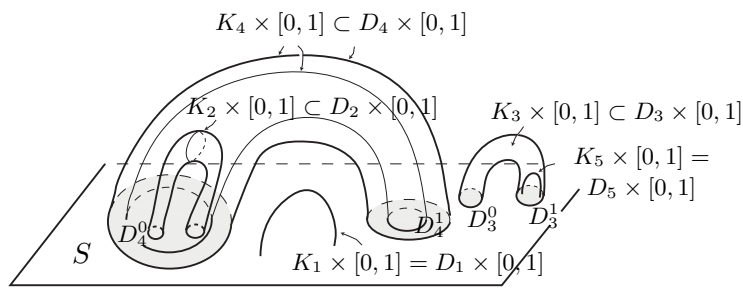}
        \caption{Product region $\bbb{P}(X)$ in $X$.}
        \label{f-product-new}
        \end{center}

\end{figure}

For now, we simply make this an assumption, which in practice will
need to be established. We  assume there is a pairwise disjoint
collection of spaces $D_j, 1\leq j\leq k'$, along with pairwise
disjoint embeddings $D_j\times [0,1]$ into $X$ so that $D_j\times
[0,1]$ is a subcomplex of $X$ and for $D_{j}^\ve = D_j\times \ve,
\ve =0,1$, we have $D_{j}^\ve$ is embedded as a subcomplex of the
induced cell structure on $S$; the inclusion of $\pi_1(D_{j}^\ve)$
into $\pi_1(S)$  is injective, $\ve = 0,1$; and, finally, $\bigcup
(K_i\times [0,1]) \subset \bigcup (D_j\times [0,1])$ and the
frontier of $\bigcup (D_j\times [0,1])$ is contained in the
frontier of $\bigcup (K_i\times [0,1])$. While  $K_{i}^0$ is
isomorphic to $K_{i}^1$ for every $i$, it will not necessarily be
the case that $D_{j}^0$ is isomorphic to $D_{j}^1$ for every $j$.
See Figure \ref{f-product-new}.

Under our assumptions, we let $\bbb{P}(X)$ be the collection of
components  of $\bigcup (D_j\times [0,1])$ and say that
$\bbb{P}(X)$ is an {\it induced product region for $X$}. We note,
in general, there is not a unique induced product region. Each
component of $\bbb{P}(\C)$ is contained in a component of
$\bbb{P}(X)$. It is possible that some $D_{j'}\times [0,1]$ is
contained in a $D_j\times [0,1], j'\ne j$; so, $D_{j'}\times
[0,1]$ is consumed and does not appear as a component of
$\bbb{P}(X)$. Also, it is possible that a number of truncated
prisms and truncated tetrahedra in the cell-decomposition $\C$ on
$X$ are also consumed into the components of $\bbb{P}(X)$. The
product structure on the components of $\bbb{P}(X)$, however,
respects the combinatorial structure of the cells of $\C$ in the
sense that each component of the frontier of $D_j\times [0,1]$ is
either an edge in $\C$ (when $D_j\times [0,1]$ is itself an edge)
or a collection of trapezoids. We assumed $\bbb{P}(\C) \ne X$; in
practice, this will give that $\bbb{P}(X) \ne X$. Of course, this
is something else we will have to establish.

If each component of $\bbb{P}(X)$ is a product $D_i\times [0,1]$
where $D_i$ is a simply connected planar complex (is cellular),
then we say the product region $\bbb{P}(X)$ for $X$ is a {\it
trivial product region for $X$}.

Now, consider the cells of $X$ of type $II$, truncated prisms.
Each cell of type II has two hexagonal faces. In the cell
decomposition of $X$, these hexagonal faces can be identified to
either a hexagonal face of a cell of type $I$ (truncated
tetrahedron) or a hexagonal face of another cell of type $II$. If
we follow a sequence of such identifications through cells of type
$II$, we trace out a well-defined arc (see Figure \ref{f-chain}
below) which either terminates at the identification of a
hexagonal face of a cell of type $II$ with one of type $I$ or it
forms a simple closed curve through cells of type $II$, ending at
the cell in which it started. We call the collection of truncated
prisms identified in this way a {\it chain}. If a chain ends in a
truncated tetrahedra, we say the chain {\it terminates};
otherwise, we call the chain a {\it cycle}.  Notice, in
particular, if for some truncated prism there is an identification
of one of its hexagonal faces with the other, then we have a cycle
of length one.

\begin{figure}[htbp]

            \psfrag{o}{\large{or}}
            \psfrag{t}{\begin{tabular}{c}
            truncated tetrahedron\\
       chain terminates\\
            \end{tabular}}
             \psfrag{p}{\begin{tabular}{c}
            truncated prism\\
       chain continues\\
            \end{tabular}}
        \vspace{0 in}
        \begin{center}
        \epsfxsize=3 in
        \epsfbox{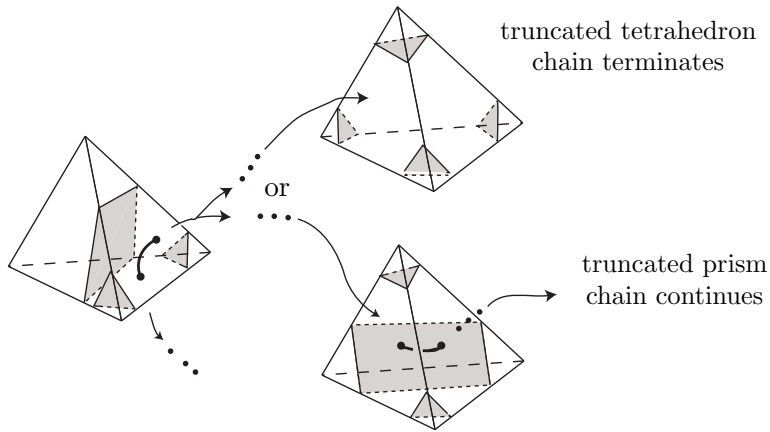}
        \caption{Chain of truncated prisms.}
        \label{f-chain}
        \end{center}

\end{figure}

We will assume there are no cycles of truncated prisms, which are
not in the induced product region of $X$, $\bbb{P}(X)$. Later, of
course, we will have the burden of proof to establish this in
applications of crushing.

\begin{figure}[htbp]

          \psfrag{o}{of prisms}
           \psfrag{c}{chain}
             \psfrag{I}{$\overline{\Delta}_i$}
            \psfrag{J}{$\overline{\Delta}_j$}
            \psfrag{i}{$\td{\Delta}_i^*$}
            \psfrag{j}{$\td{\Delta}_j^*$}
            \psfrag{1}{$\overline{\sigma}_i$}
            \psfrag{2}{$\overline{\sigma}_j$}
            \psfrag{3}{$\td{\sigma}_j$}
            \psfrag{4}{$\td{\sigma}_i$}
            \psfrag{a}{\begin{tabular}{c}
            after\\
        crush\\
            \end{tabular}}
            \psfrag{f}{\begin{tabular}{c}
            faces\\      identified\\
            \end{tabular}}
       \vspace{0 in}
        \begin{center}
        \epsfxsize = 3.5 in
        \epsfbox{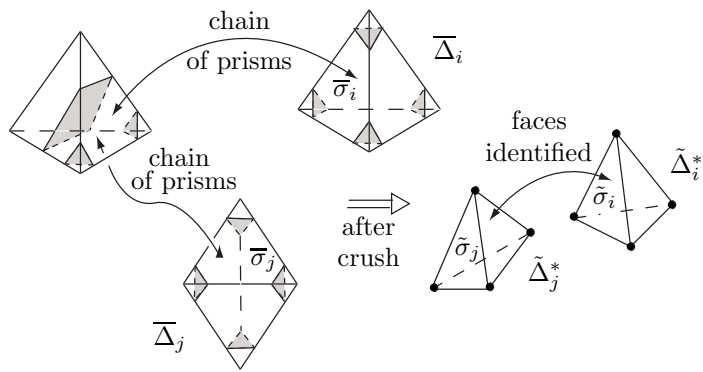}
        \caption{Face identifications induced through a chain of truncated
        prisms.}
        \label{f-ident-prisms}
        \end{center}

\end{figure}

By our assumptions, there are truncated tetrahedra in $X$, which
are not in $\bbb{P}(X)$ ($X\ne \bbb{P}(X)$), the frontier of
$\bbb{P}(X)$ consists of trapezoids, and there is no cycle of
truncated prisms, which is not in $\bbb{P}(X)$). Let
$\{\overline{\Delta}_1,\ldots,\overline{\Delta}_n\}$ be the
collection of truncated tetrahedra in $X$, which are not in
$\bbb{P}(X)$. Now, a face of a truncated tetrahedron in $\C$ is
either shared by two truncated tetrahedra in $\C$ or is also a
hexagonal face of a truncated prism in $\C$. In the last case,
suppose we have a hexagonal face, say $\overline{\sigma}_i$ of the
truncated tetrahedron $\overline{\Delta}_i$, which also is a
hexagonal face in the first truncated prism in a chain of
truncated prisms. If we follow the chain of truncated prisms,
there is a last hexagonal face, say $\overline{\sigma}_j$, which
also is in a truncated tetrahedron, $\overline{\Delta}_j$,
possibly, $i = j$. Hence, there is an induced identification of
the face $\overline{\sigma}_i$ of $\overline{\Delta}_i$ with the
face $\overline{\sigma}_j$ of $\overline{\Delta}_j$, through the
chain of truncated prisms. So, the faces of the truncated
tetrahedra in $\{\overline{\Delta}_1,\ldots,\overline{\Delta}_n\}$
have an induced pairing. See Figure \ref{f-ident-prisms}. Now,
notice that each truncated tetrahedron in $X$ has its triangular
faces in $S$. We can identify each such triangular face to a point
(distinct points for each triangular face) and we get tetrahedra.
We will now use the notation $\td{\Delta}_i^*$ for the tetrahedron
coming from the truncated tetrahedron $\overline{\Delta}_i$ by
identifying each of the triangular faces of $\overline{\Delta}_i$
to a point (distinct points for each triangular face) and use
$\td\sigma_i$ for the triangular face coming from the hexagonal
face $\overline{\sigma}_i$. Then $\boldsymbol{\Delta}^*(X)
=\{\td{\Delta}^*_1,\ldots, \td{\Delta}^*_n\}$ is a collection of
tetrahedra with orientation induced by that on $\T$ and the
induced pairings between the triangular faces $\td\sigma_i$ and
$\td\sigma_j$ described above is a family $\boldsymbol{\Phi}^*$ of
orientation reversing affine isomorphisms. Hence, we get a
$3$--complex $\boldsymbol{\Delta}^*(X)/\boldsymbol{\Phi}^*$, which
is a $3$--manifold except, possibly, at its vertices. We will
denote the associated ideal triangulation by $\T^*$ and the
underlying point set by $\abs{\T^*}$. We call $\T^*$ the ideal
triangulation obtained by {\it crushing the triangulation $\T$
along $S$}. We denote the image of a tetrahedron $\td{\Delta}^*_i$
by $\Delta^*_i$ and, as above, call  $\td{\Delta}^*_i$ the lift of
$\Delta^*_i$.

We have the following theorem.

\begin{thm}\label{crush} Suppose $\T$ is a triangulation of a closed, orientable $3$--manifold
or an ideal triangulation of the interior of a compact, orientable
$3$--manifold $M$. Suppose $S$ is a normal surface embedded in
$M$, $X$ is the closure of a component of the complement of $S$
and $X$ does not contain any vertices of $\T$. Suppose there is an
induced product region, $\bbb{P}(X)$, for $X$. If
\begin{enumerate}
\item[i)] $X\ne \bbb{P}(X)$, \item[ii)] $\bbb{P}(X)$ is a trivial
product region for $X$, and \item[iii)] there are no cycles of
truncated prisms in $X$, which are not in $\bbb{P}(X)$,
\end{enumerate}
then the triangulation $\T$ can be crushed along $S$ and $\T^*$ is
an ideal triangulation of $\open{X}$.
\end{thm}

\begin{figure}[htbp]

            \psfrag{e}{$e$}
            \psfrag{x}{$D_1\times I$}
            \psfrag{z}{$D_i\times I$}
            \psfrag{y}{$D_n\times I$}
            \psfrag{1}{$e_1$}
            \psfrag{i}{$e_i$}
            \psfrag{n}{$e_n$}
            \psfrag{v}{$v$}
            \psfrag{c}{``crush"}
            \psfrag{S}{\large{$S$}}
            \psfrag{s}{\begin{tabular}{c}
            $S \rightarrow v$\\
       and\\
       $D_j\times I\rightarrow e_j$
            \end{tabular}}
        \vspace{0 in}
        \begin{center}
        \epsfxsize = 4 in
        \epsfbox{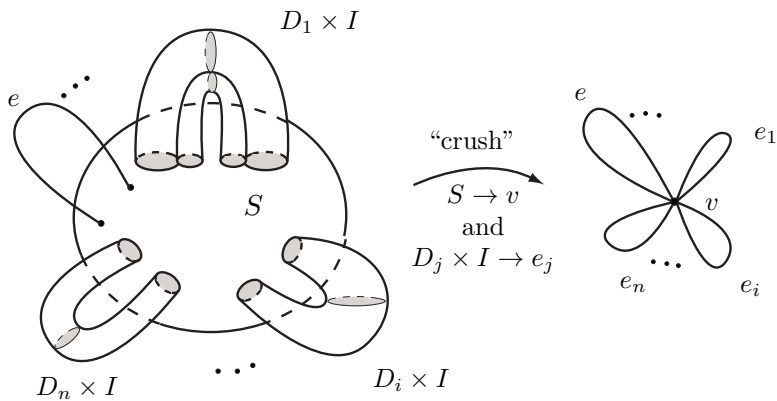}
        \caption{Components of trivial product region ``crushed" to edges.}
        \label{f-crush-product}
        \end{center}

\end{figure}

\begin{proof} The underlying point set $\abs{\T^*}$ is obtained from $X$ by identifying each
component of $S$ to a point (distinct points for distinct
components), identifying each component $D_i\times [0,1]$ of
$\bbb{P}(X)$, the product region for $X$, to an edge $e_i =
[0,1]_i$ (distinct edges for distinct components; see Figure
\ref{f-crush-product}), and identifying each chain of truncated
prisms to a face (see Figure \ref{f-ident-prisms}). If we look at
this identification map we have the inverse image of a point in
the interior of a tetrahedron $\Delta^*_i$ is just a point in the
interior of the truncated tetrahedron $\overline{\Delta}_i$; the
inverse image of a point in the interior of a face is either a
point or an arc, the latter in the case a chain of truncated
prisms is identified to a face; and the inverse image of a point
in the interior of an edge is a copy $D_j\times x$ for some $j$
and $x\in [0,1]$. Notice that in the identification of a chain of
truncated prisms to a face; the associate identification of the
edges is through a band of trapezoids and so there are no new
identifications not already made in $D_j\times [0,1]$ for some
$j$. Thus the identification map on $\open{X}$ is a cell-like map.
It follows by \cite{arm, sie}, that $\T^*$ is an ideal
triangulation of $\open{X}$.
\end{proof}

We end this section with a couple of observations. First, the
tetrahedra of the ideal triangulation $\T^*$ are in one-one
correspondence with the truncated tetrahedra of the
cell-decomposition $\C$, which are not in the induced product
region $\bbb{P}(X)$. Hence, if $t$ is the number of tetrahedra of
$\T$ and $t^*$ is the number of tetrahedra of $\T^*$, then $t^*\le
t$. The inequality is strict unless $S$ is  vertex-linking; in
fact,  $\T^* = \T$ if and only if $S$ is vertex-linking. If each
component of $S$ is a $2$--sphere and $\widehat{X}$ is the
$3$--manifold obtained by capping off each $2$--sphere in $\bdy
X$, then $\T^*$ is a triangulation of $\widehat{X}$. More in the
spirit of crushing, if each component of $S$ is a $2$--sphere, we
can think of $\abs{\T^*}$ as the the manifold obtained from $X$ by
identifying each $2$--sphere in $S$ to a point (distinct points
for distinct $2$--spheres). This is the same manifold as
$\widehat{X}$.

\section{$0$--efficient triangulations}\label{sect-0-eff}

In this section we develop the concepts of $0$--efficient
triangulations of closed, orientable $3$--manifolds and compact,
orientable $3$--manifolds with boundary. In a later section we
define and study $0$--efficient ideal triangulations. A
$0$--efficient triangulation severally limits the nature of
embedded normal $2$--spheres (closed manifolds and ideal
triangulations) and normal disks (bounded manifolds). Also, in
this section, and in a later section on minimal triangulations, we
show that with few exceptions, $0$--efficient triangulations do
not have edges of order less than four.

\subsection{$\mathbf{0}$--efficient triangulations for closed $\mathbf
3$--manifolds} 
\begin{figure}[htbp]

            \psfrag{i}{identified}
            \psfrag{L}{$L(3,1)$}
            \psfrag{B}{\begin{tabular}{c}
            (B) Solid torus and $3$--cell;\\
     not $0$-efficient\\
            \end{tabular}}
             \psfrag{A}{\begin{tabular}{c}
            (A) Two solid tori;\\
       $0$--efficient\\
            \end{tabular}}
        \vspace{0 in}
        \begin{center}
        \epsfxsize=3.5 in
        \epsfbox{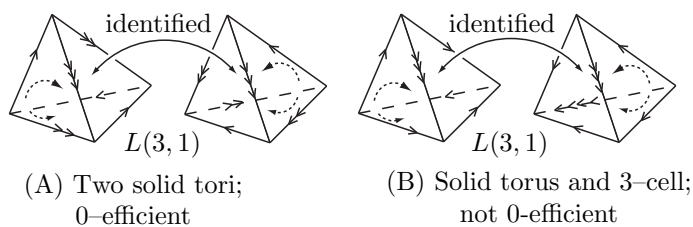}
        \caption{Two tetrahedron, one-vertex triangulations of $L(3,1)$.
(A) is $0$--efficient. (B) is not $0$--efficient.}
        \label{f-two-L3_1}
        \end{center}

\end{figure}

\begin{figure}[htbp]

            \psfrag{A}{A}
            \psfrag{B}{B}
            \psfrag{R}{$\mathbb{R}P^3$}
            \psfrag{P}{$2\mathbb{R}P^2 = S^2$}
            \psfrag{S}{$S^2$}
            \psfrag{L}{$L(3,1)$}
        \vspace{0 in}
        \begin{center}
        \epsfxsize = 3 in
        \epsfbox{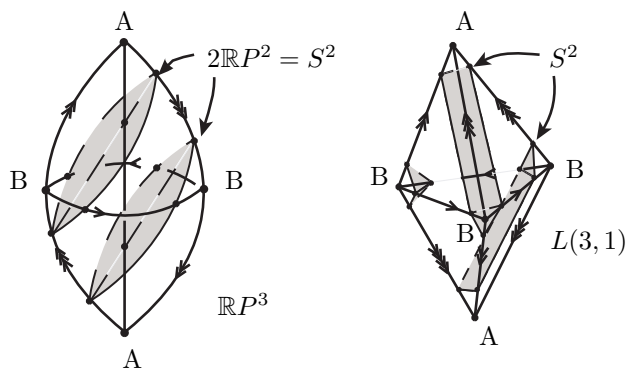}
        \caption{Standard lens space representations giving two-tetrahedron,
         two-vertex triangulations of $\rp$ and
$L(3,1)$, neither of which is $0$--efficient.}
        \label{f-not-0-eff}
        \end{center}

\end{figure}

A triangulation of a closed $3$--manifold is said
to be {\it $0$--efficient} if and only if the only embedded,
normal $2$--spheres are vertex-linking. For example, both of the
one-tetrahedron triangulations of the $3$--sphere are
$0$--efficient (see Figure \ref{f-tetra} (4) and (5)), while only
one of the two-tetrahedra, one-vertex triangulations of $L(3,1)$
is $0$--efficient, see  Figure \ref{f-two-L3_1}. Neither of the
two-tetrahedra, two-vertex triangulations of $L(3,1)$ are
$0$--efficient, which includes the two tetrahedron standard lens
space presentations of
 $L(3,1)$ (see Figure \ref{f-not-0-eff}). No triangulation of
$\rp$ can be $0$--efficient; there are two two-tetrahedra
triangulations of $\rp$, one has one vertex and the other is the
standard lens space presentation of $\rp$ which has two vertices
(again, see Figure \ref{f-not-0-eff}). In \cite {jac-let-rub1}, it
is shown that $S^3$ and all lens spaces, except $\rp$, admit
infinitely many one-vertex $0$--efficient triangulations. In fact,
it can be shown that if there is a $0$--efficient triangulation of
the $3$--manifold $M$ which contains a layered solid torus, then
there are infinitely many $0$--efficient triangulations of $M$.
However, there are manifolds for which the only $0$--efficient
triangulation we know does not have a layered solid torus as a
subcomplex. We suspect that any manifold which admits a
$0$--efficient triangulation admits infinitely many such
triangulations but we have not been able to confirm this in
general.

The following proposition is quite useful and provides
some insight to
$0$--efficient triangulations. We wish to thank Eric Sedgwick for helpful
comments regarding the proof of this proposition.

\begin{prop} \label{p-0-eff} Suppose $M$ is a closed, orientable
$3$--manifold. If
$M$ has a
$0$--efficient triangulation, then $M$
is irreducible and
$M\neq\mathbb{R}P^3$. Furthermore, either the triangulation has
one vertex or $M$ is $S^3$ and the triangulation has precisely two
vertices. \end{prop}
\begin{proof} Suppose $\T$ is a 0-efficient triangulation of $M$. If $M$ were not
irreducible, then
$M$ would contain an essential
$2$--sphere, a $2$--sphere that does not bound a $3$--cell. If this is the case,
then as in Theorem \ref{normalsphere}, also \cite{kne, haken1, jac-rub4}, for
any triangulation, in particular $\T$, $M$ has an embedded, essential, normal
$2$--sphere. However, an embedded, essential, normal $2$--sphere can not be
vertex-linking. So, $M$ must be irreducible.

Now, if $M$ contains an embedded $\rpp$, then, similarly, we have
for any triangulation of $M$ there is an embedded, normal $\rpp$.
Since an $\rpp$ can not be embedded in a $3$--cell, such an $\rpp$
must contain a normal quadrilateral. Hence, its double, which is
an embedded normal $2$--sphere, can not be vertex-linking.  See
Figure \ref{f-not-0-eff} where there is an embedded normal $\rpp$
with two normal quads and its double, which is an embedded, normal
$2$--sphere, has four normal quads.  This completes the proof of
the first part of the proposition.

So, suppose $\T$ is $0$--efficient and has more than one vertex. In this
case, we first show
$M$ is the
$3$--sphere.

Since we are assuming $\T$ has more than one vertex, there is an
edge $e$ in $\T$ which has distinct vertices. Let $N(e)$ denote a
small neighborhood of $e$ and let $S$ denote $\bdy N(e)$; then
$N(e)$ is a $3$--cell and $S$ is a $2$--sphere. If $S$ is normal,
then $S$ is vertex-linking and so, must bound a $3$--cell
complementary to the $3$--cell $N(e)$, which contains the edge
$e$. It follows in this case that $M$ is the $3$--sphere. So,
suppose $S$ is not normal. Then $S$ forms a barrier surface in the
component of its complement not containing $e$; so, we can shrink
$S$. In shrinking $S$ we end up with a punctured $3$--cell
containing the edge $e$ and its boundary is made up of normal
$2$--spheres and $0$--weight $2$--spheres contained entirely in
the interiors of tetrahedra of $\T$. Let $P$ denote this punctured
$3$--cell containing $e$. Any $2$--sphere in $\bdy P$, which is
embedded in a tetrahedron, bounds a $3$--cell in that tetrahedra
which does not meet $e$, which can be added to $P$, still giving a
punctured $3$--cell; and any normal $2$--sphere in $\bdy P$, which
is vertex-linking, bounds a $3$--cell not meeting $e$, which can
be added to $P$, giving a punctured $3$--cell. Thus each boundary
component of $P$ bounds a $3$--cell in the complement of $P$. It
follows that $M$ is $S^3$.

We only have left to prove that a $0$--efficient triangulation of
$S^3$, not having just one vertex, has precisely two. We establish
this through a series of claims.

\begin{figure}[htbp]

            \psfrag{0}{$v_0$}
            \psfrag{1}{$v_1$}
            \psfrag{2}{$v_2$}
            \psfrag{3}{$v_3$}
            \psfrag{N}{$\bdy N(e) = S^2$}
            \psfrag{O}{OR}
        \vspace{0 in}
        \begin{center}
        \epsfxsize = 3 in
        \epsfbox{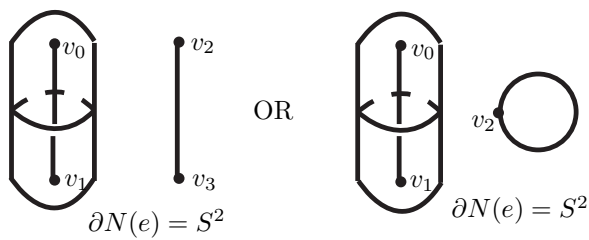}
        \caption{No disjoint edges in $\T$ with one an arc.}
        \label{f-0-eff1}
        \end{center}

\end{figure}

Suppose $\T$ is a $0$--efficient triangulation of $S^3$ with more than one
vertex.

\vspace{.15 in} \noindent {\it Claim. The $1$--skeleton of $\T$
can not have disjoint edges, one of which is an arc}. See Figure
\ref{f-0-eff1}. \vspace{.1 in}

Suppose $e$ and $e'$ are disjoint edges of $\T$. By Proposition
\ref{doublebarrier}, $e$ and $e'$ can be separated by a normal
surface; furthermore, since $e$, say, is an arc, then such a
normal surface can be taken to be a normal $2$--sphere. However
such a normal $2$--sphere can not be vertex-linking. Thus we
arrive at a contradiction.

\vspace{.15 in} \noindent {\it Claim. The $1$--skeleton of $\T$
can not have a cycle of three embedded edges between three
distinct vertices; i.e., there are no triads in the $1$--skeleton
of $\T$}. See
 Figure \ref{f-0-eff2}.
\vspace{.1 in}

Suppose $v_0,v_1,v_2$ are three distinct vertices determining the
cycle of distinct edges $\overline{v_0v_1}$, $\overline{v_1v_2}$,
$\overline{v_2v_0 }$. If there were more than three vertices, then
there is an edge $e'$ from a vertex $v$, distinct from each of the
vertices $v_0,v_1,v_2$, to one of these vertices, say $v_0$. Then
$e'$ and $e = \overline{v_1v_2}$ are disjoint edges, both of which
are arcs. This contradicts our previous observation (see Figure
\ref{f-0-eff2}). So, we may suppose that if there is an embedded
triad in the $1$--skeleton of $\T$, then $\T$ has at most three
vertices. If there were only three vertices, then since a
tetrahedron has four vertices, there must be an edge $e'$ of $\T$
having both its vertices at the same vertex of the triad with
vertices $v_0,v_1,v_2$, say at $v_0$. Again, we have an edge $e'$
disjoint from the edge $e = \overline{v_1v_2}$ which is an arc
and, hence, gives a contradiction to the previous Claim (again,
see Figure \ref{f-0-eff2}). This shows that there are no embedded
triads in the $1$--skeleton of $\T$.

\begin{figure}[htbp]

            \psfrag{0}{$v_0$}
            \psfrag{1}{$v_1$}
            \psfrag{2}{$v_2$}
            \psfrag{e}{$e$}
            \psfrag{d}{$e'$}
            \psfrag{v}{$v$}
            \psfrag{D}{\begin{tabular}{c}
            Disjoint edges;\\
    both are arcs\\
            \end{tabular}}
             \psfrag{E}{\begin{tabular}{c}
            Disjoint edges;\\
       only one an arc\\
            \end{tabular}}
        \vspace{0 in}
        \begin{center}
        \epsfxsize=3 in
        \epsfbox{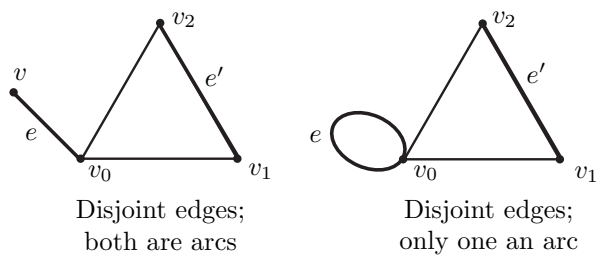}
        \caption{No embedded triads in $1$-skeleton of $\T$.}
        \label{f-0-eff2}
        \end{center}

\end{figure}

\vspace{.15 in} \noindent {\it Claim. If $\T$ has more than two
vertices, then there is a vertex $v_0$ of $\T$ so that every edge
has $v_0$ as a vertex}. See Figure \ref{f-0-eff3}. \vspace{.1 in}

Since we are assuming $\T$ has at least three vertices and we have
shown that no two edges, one of which is an arc, are disjoint, we
can choose notation so there are distinct vertices $v_0, v_1, v_2$
and edges $\overline{v_{0}v_1}$ and $\overline{v_{0}v_2}$, meeting
only in $v_0$. Suppose some edge $e'$ does not meet $v_0$. Then
$e'$ can not have both $v_1$ and $v_2$ as vertices, since this
would form an embedded triad in the $1$--skeleton of $\T$. Hence,
either $e'$ is disjoint from the two edges $\overline{v_{0}v_1}$
and $\overline{v_{0}v_2}$ or $e'$ meets just one of the vertices
$v_1$ or $v_2$; $e'$ may be a loop. In either case, we have two
disjoint edges, $e'$ and one of $\overline{v_{0}v_1}$ or
$\overline{v_{0}v_2}$, each of which is an arc. This contradicts
an earlier claim.

\begin{figure}[htbp]

            \psfrag{0}{$v_0$}
            \psfrag{1}{$v_1$}
            \psfrag{2}{$v_2$}
            \psfrag{n}{$v_n$}
        \vspace{0 in}
        \begin{center}
        \epsfxsize=1.5 in
        \epsfbox{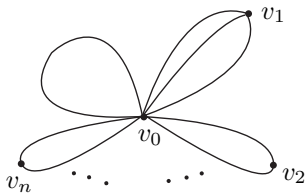}
        \caption{All edges of $\T$ meet at $v_0$.}
        \label{f-0-eff3}
        \end{center}

\end{figure}

\vspace{.15 in} \noindent {\it Claim. If every edge in the
$1$--skeleton of $\T$ meets a single vertex and there is more than
one vertex, then $\T$ is not $0$--efficient}. See Figure
\ref{f-0-eff4}. \vspace{.1 in}

Suppose every edge meets the vertex $v_0$ and there are at least two vertices,
say
$v_1$ is a second vertex. Since each edge has a vertex at $v_0$, there is an
edge $e$ having
$v_0$ and
$v_1$ as vertices. Furthermore, every tetrahedron having $v_1$ as a vertex has
all its other vertices at $v_0$.

There are two obstructions to building an edge-linking normal
$2$--sphere about an edge in a triangulation of a $3$--manifold.
One is that the ends of the edge are the same vertex; so, we
choose an edge from $v_1$ to $v_0$, say $e = \overline{v_1v_0}$.
The other obstruction is that, say for the edge $e$, there is a
tetrahedron $\tau$ in $\T$ having $e$ as an edge and $\tau$ has an
adjacent edge $e'$, also from $v_1$ to $v_0$,  identified with
$e$; this prevents a small neighborhood about the edge $e$ from
being normal. However, if this is the case, we will show how to
use the induced triangulation of the vertex-linking normal
$2$--sphere, $S_1$, about $v_1$ to find another edge, say $e^*$,
from $v_1$ to $v_0$ that does not have this obstruction. See
Figure \ref{f-0-eff4}. In the identification of the tetrahedron
$\tau$ we assume that the edge $e' = e$. Let $\sigma$ denote the
normal triangle at the vertex $v_1$ in $S_1$. The edge $a$ from
$e$ to $e'$ in this triangle is a loop in the induced
triangulation on $S_1$. But since we are on a $2$--sphere, there
must always be an ``innermost" loop, which is also formed by two
adjacent edges in a tetrahedron from $v_1$ to $v_0$ being
identified. We may just assume this is the situation for $e$ and
$e'$ in $\tau$; i.e., $a$ is an ``innermost" loop on $S_1$. We
will also continue to use $\sigma$ for the triangle in the
vertex-linking $2$--sphere with vertices on $e$ and $e'$. In this
situation, an edge $e^*$ going through any vertex in the interior
of the loop determined by $a$ on the vertex-linking $2$--sphere
does not have an edge in $S_1$ a loop and so does not have an
obstruction to building an embedded normal sphere about that edge.
So, there is an embedded, non vertex-linking, normal $2$--sphere.
This gives the desired contradiction and completes the proof of
the Claim.

\begin{figure}[htbp]

            \psfrag{0}{$v_0$}
            \psfrag{1}{$v_1$}
            \psfrag{e}{$e$}
            \psfrag{d}{$e'$}
            \psfrag{f}{$e^*$}
            \psfrag{g}{$e = e'$}
            \psfrag{t}{$\tau$}
            \psfrag{S}{$S_1$}
            \psfrag{s}{$\sigma$}
            \psfrag{a}{$\alpha$}
            \psfrag{A}{$e'\ne e$ for all $\tau$}
            \psfrag{B}{$e' = e$ for all $\tau$}
            \psfrag{h}{$e\ne e'$}
        \vspace{0 in}
        \begin{center}
        \epsfxsize=4 in
        \epsfbox{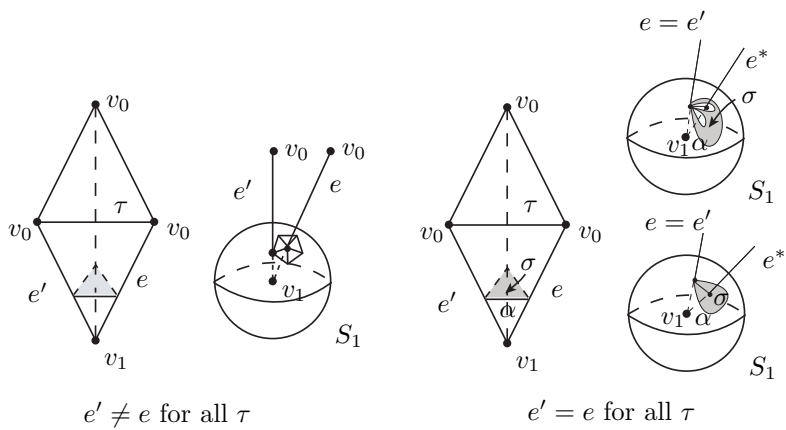}
        \caption{If all edges meet at $v_0$ and there is more than one
vertex, then $\T$ is not $0$--efficient.}
        \label{f-0-eff4}
        \end{center}

\end{figure}

We have shown that if $\T$ has more than two vertices, then the triangulation
must have a special vertex which meets every edge; however, if every edge meets
a vertex and there is more than one vertex, then it leads to a contradiction of
$0$--efficiency.

This completes the proof of the Proposition. \end{proof}

Ben Burton has constructed an infinite family of two vertex,
$0$--efficient triangulations of $S^3$ (see \cite{jac-let-rub1}).
The one-tetrahedron, two-vertex triangulation of $S^3$ is dual to
the natural cell complex coming from the genus one Heegaard
splitting of $S^3$ with the two meridian disks attached to the
Heegaard surface. Triangulations dual to cell subdivisions induced
by a Heegaard surface with a system of meridional disks attached
have precisely two vertices. By Proposition \ref{p-0-eff}, except
possibly for $S^3$, these are never $0$--efficient.

 There is  almost a converse for two-vertex,
$0$--efficient triangulations of $S^3$. Suppose $\T$ is a
two-vertex triangulation of the $3$--manifold $M$ with vertices
$v_1$ and $v_2$. For each tetrahedron $\Delta$ of $\T$ which has
both $v_1$ and $v_2$ as vertices there is either a uniquely
determined normal triangle or normal quadrilateral in $\Delta$,
separating the vertices identified with $v_1$ from those
identified with $v_2$. The collection of all such normal triangles
and quadrilaterals forms a normal surface $S$. $S$ separates the
vertex $v_1$ from the vertex $v_2$ and if $S$ meets an edge of
$\T$ at all, then it meets it in at most one point. This
construction works no matter how many vertices, as long as we
divide the vertices into two disjoint sets. However, it is
probably most interesting in the case of two vertices. In fact, if
$K$ and $L$ are disjoint subcomplexes, then as observed in
Proposition \ref{doublebarrier}, there is a normal surface in $M$
separating $K$ from $L$. There is not, in general, a unique such
surface as one gets in the above situation for a triangulation
with two vertices. We call $S$ a {\it vertex-splitting surface}.

\begin{lem} If $\T$ is a two-vertex, $0$--efficient triangulation of $S^3$, then the
vertex-splitting surface
$S$ determines a Heegaard splitting of $S^3$.
\end{lem}

\begin{proof} Choose an edge $e$ of $\T$ joining $v_1$ and $v_2$. As observed above, $e$
meets $S$ precisely once. A small  regular neighborhood of $S\cup
e$ has two boundary components, each parallel to $S$ and each
barrier surfaces in the components of their complements not
meeting $S\cup e$. So, each can be shrunk in these components of
their complements, resulting in two compression bodies, $H$ and
$H'$, where $\bdy_+H = \bdy_+H' = S$, $\bdy_-H$ and $\bdy_-H'$ are
normal surfaces along with, possibly, some $0$--weight
$2$--spheres contained entirely in the interior of tetrahedra of
$\T$; furthermore, $e$ is contained in $H\cup H'$.

We claim that there are no components in either $\bdy_-H$ or in
$\bdy_-H'$ which are normal. For suppose there were normal
components. Let $N(e)$ be a small $3$--cell neighborhood of $e$
(disjoint from $\bdy_-H\cup\bdy_-H'$) with boundary the
$2$--sphere $\Sigma$. Now, assuming there are normal surfaces in
$\bdy_-H$ or in $\bdy_-H'$, $\Sigma$  separates these normal
surfaces from $e$. Furthermore, $\Sigma$ along with
$\bdy_-H\cup\bdy_-H'$ are barrier surfaces. Thus we can shrink
$\Sigma$ in the component of the complement of
$\Sigma\cup\bdy_-H\cup\bdy_-H'$ whose closure contains $\Sigma$
but not $e$. But we then must have a punctured $3$--cell
containing $e$ and having boundary a collection of normal
$2$--spheres and $0$--weight $2$--spheres entirely contained in
the interior of tetrahedra of $\T$, which separate $e$ from the
normal components of $\bdy_-H\cup\bdy_-H'$. It follows that if
$\bdy_-H\cup\bdy_-H'$ had any normal surfaces, then there must be
a normal $2$--sphere in the boundary of this punctured $3$--cell.
But this is impossible for $\T$ a $0$--efficient triangulation.

It follows that  both $\bdy_-H\cup \bdy_-H' $ is either empty or consists entirely of
$2$--spheres contained entirely in the interior of tetrahedra of $\T$. Any such $2$--sphere
bounds a $3$--cell complementary to $H$ and $H'$. Hence, we may assume that we have filled any
such boundary components and therefore
$H$ and
$H'$ are handlebodies. Hence, $S$ is a Heegaard surface, as claimed.\end{proof}

We point out that in the case, cited above, of the Burton examples of two-vertex,
$0$--efficient triangulations, the vertex-splitting surface determines a genus one Heegaard
splitting in every case. We have not tried to get two-vertex, $0$--efficient triangulations of
$S^3$ with higher genus splitting surfaces nor have we tried to understand those conditions on
a Heegaard splitting of $S^3$ which guarantee that the two-vertex triangulation dual to the
cell decomposition coming from the Heegaard surface along with a complete system of meridional
disks is $0$--efficient.

We have the following observations concerning  edges in $0$--efficient
triangulations.

\begin{prop} Suppose the closed, orientable $3$--manifold $M$ has a $0$--efficient triangulation
$\T$. If
$\T$ has an edge bounding an embedded disk in $M$, then $M = S^3$.
\end{prop}
\begin{proof} Suppose the edge $e$ bounds an embedded disk $D$ (in particular, the edge has
just one vertex and is an embedded simple closed curve). Let $S$
be the boundary of a small regular neighborhood of $D$. Then $S$
is a $2$--sphere and bounds a $3$--cell containing the edge $e$.
Let $T$ be the boundary of a small regular neighborhood of $e$
which is contained entirely in the $3$--cell bounded by $S$. Then
$T$ is a barrier surface in the component of its complement not
containing $e$; hence, in the component containing $S$. So, we can
shrink $S$ obtaining a punctured $3$--cell, which contains the
edge $e$ and has boundary a number of $2$--spheres, each of which
is either normal or contained in the interior of a tetrahedron. We
may fill in those in tetrahedra with $3$--cells missing $e$. On
the other hand, if there is a normal $2$--sphere then it must be
vertex linking. It can not be linking the vertex of $e$; so, we
conclude that it bounds a $3$--cell not meeting $e$. Hence $M$ is
$S^3$. Of course, if there were a normal $2$--sphere after
shrinking, then there would necessarily be more than one vertex
and we have from Proposition \ref{p-0-eff} that $M$ is $S^3$.
\end{proof}

\begin{cor}\label{nodisk} Suppose the closed, orientable $3$--manifold $M$ has a $0$--efficient
triangulation
$\T$.
\begin{enumerate}
\item If $\T$ has an edge of order one, then $M = S^3$.
\item If $\T$ has a face which is a cone, then $M = S^3$.
\end{enumerate}
\end{cor}

\begin{proof} If there is an edge, say $e$, of order one, then there is a single tetrahedron
$\td{\Delta}$ in our triangulation and a single edge $\td{e}$ in $\td{\Delta}$, which is
the full collection of edges identified to $e$. Let $\td{e}'$ be the edge in $\td{\Delta}$ dual
to
$\td{e}$; i.e.,
$\td{\Delta} = \td{e}
\ast \td{e}'$. Then if $e'$ is the image of $\td{e}'$, $e'$ bounds an embedded disk in the
image of
$\td{\Delta}$. Hence, by the previous corollary, $M$ is $S^3$. This proves Part (1).

If there is a face which is a cone, then there are two
possibilities (see Figure \ref{f-faces} (4) and (5)). There is a
tetrahedron $\td{\Delta}$ in the triangulation $\T$ and a face,
$\td{\sigma}$, with vertices $a, b, c$ where the edge
$\overline{ac}$ is identified with the edge $\overline{bc},
a\leftrightarrow b, c\leftrightarrow c$. In one case, Figure
\ref{f-faces} (4), $\sigma$, the image of $\td{\sigma}$,  is an
embedded disk (cone) having boundary the edge $\overline{ab}$. The
other is similar, except the vertex $c$ is also identified with $a
= b$. The cone is not embedded but we can replace a small open
disk in a neighborhood of $c$ in the cone by a disk meeting the
cone only in its boundary, still giving a disk bounded by the edge
$\overline{ab}$. Hence, an edge bounds an embedded disk; by the
previous proposition, $M = S^3$. This completes the proof of Part
(2). \end{proof}

The two vertex, $0$--efficient triangulations of $S^3$ constructed
by Ben Burton feature these anomalies. Later, we show if the
closed, orientable $3$--manifold $M$ has a $0$--efficient
triangulation $\T$, which does not allow certain local reductions
in the number of tetrahedra, then we can make similar conclusions
regarding edges of order two and order three. See Section
\ref{min-triang} on minimal triangulations.

We now show that all closed, orientable, irreducible
$3$--manifolds, except $\rp$, admit a $0$--efficient
triangulation. This is an existence theorem. However, the
techniques provide an algorithm to modify a given triangulation of
a closed, orientable, irreducible $3$--manifold to a
$0$--efficient triangulation of the manifold or along the way
conclude that the given $3$--manifold is $S^3, \rp$ or $L(3,1)$.
The proof of this uses the $3$--sphere recognition algorithm. In
practice, this is not very reasonable; furthermore, one of our
main objectives, leading to $0$--efficient triangulations, is to
get a practical implementation of the $3$--sphere recognition
algorithm. So, we show that given any triangulation of a closed,
orientable $3$--manifold, there is an algorithm to construct a
connected sum decomposition of the $3$--manifold in which each
factor has a $0$--efficient triangulation or is known to be $S^3,
S^2\times S^1, \rp$, or the lens space L(3,1). This decomposition
also has the feature that if the manifold is given via the
triangulation $\T$ and $\T$ is not $0$--efficient, then the number
of tetrahedra needed in the triangulations of all of the factors
in this connected sum decomposition having a $0$--efficient
triangulation is strictly less than the number of tetrahedra in
$\T$. This method is implemented in the computer program ``\small
{REGINA}" by Ben Burton and David Letscher. It can be used to
decide if a given $3$--manifold is the $3$--sphere, the
$3$--sphere recognition algorithm \cite{rub, tho}. We point out
here that the $3$--sphere recognition algorithm really wants a
$0$--efficient triangulation. First, if there are non
vertex-linking, normal $2$--spheres, then standard methods do not
provide an algorithm to find an almost normal $2$--sphere. (In the
original proof of the $3$--sphere recognition algorithm, this was
circumvented by showing that one can construct a maximal
collection of pairwise disjoint, normal $2$--spheres. Then the
closure of a component of the complement of such a collection is
either a  small regular neighborhood of a vertex in the
triangulation, or has multiple (more than one) boundary components
and is a punctured $3$--cell, or has just one boundary component,
is not of the first kind,  and every normal $2$--sphere is
parallel to the boundary component. In the last situation, one has
a cell decomposition which has some of the characteristics of a
$0$--efficient triangulation.) Also, in the presences of non
vertex-linking normal $2$--spheres, the existence of an almost
normal $2$--sphere does not tell one much about the topology of
the given $3$--manifold. Other aspects and implementation of these
algorithms are discussed in \cite{jac-let-rub1, jac-let-rub2}.

\begin{thm} \label{0-eff-exists} A closed, orientable,
irreducible
$3$--manifold distinct from $\rp$ has a
$0$--efficient triangulation.
\end{thm}

\begin{proof}  Suppose $M$ is a closed, orientable irreducible $3$--manifold. Let $\T$
be any triangulation of $M$. If $\T$ is $0$--efficient, there is nothing to prove; so,
we assume $\T$ is not $0$--efficient. Hence, there is a non vertex-linking, normal
$2$--sphere, say $S$, in $M$. Since $M$ is irreducible, $S$ separates and bounds a
$3$--cell in $M$.

Let $S$ be a maximal, non-vertex-linking, normal $2$--sphere in
$M$. Here we are using maximal in the sense that if $S'$ is a non
vertex-linking, normal $2$--sphere and $S'$ bounds a $3$--cell
containing $S$, then $S' = S$. In this situation, such a maximal
non vertex-linking normal $2$--sphere exists by Kneser's Theorem,
Theorem \ref{kneser}. Furthermore, if $S$ is maximal and we denote
the $3$--cell that $S$ bounds by $E$, then by Proposition
\ref{engulf-sphere}, $E$ contains all the vertices of $\T$ or $M$
is the $3$--sphere $S^3$. So, either we have a maximal, non
vertex-linking, normal $2$--sphere $S$ bounding a $3$--cell $E$ in
$M$, which contains all the vertices of $\T$, or $M$ is $S^3$.

 Let
$X$ be the closure of the complement of $E$ in $M$. Then $X$ is
homeomorphic to $M$ with an open $3$--cell removed. We will use
Theorem \ref{crush} to crush the triangulation $\T$ along $S$. By
the comments following Theorem \ref{crush}, if we do this, then we
will have a triangulation $\T^*$, which is an ideal triangulation
of $\open{X}$ and since the index of the ideal vertex is zero,
$\abs{\T^*}$ is homeomorphic to $\hat{X}$, which is homeomorphic
to $M$ and, therefore, $\T^*$ is a triangulation of $M$. Let $\C$
be the induced cell-decomposition of $X$. Then since all the
vertices of $\T$ are in $E$, the cells of $\C$ are of type I, II,
III or IV (see Figure \ref{f-cell-decomp}). To apply Theorem
\ref{crush}, we have to establish the existence of an induced
product region on $X$ and then verify the three conditions in the
hypothesis of Theorem \ref{crush}. We will do this through a
sequence of claims.

\vspace{.15 in} \noindent {\it Claim. $\bbb{P}(\C)$ is a product
$I$--bundle and $\bbb{P}(\C)\ne X$.} \vspace{.1 in}

\noindent {\it Proof of Claim}. Recall that $\bbb{P}(\C) =
\{$edges of $\C$ not in $S\}\cup$ $\{$cells of type III and type
IV in $\C\} \cup$ $\{$all trapezoidal faces of $\C\}$  and is an
$I$--bundle. If $\bbb{P}(\C) = X$, then $X$ is an $I$--bundle with
$2$--sphere boundary; hence, $X$ is a twisted $I$--bundle over
$\rpp$ and $M = \rp$. This contradicts our hypothesis. Similarly,
if a component of $\bbb{P}(\C)$ were not a product $I$--bundle,
then there would be a M\"obius band properly embedded in $X$ (its
boundary in the $2$--sphere $S$). However, if this were the case,
then there would be an $\rpp$ embedded in $M$ and since $M$ is
irreducible, then $M = \rp$. However, again, this contradicts our
hypothesis.

This completes the proof of the {\it Claim}.

\vspace{.15 in} \noindent {\it Claim. There is a trivial induced
product region $\bbb{P}(X)$ for $X$. } \vspace{.15 in}

Before proving this claim, we will establish some notation and
terminology which we will use through the remainder of this work.

 If $K_i\times [0,1]$ is a component of $\bbb{P}(\C)$, then
each component of the complement of $K_{i}^\ve, \ve = 0$ or $1$,
in the $2$--sphere $S$ is simply connected. For each $K_i\times
[0,1]$, let $D_{i}^0$ denote the union of $K_{i}^0$  along with
all components of the complement of $K_{i}^0$ in $S$,  which do
not meet $K_{i}^1$; similarly, let $D_{i}^1$  denote the union of
$K_{i}^1$ along with all components of the complement of $K_{i}^1$
in $S$, which do not meet $K_{i}^0$. For a fixed  $i, 1\le i\le
k$, we have that $D_{i}^0$ and $ D_{i}^1 $ are  disjoint and
simply connected. Furthermore, because of the product $K_i\times
[0,1]$, we have that $D_i^0$ and $D_i^1$ are homeomorphic
subcomplexes in the induced cell structure on $S$. Furthermore,
 $K_{i}^\ve\subseteq D_{i}^\ve, \ve = 0,1$ and we have equality if and only if
$K_i$ is, itself, simply connected. Note, while $D_{i}^0$ is
homeomorphic with $D_{i}^1$, it may be the case that $D_{i}^0$ has
a cell structure induced by the cell structure on $S$, which is
quite different from the induced cell structure on $D_{i}^1$, even
though $K_{i}^0$ and $K_{i}^1$ have isomorphic induced cell
decompositions. Let $N_i = N(K_i\times [0,1]\cup D_{i}^0\cup
D_{i}^1)$ be a small regular neighborhood of the subcomplex
$K_i\times [0,1]\cup D_{i}^0\cup D_{i}^1$ in $X$. The frontier of
$N_i$ consists of one properly embedded annulus component along
with a number of $2$--sphere components (these $2$--sphere
components exists if and only if $K_i$ is not simply connected).
In our current situation ($M$ is irreducible), each of these
$2$--spheres separates, bounds a $3$--cell and is a barrier
surface in the closure of its complement in $X$ not meeting
$K_i\times [0,1]\cup D_{i}^0\cup D_{i}^1$. While we know each
$2$--sphere bounds a $3$--cell, we want to show that such a
$2$--sphere bounds a $3$--cell not meeting $K_i\times [0,1]\cup
S$. We will denote the closure of a component of the complement of
$N_i$ in $X$, which has a $2$--sphere boundary component, by
$N_{i,j}$ and its boundary $2$--sphere (also in the boundary of
$N_i$) by $S_{i,j}$. We call $N_{i,j}$ a {\it plug for $K_i\times
[0,1]$}. See Figure \ref{f-product-fill}.

\begin{figure}[htbp]

            \psfrag{X}{$X$}
            \psfrag{S}{$S$}
            \psfrag{s1}{$S_{i,1}$}
            \psfrag{s2}{$S_{i,2}$}
            \psfrag{L}{$K_i\times [0,1]$}
            \psfrag{k1}{$K_i^1$}
            \psfrag{k0}{$K_i^0$}
            \psfrag{N1}{$N_{i,1}$}
            \psfrag{N2}{$N_{i,2}$}
            \psfrag{D1}{$D_i^1$}
            \psfrag{D0}{$D_i^0$}
            \psfrag{P}{$K_i\times [0,1]\cup D_i^0\cup D_i^1$}
            \psfrag{w}{\begin{tabular}{c}
            new normal\\
       $2$--spheres\\
            \end{tabular}}
        \vspace{0 in}
        \begin{center}
        \epsfxsize=5 in
        \epsfbox{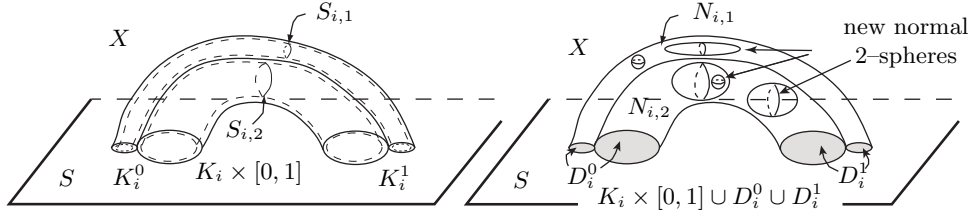}
        \caption{Plugs for product regions in $X$ and new normal
$2$--spheres, which  appear in shrinking $S_{i,j}$.}
        \label{f-product-fill}
        \end{center}

\end{figure}

Using this terminology, we will establish in this claim that a
plug for $K_i\times [0,1]$ is a $3$--cell. We remark that there
are a number of alternate ways to approach our problem at this
juncture. We have chosen the following as the most straight
forward and natural for this proof and for similar situations we
encounter later in this work.

\vspace{.15 in}\noindent{\it Proof of Claim}: Chose some order for
the components of $\bbb{P}(\C)$, say $K_1\times [0,1],\ldots,
K_n\times [0,1]$, where $n$ is the number of components of
$\bbb{P}(\C)$. If $K_1$ is simply connected, then let $D_1 = K_1$.
If $K_1$ is not simply connected, let $D_{1}^0$ and $D_{1}^1$ be
as defined above. Then, using the above notation, we let $N_1$ be
a small regular neighborhood of $K_1\times [0,1]\cup D_{1}^0\cup
D_{1}^1$. Associated with $N_1$, we have some set of plugs
$N_{1,j}$ and boundary $2$--spheres $S_{1,j}$, each of the latter
being a barrier surface in the component of its complement in
$N_{1,j}$ not meeting $K_1\times [0,1]\cup D_{1}^0\cup D_{1}^1$.
Shrink each $S_{1,j}$ in $N_{1,j}$. Each $S_{1,j}$ shrinks to a
collection of normal $2$--spheres and possibly some $0$--weight
$2$--spheres contained entirely in the interior of tetrahedra.
Again, see Figure \ref{f-product-fill}. We can fill each of the
$0$--weight $2$--spheres in with a $3$--cell which does not meet
$K_1\times [0,1]\cup S$. Now, suppose in shrinking some $S_{1,j}$
we have a normal $2$--sphere. Since $M$ is irreducible, any such
$2$--sphere bounds a $3$--cell in $M$. However, by the maximality
of $S$, such a $3$--cell does not meet $K_1\times [0,1]\cup S$. It
follows that if  such a normal $2$--sphere occurs, then it bounds
a $3$--cell in $X$ missing $K_1\times [0,1]\cup S$. Thus each plug
for $N_1$ is a $3$--cell.

Since the frontier of each of the regions complementary to
$K_1\times [0,1]$ is a union of trapezoids, there is a simply
connected planar complex $D_1$ and an embedding of $D_1\times
[0,1]$ into $X$ with $D_1\times\ve = D_{1}^\ve, \ve = 0,1$,
$K_1\times [0,1]\subseteq D_1\times [0,1]$ and the frontier of
$D_1\times [0,1]$ is contained in the frontier of $K_1\times
[0,1]$. So, we have replaced the product $K_1\times [0,1]$ with a
product $D_1\times [0,1]$ where we have that $D_1$ is a simply
connected planar complex.

Now, suppose for the components $K_1\times [0,1],\ldots, K_k\times
[0,1]$ of $\bbb{P}(\C)$, we have simply connected planar complexes
$D_1,\ldots,D_{k'}, k'\le k$ and embeddings $D_j\times [0,1], 1\le
j\le k'$ into $X$ so that $D_j\times\ve = D_{j}^\ve, \ve = 0,1$.
Furthermore, suppose for $j\ne j'$, either $(D_j\times [0,1])\cap
(D_{j'}\times [0,1]) = \emptyset$ or $D_j\times [0,1]\subset
D_{j'}\times [0,1]$ or vice versa; $\bigcup_1^k (K_i\times
[0,1])\subset \bigcup_1^{k'}(D_j\times [0,1])$ and the frontier of
$\bigcup_1^{k'}(D_j\times [0,1])$ is contained in the frontier of
$\bigcup_1^k (K_i\times [0,1])$. So, we have replaced a number of
the $K_i\times [0,1]$ with trivial products. See Figure
\ref{f-product-1}.

If $k< n$, consider the component $K_{k+1}\times [0,1]$ of
$\bbb{P}(\C)$. If  $K_{k+1}\times [0,1] \subset
\bigcup_1^{k'}(D_j\times [0,1])$, then $\bigcup_1^{k+1} (K_i\times
[0,1])\subset \bigcup_1^{k'}(D_j\times [0,1])$ and the frontier of
$\bigcup_1^{k'}(D_j\times [0,1])$ is contained in the frontier of
$\bigcup_1^{k+1} (K_i\times [0,1])$ and there is nothing to do.
So, suppose $K_{k+1}\times [0,1]\not\subset
\bigcup_1^{k'}(D_j\times [0,1])$. If $K_{k+1}$ is simply
connected, we set $D_{k'+1} = K_{k+1}$; $D_{k'+1}\times [0,1]$ is
disjoint from each $D_j\times [0,1], 1\le j\le k'$. So, suppose
$K_{k+1}$ is not simply connected. We let $D_{k+1}^\ve, \ve = 0,1$
be defined as above and let $N_{k+1} = N(K_{k+1}\times [0,1]\cup
D_{k+1}^0\cup D_{k+1}^1)$ denote a small regular neighborhood of
$K_{k+1}\times [0,1]\cup D_{k+1}^0\cup D_{k+1}^1$. Then just as in
the above for the case of $K_1\times [0,1]$ and $N_1$, we can show
that each plug for $N_{k+1}$ is a $3$--cell.

Thus there is a simply connected planar complex $D_{k'+1}$ and an
embedding of $D_{k'+1}\times [0,1]$ into $X$ so that
$D_{k'+1}\times\ve = D_{k+1}^\ve, \ve = 0,1$. Furthermore, for
$j\le k'$, either $(D_j\times [0,1])\cap (D_{k'+1}\times [0,1] )=
\emptyset$ or $D_j\times [0,1]\subset D_{k'+1}\times [0,1]$; and
finally,  $\bigcup_1^{k+1} (K_i\times [0,1])\subset
\bigcup_1^{k'+1}(D_j\times [0,1])$ and the frontier of
$\bigcup_1^{k'+1}(D_j\times [0,1])$ is contained in the frontier
of $\bigcup_1^{k+1} (K_i\times [0,1])$.  Hence, we can enlarge our
collection  $D_1,\ldots,D_{k'}, k'\le k$ to a collection including
$D_{k'+1}$, which incorporates the product $K_{k+1}$. See Figure
\ref{f-product-1}.

Let $\bbb{P}(X)$ be the union of the components of  $\bigcup_i
(D_i\times [0,1])$. As we observed above, it is possible that some
$D_j'\times [0,1]$ may actually be embedded in a $D_{j}\times
[0,1]$, $j'\neq j$. For example, in Figure \ref{f-product-1}, we
have $D_1\times [0,1]\subset D_2\times [0,1]$ and $D_3\times
[0,1]\subset D_4\times [0,1]$. We also see that $K_4\times
[0,1]\subset D_2\times [0,1]$ and therefore, is dropped.  In this
example, $\bbb{P}(X)$ has three components, $D_2\times [0,1],
D_4\times [0,1]$ and $D_5\times [0,1]$. Also, this example shows
that there is not a unique product region; had we ordered the
components $K_1\times [0,1],\ldots, K_n\times [0,1]$ of
$\bbb{P}(\C)$ differently, say with $K_4\times [0,1]$ first, then
all the other components of $\bbb{P}(\C)$ would have been dropped
and we would have only one component for $\bbb{P}(X)$. In those
cases where $M$ is irreducible, if something like this happened,
then it would follow that $M$ is $S^3$.

\begin{figure}[htbp]

            \psfrag{S}{$S$}
            \psfrag{1}{\begin{tabular}{c}
           {\scriptsize$K_1\times [0,1]$}\\
       {\scriptsize$\subset D_1\times
            [0,1]$}\\
            \end{tabular}}
            \psfrag{2}{\scriptsize{$K_2\times [0,1]\subset D_2\times
            [0,1]$}}
            \psfrag{3}{\scriptsize{$K_3\times [0,1] = D_3\times
            [0,1]$}}
            \psfrag{4}{\scriptsize{$K_4\times [0,1]$}}
            \psfrag{5}{\scriptsize{$K_5\times [0,1]\subset D_4\times
            [0,1]$}}
            \psfrag{6}{\scriptsize{$K_6\times [0,1]= D_5\times [0,1]$}}
            \psfrag{a1}{\tiny$D_1^0$}
            \psfrag{b1}{\tiny$D_1^1$}
            \psfrag{a2}{\scriptsize$D_2^0$}
            \psfrag{b2}{\scriptsize$D_2^1$}
            \psfrag{a4}{\scriptsize$D_4^0$}
            \psfrag{b4}{\scriptsize$D_4^1$}
        \vspace{0 in}
        \begin{center}
        \epsfxsize=3 in
        \epsfbox{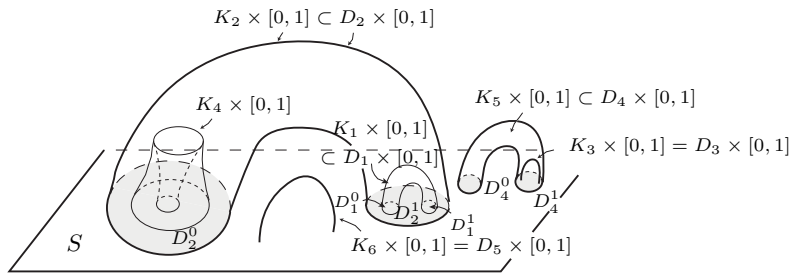}
        \caption{Building a trivial product region, $\bbb{P}(X)$, from the combinatorial
products $K_i\times [0,1]$.}
        \label{f-product-1}
        \end{center}

\end{figure}

This completes the proof of the {\it Claim}.

\vspace{.15 in} \noindent {\it Claim. Either there is no cycle of
truncated prisms in $X$, which is not in $\bbb{P}(X)$, or $M$ is
the manifold $S^3$ or the manifold $L(3,1)$.} \vspace{.1 in}

\noindent {\it Proof of Claim}. There are two types of cycles of
truncated prisms: one is a cycle about an edge $e$ of $\T$ (see
Figure \ref{f-cycle}(A)) and the other cycles about more than one
edge of $\T$ (see Figure \ref{f-cycle}(B)).

\begin{figure}[htbp]

            \psfrag{d}{$e'$}
            \psfrag{f}{$e''$}
            \psfrag{e}{$e$}
            \psfrag{A}{\begin{tabular}{c}
            (A) Cycle about\\
    a single edge\\
            \end{tabular}}
             \psfrag{B}{\begin{tabular}{c}
            (B) Cycle about\\
      more than one edge\\
            \end{tabular}}
        \vspace{0 in}
        \begin{center}
        \epsfxsize=4 in
        \epsfbox{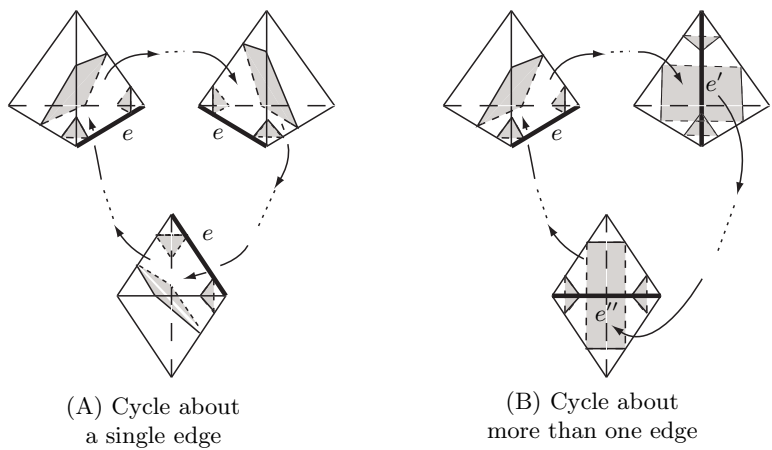}
        \caption{Cycles of truncated prisms.}
        \label{f-cycle}
        \end{center}

\end{figure}

If there is a complete cycle about an edge $e$ as in Figure
\ref{f-cycle}(A) (i.e., the $2$--sphere $S$ contains a thin tube
of elementary quads about the edge $e$), then there is a properly
embedded disk $D$ in $X$ meeting $e$ in precisely one point and
meeting $S$ in $\bdy D$. A surgery on $S$ at $D$ gives two {\it
normal} $2$--spheres $S_{0}$ and $S_1$, neither of which is
vertex-linking and together with $S$ bound a punctured $3$--cell.
See Figure \ref{f-surgery-cycle}. Since $M$ is irreducible, each
of these $2$--spheres bounds a $3$--cell. If either bounds a
$3$--cell containing $S$ (hence, the $3$--cell $E$ and all the
vertices of $\T$), then we have a contradiction to the maximality
of the $2$--sphere $S$. So, the only possibility is that they
bound $3$--cells not containing $E$. Hence, $M$ must be $S^3$. It
follows that having chosen $S$ a  maximal $2$--sphere, if there is
a complete cycle of truncated prisms in the induced cell structure
on $X$ about a single edge, then we have that the manifold $M$ is
homeomorphic with $S^3$.

\begin{figure}[htbp]

            \psfrag{D}{$D$}
            \psfrag{s}{surgery}
            \psfrag{1}{$S_1$}
            \psfrag{e}{$e$}
            \psfrag{0}{$S_0$}
        \vspace{0 in}
        \begin{center}
        \epsfxsize=3.8 in
        \epsfbox{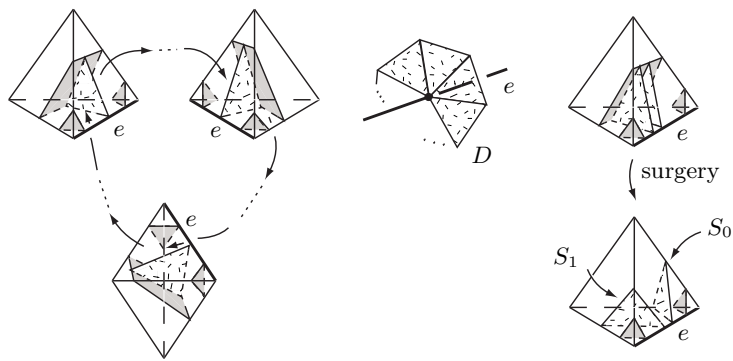}
        \caption{Cycle about a single edge gives a new normal $2$--sphere.}
        \label{f-surgery-cycle}
        \end{center}

\end{figure}

If there is a complete cycle about more than one edge (see Figure
\ref{f-cycle}(B)), then the collection of cells of type II
(truncated prisms), form a solid torus with, possibly, some self
identifications in its boundary (possibly, some of the trapezoidal
faces in the boundary are identified); furthermore, each hexagonal
face of a truncated prism in the cycle of truncated prisms is a
meridional disk for the solid torus. Because of the possible
singularities, we distinguish between the cycle of truncated
prisms, which we denote $\hat{\tau}$, and the cycle of truncated
prisms minus the bands of trapezoids, which we denote by $\tau$.
We have that $\tau\cap S$ is either three open annuli, each
meeting a meridional disk of $\tau$ precisely once, or a single
open annulus, meeting a meridional disk of $\tau$ three times. See
Figure \ref{f-cycle-torus} where we also show the trapezoidal
annuli slightly shrunken into the torus $\tau$.

\begin{figure}[htbp]

            \psfrag{1}{{\footnotesize $A_1$}}
            \psfrag{2}{{\footnotesize $A_2$}}
            \psfrag{3}{{\footnotesize $A_3$}}
            \psfrag{a}{{\footnotesize $A_1'$}}
            \psfrag{b}{{\footnotesize $A_2'$}}
            \psfrag{c}{{\footnotesize $A_3'$}}
            \psfrag{t}{\footnotesize $t$}
                \psfrag{x}{$\hat{\tau}$}
                    \psfrag{S}{{\footnotesize in $S$}}
                    \psfrag{s}{{\footnotesize $t = $}\hspace{.1 in}{\footnotesize trapezoid}}
            \psfrag{A}{\begin{tabular}{c}
           {\footnotesize $\hat{\tau}$, a cycle of}\\
   {\footnotesize truncated prisms}\\
            \end{tabular}}
             \psfrag{C}{\begin{tabular}{c}
            {\footnotesize $\overline{\tau}$, a cycle of}\\
      {\footnotesize (slightly shrunken)}\\
      {\footnotesize truncated prisms}\\
            \end{tabular}}
             \psfrag{D}{\begin{tabular}{c}
            {\footnotesize solid torus of}\\
      {\footnotesize (hexagonal cubes)}\\
      {\footnotesize truncated prisms}\\
            \end{tabular}}
             \psfrag{E}{\begin{tabular}{c}
            {\footnotesize one annulus}\\
      {\footnotesize in $S\cap\overline{\tau}$ (or $\tau$)}\\
            \end{tabular}}
             \psfrag{F}{\begin{tabular}{c}
            {\footnotesize three annuli}\\
      {\footnotesize in $S\cap\overline{\tau}$ (or $\tau$)}\\
            \end{tabular}}
        \vspace{0 in}
        \begin{center}
        \epsfxsize=3.5 in
        \epsfbox{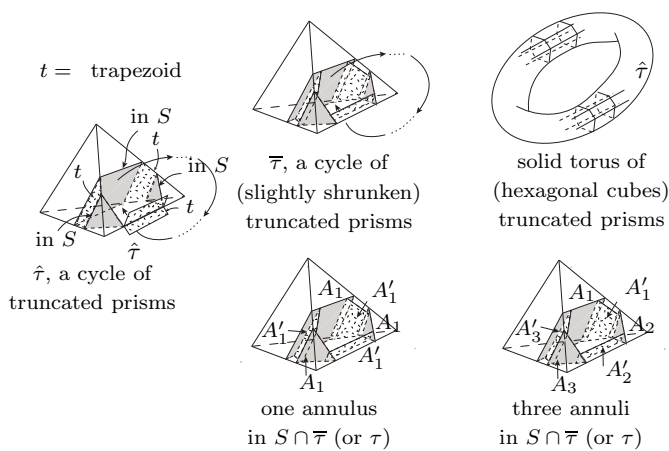}
        \caption{Cycle of truncated prisms determine a solid torus.}
        \label{f-cycle-torus}
        \end{center}

\end{figure}

First, we consider the case we have three annuli in $\tau\cap S$,
say $A_1, A_2$ and $A_3$, each meeting a hexagonal face of a
truncated prism in $\hat{\tau}$ precisely once. The frontier of
each $A_i, i = 1,2,3$ in $S$ is in the collection of trapezoids in
the faces of the cycle of truncated prisms $\hat{\tau}$; thus in
the induced $I$--bundle region $\bbb{P}(\C)$.  Furthermore, there
are two components of the frontier of each $A_i$ and each
component of the frontier of $A_i$ separates $S$. Consider $A_1$
and denote the two components of the frontier of $A_1$ by $a_1$
and $a_1'$. Then both $a_1$ and $a_1'$ are in trapezoids in
$\hat{\tau}$ and so in $\bbb{P}(\C)$. By our above construction of
the induced product region for $X$, $\bbb{P}(X)$, it follows that
$a_1$ and $a_1'$ are each in a simply connected region of $S$
common to $\bbb{P}(X)$. Hence, we either have $A_1$, and therefore
$\hat{\tau}$ in $\bbb{P}(X)$ or $\bbb{P}(X)$ meets $S$ in the
complement of $A_1$. But the latter is impossible, since $A_2$ and
$A_3$ are in regions of $S$ complementary to $A_1$. So, the only
possibility is that $A_1$ is in $\bbb{P}(X)$; that is, such a
cycle of truncated prisms, $\hat{\tau}$, is in the induced product
region for $X$, $\bbb{P}(X)$. See Figure \ref{f-annuli}.

\begin{figure}[htbp]

            \psfrag{S}{\large$S$}
            \psfrag{a}{\small$a_1$}
            \psfrag{b}{\small$a_1'$}
            \psfrag{1}{\small$A_1$}
            \psfrag{2}{\small$A_2$}
            \psfrag{3}{\small$A_3$}
        \vspace{0 in}
        \begin{center}
        \epsfxsize=3 in
        \epsfbox{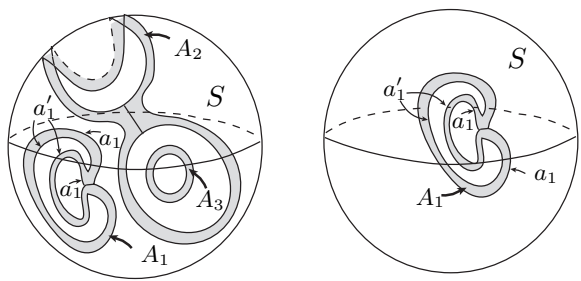}
        \caption{The (possibly singular) torus $\hat{\tau}$ meets $S$ in
(possibly singular) annuli.}
        \label{f-annuli}
        \end{center}

\end{figure}

While the previous observations take care of this case, it can be
shown that if $S$ is maximal, $M$ is irreducible and there is a
cycle of truncated prisms about more than one edge and there are
three annuli in $\tau\cap S$, then $M$ is $S^3$.

If we have just one annulus, say $A_1$, in $\tau\cap S$, then
$A_1$ meets each hexagonal face in the chain of truncated prisms,
$\hat{\tau}$, three times. We can not make a conclusion similar to
that above when there were three annuli, since it is possible that
the components of $S$ complementary to $A_1$ are in the induced
product region on $X$, $\bbb{P}(X)$; in fact, it is necessary.
However,  the core of the annulus $A$ bounds a disk in $S$; so, it
follows that $M$ is a connected sum with the lens space $L(3,1)$.
But $M$ is irreducible so is, itself, $L(3,1)$.  This completes
the proof of the Claim.

Thus we have shown that the conditions in the hypothesis of
Theorem \ref{crush} are satisfied; furthermore, we may assume
there are no cycles of truncated prisms in $X$, which are not in
$\bbb{P}(X)$, since these lead to the conclusion that $M$ is
either $S^3$ or $L(3,1)$. So, from the truncated tetrahedra in
$\C$, which are not in $\bbb{P}(X)$ (and there must be some as
$\bbb{P}(X) \neq X$ and there are no complete cycles of truncated
prisms, which are not in $\bbb{P}(X)$), we get an ideal
triangulation $\T^*$ of $\open{X}$. However, since $S$ is a
$2$--sphere the ideal vertex has index zero and so, $\abs{\T^*}$
is homeomorphic to $M$ and $\T^*$ is a triangulation of $M$. Since
$S$ was chosen to be non vertex-linking, there is at least one
truncated prism in $X$ and so there are strictly fewer truncated
tetrahedra in $\C$ than there are in $\T$. It follows that the
triangulation $\T^*$ of $M$ has fewer tetrahedra than the
triangulation $\T$.

 Hence, by iterating the process if necessary, we must
terminate in a desired  $0$--efficient triangulation of $M$ or we
have that $M$ is homeomorphic with $S^3$ or the lens space
$L(3,1)$. Each of these latter two manifolds admits a
$0$--efficient triangulation. See Figure \ref{f-tetra}(4) and (5)
and Figure \ref{f-two-L3_1}(A).

This completes the proof of the Theorem. \end{proof}

Now, we have a series of results which follow from variations in the proof of Theorem
\ref{0-eff-exists}. Other variations and generalizations
 appear in \cite{jac-let-rub1}.

\begin{thm} Suppose $M$ is a closed, orientable, irreducible $3$--manifold. Then $M$
admits a one-vertex triangulation.
\end{thm}

\begin{figure}[htbp]

            \psfrag{a}{$a$}
            \psfrag{b}{$b$}
            \psfrag{c}{$c$}
            \psfrag{d}{$d$}
        \vspace{0 in}
        \begin{center}
        \epsfxsize=4 in
        \epsfbox{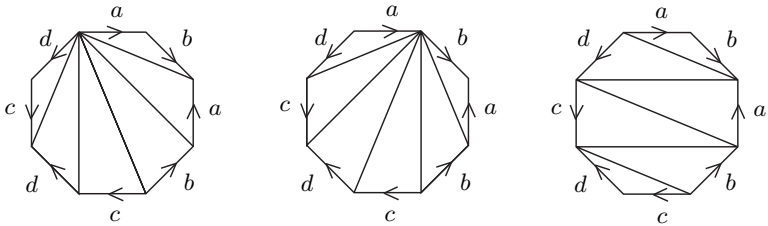}
        \caption{Distinct one-vertex triangulations of a genus $2$
surface.}
        \label{f-one-vertex-surface}
        \end{center}

\end{figure}

The standard presentation of closed $2$--manifolds as an
identification space of a planar $n$--gon, leads to easy
one-vertex triangulations of closed $2$--manifolds distinct from
$S^2$ and $\rpp$. See Figure \ref{f-one-vertex-surface} for three
distinct one-vertex triangulations of a genus two surface.

When we first started this study of triangulations of
$3$--manifolds, it was not exactly clear how to get a one-vertex
triangulation of an arbitrary $3$--manifold (we now know several
ways, see \cite{jac-let-rub1}). Theorem \ref{0-eff-exists},
Proposition \ref{p-0-eff} and one-vertex triangulations of $S^3,
\rp$ and $L(3,1)$ provide an argument showing that any closed,
orientable, irreducible $3$--manifold admits a one-vertex
triangulation. However, the methods in the proof of Theorem
\ref{0-eff-exists} provide a direct proof of existence of
one-vertex triangulations of a closed, orientable, irreducible
$3$--manifold. Furthermore, these methods enable us to do what
seems obvious; namely, crush a maximal tree in the $1$--skeleton
to a point. They do not require the iteration necessary in the
proof of Theorem \ref{0-eff-exists}. Note, these methods do not
seem to adapt to show that any closed $3$--manifold admits a
one-vertex triangulation; however, it can be shown, using other
methods, that this is true (see \cite{jac-let-rub1}).

\begin{proof}
 As we have said, the intuitive thing to do is to ``crush" a maximal tree in the
$1$--skeleton of a triangulation to a point. While this seems
obvious; it is not so obvious how to organize a proof without the
methods used in the proof of Theorem \ref{0-eff-exists}.  To  this
end, suppose we have a triangulation $\T$ of $M$. If $\T$ has just
one vertex, then there is nothing to prove.  So assume $\T$ has
more than one vertex. Let $S'$ be a vertex-linking normal
$2$--sphere. Then by Proposition \ref{engulf-sphere}, there is a
normal $2$--sphere $S$ bounding a ball $E$, which contains all the
vertices of $\T$. As in the proof of Theorem \ref{0-eff-exists},
either there is a maximal normal $2$--sphere bounding a $3$--cell
and containing all the vertices of $\T$ or we have that $M$ is
$S^3$. We now use $S$ to denote such a maximal $2$--sphere. We
wish to crush the triangulation $\T$ along $S$.

As in the proof of Theorem \ref{0-eff-exists}, let $X$ denote the
closure of the complement of $E$ and let $\C$ be the
cell-structure induced on $X$ by $\T$. Now, if  $\bbb{P}(\C) =
\{$edges of $\C$ not in $S\}\cup$ $\{$cells of type III and type
IV in $\C\} \cup$ $\{$all trapezoidal faces of $\C\}$, then each
component of $\bbb{P}(\C)$ is an $I$-bundle. If $\bbb{P}(\C)$ is
not a product $I$--bundle, then there is a M\"obius band in
$\bbb{P}(\C)$ with its boundary in $S$. It follows, in this case,
there is an $\rpp$ embedded in $M$. Similarly, if $\bbb{P}(\C) =
X$, then since $S$ is a $2$--sphere, we have the $I$--bundle
$\bbb{P}(\C)$ is a twisted $I$--bundle with boundary a
$2$--sphere; so, is a twisted $I$--bundle over $\rpp$. It follows
in both situations that $M = \rp$, which admits a one-vertex
triangulation. The proof now follows the proof of Theorem
\ref{0-eff-exists} and either we can crush the triangulation $\T$
along $S$ or we have that $M$ is one of $S^3$ or $L(3,1)$. Both of
which admit one-vertex triangulations.

Let $\T^*$ denote the ideal triangulation of $\open{X}$ we achieve
by  crushing $\T$ along $S$. Then, just as above, $\abs{\T^*}$ is
homeomorphic to $M$. $\T^*$ is a one-vertex triangulation of $M$.
In this case, we are done; we do not need to iterate, as we may
have needed to do in the proof of Theorem
\ref{0-eff-exists}.\end{proof}

\begin{figure}[htbp]

            \psfrag{A}{(A) $\rpp\subset \rp$}
            \psfrag{B}{(B) $S^2\subset S^2\times S^1$}
        \vspace{0 in}
        \begin{center}
        \epsfxsize=4 in
        \epsfbox{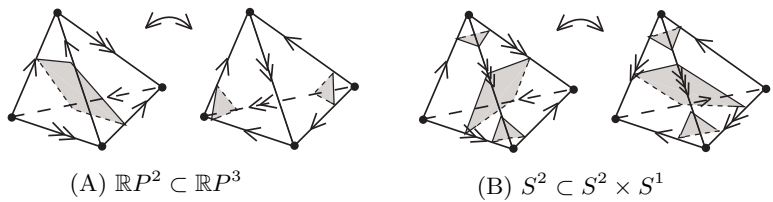}
        \caption{Minimal (one-vertex) triangulations of $\rp$ and
$S^2\times S^1$.}
        \label{f-onevertex-RP3-S2xS1}
        \end{center}

\end{figure}

One of the motivating features of $0$--efficient triangulations is
a more effective implementation of the $3$--sphere recognition
algorithm. Our next result provides the environment to apply the
$3$--sphere recognition algorithm. In particular, we show that
given a triangulation of an arbitrary closed, orientable
$3$--manifold $M$, we can produce a connected sum decomposition of
$M$ into factors, each of which either has a $0$--efficient
triangulation or is known to be one of the manifolds $S^3,
S^2\times S^1, \rp$ or the lens space $L(3,1)$. We single out
$\rp$ and $S^2\times S^1$ because they get split off in our
construction and do not admit $0$--efficient triangulations. We
include $S^3$ and $L(3,1)$, even though both admit $0$--efficient
triangulations, because from our construction and like $S^2\times
S^1$ and $\rp$, we know precisely the homeomorphism type of some
of the factors without resorting to a recognition algorithm.
Finally, the total number of tetrahedra needed in our
triangulations of the components in this connected sum
decomposition, which have $0$--efficient triangulations, is
strictly smaller than the number of tetrahedra in $\T$. It may be
that the total number of tetrahedra needed to give our
triangulations of all the factors, including those which are $S^3,
S^2\times S^1, \rp$ and $L(3,1)$, is, in general, a smaller number
than the number of tetrahedra for $\T$. This seems likely if we
use a one-tetrahedron triangulation of $S^3$ (or toss $S^3$ out of
the connected sum) or the two tetrahedra, one-vertex
triangulations for $\rp, S^2\times S^1$ and $L(3,1)$. See Figure
\ref {f-onevertex-RP3-S2xS1}, which gives the unique minimal
triangulation of $S^2\times S^1$ and the two tetrahedra,
one-vertex triangulation for
 $\rp$. See Figure \ref {f-two-L3_1}(A), which gives a two-tetrahedron,
 one-vertex, $0$--efficient
 triangulation of $L(3,1)$ and Figure \ref{f-tetra}(4) and (5) for one tetrahedra,
 $0$--efficient triangulations of $S^3$.
 However, this is not clear and seems more a curiosity than useful.

We do not construct normal $2$--spheres in $M$ giving our
decomposition; however, it is possible to do this from our
methods. It also is possible to construct an irreducible
decomposition. However, we must wait and apply the $3$--sphere
recognition algorithm to get a prime decomposition, because there
may be factors which are $S^3$ but it is unknown at this stage.
There may also be factors which have a $0$--efficient
triangulation and are $L(3,1)$, whose homeomorphism type is also
not known at this time. We give an alternate method to construct
an irreducible decomposition of a given $3$--manifold. This is by
construction of a maximal pairwise disjoint collection of distinct
normal $2$--spheres. The methods we use provide algorithms and the
different techniques may lead to one being preferred to the other
when we consider issues of complexity. Our proof of the
decomposition theorem, Theorem \ref{reduce} below, provides the
steps for an algorithm, which has been coded  by David Letscher
and Ben Burton, to construct $0$--efficient triangulations and
implement the $3$--sphere recognition algorithm. It seems to be
the preferred algorithm for arriving at an irreducible
decomposition of a given $3$--manifold. The philosophy behind this
preference is that in practice the crushing of a triangulations
along normal (non vertex-linking) $2$--sphere greatly reduces the
number of tetrahedra; so, one should crush whenever one can. A
construction of a maximal collection of normal $2$--spheres is a
result obtained by the authors in 1988 and is used by the second
author in the proof of the $3$--sphere recognition algorithm.
Quite different algorithms, which produce a family of normal
$2$--spheres giving an irreducible or a prime decomposition of a
given $3$--manifold, appear in \cite{jac-tol} and \cite{jac-ree}.

First, we give some results necessary for our proof. the next
Proposition also appears in \cite{jac-ree}.

\begin{prop}\label{0-decide} Given a triangulation $\T$ of a closed, orientable
$3$--manifold $M$, there is an algorithm to decide if $\T$ is
$0$--efficient; furthermore, the algorithm will construct a non
vertex-linking $2$--sphere, if one exists.
\end{prop}

\begin{proof} If there is a  non vertex-linking, normal $2$--sphere, we show
there must be one whose projective class is a vertex of the
projective solution space; i.e., a vertex solution. Since a normal
$2$--sphere is not vertex-linking if and only if it meets some
tetrahedron in a quadrilateral, it can be decided if a normal
$2$--sphere is non vertex-linking. To this end suppose $\Sigma$ is
a non vertex-linking, normal $2$-sphere, $\mathcal{C}(\Sigma)$ is
the carrier of $\Sigma$ and $\Sigma$ has been chosen so that among
all non vertex-linking, normal $2$--spheres in $M$, the dimension
of $\mathcal{C}(\Sigma)$ is a minimum. If $\mathcal{C}(\Sigma)$ is
not a vertex of the projective solution space, then there are
normal surfaces $X$ and $Y$ carried by proper faces of
$\mathcal{C}(\Sigma)$ and nonnegative integers $k, n,$ and $m$
with $k\Sigma = nX + mY$. Since, $\chi(\Sigma) > 0$, we must have
either $\chi(X) > 0$ or $\chi(Y) > 0$ and so a component of one
with positive Euler characteristic, say a component $X'$ of $X$.
The carrier of $X'$ is also a proper face of
$\mathcal{C}(\Sigma)$. $X'$ can not be a projective plane; for
then, its double would be a non vertex-linking, normal $2$--sphere
carried by a proper face of $\mathcal{C}(\Sigma)$ and so, as a
proper face of $\mathcal{C}(\Sigma)$, has dimension less than that
of $\mathcal{C}(\Sigma)$. The only other possibility is that $X'$
is a vertex-linking $2$--sphere; however, then a component of
$k\Sigma$ would be vertex-linking and so $\Sigma$ would be
vertex-linking. This is a contradiction to our choice of $\Sigma$.
Thus $\Sigma$ has to be carried by a vertex of the projective
solution space. The vertices of  the projective solution space can
be constructed and as we remarked, we can decide if a normal
surface is a vertex-linking $2$--sphere.
\end{proof}

\begin{lem}\label{oneboundary} Suppose $M$ is a closed, orientable $3$--manifold, $\T$ is a
triangulation of $M$ and $\mathcal{S}$ is a collection of pairwise
disjoint, normal $2$--spheres embedded in $M$.  If $X$ is the
closure of a component of the complement of the $2$--spheres in
$\mathcal{S}$ and $\bdy X$ is not connected, then there is a
punctured $3$--sphere, $P_X$,  embedded in $X$, each component of
$\bdy P_X$ is a normal $2$--sphere, $\bdy X\subset \bdy P_X$ and
no component of the frontier of $P_X$ in $X$ is normally isotopic
to a component of $\bdy X$.
\end{lem}

\begin{proof} We have that
$X$ is the closure of a component of the complement of the
collection of normal $2$--spheres $\mathcal{S}$ and $X$ has more
than one boundary component. Hence, we can find a pairwise
disjoint collection, $\Lambda_X$, of arcs in the $1$--skeleton of
the induced cell-decomposition on  $X$, each such arc has its end
points in  distinct components of $\bdy X$ and the complex
$\Lambda_X \cup \bdy X$ is connected and simply connected. A small
regular neighborhood of the  complex $\Lambda_X \cup \bdy X$ is a
punctured $3$--sphere and its frontier, which is a single
$2$--sphere, is a barrier surface in the component of its
complement not containing the complex $\Lambda_X \cup \bdy X$.
Thus we can shrink this frontier $2$--sphere. We obtain a
punctured $3$--sphere, containing $\Lambda_X \cup \bdy X$ and
whose boundary consists of all the $2$--spheres in $\bdy X$, along
with (possibly) some normal $2$--spheres in the interior of $X$
and some $0$--weight $2$--spheres contained entirely in the
interior of cells in $X$. The $2$--spheres in the interior of
cells in $X$ can be filled in with $3$--cells in $X$; we do this
and forget them. None of the normal $2$--spheres in the interior
of $X$ are normally isotopic to a component of $\bdy X$, since
there is at least one edge of the collection of edges $\Lambda$
meeting each boundary component and the complex $\Lambda_X\cup
\bdy X$ is a barrier surface. It is possible, of course, that we
have no frontier of $X$  and so $X$ is itself a punctured
$3$--sphere. This completes the proof.\end{proof}

Note Lemma \ref{oneboundary} can be modified to be applicable to
the situation where we have a manifold with boundary and some of
the $2$--spheres in $\mathcal{S}$ are in the boundary of the
manifold.

\begin{thm}\label{reduce} Given a closed, orientable $3$--manifold $M$, there is an algorithm
to construct a finite family of $3$--manifolds, $M_1,\ldots,M_n$,
so that $M = M_1 \#\ldots\# M_n$, where $M_i, i = 1,\ldots,n$,
either has a $0$--efficient triangulation or can be shown to be
homeomorphic with one of $S^3, S^2\times S^1, \rp$, or the lens
space, $L(3,1)$.
\end{thm}

\begin{proof} We shall first outline our approach. If $\T$ is not
$0$--efficient, then there is a non vertex-linking normal
$2$--sphere and by Proposition \ref{0-decide}, we have an
algorithm to construct one. Our next step is to construct a
punctured $3$--sphere, say $P$, which contains the vertices of
$\T$ and has non vertex-linking normal $2$--spheres in its
boundary or we can conclude that $M = S^3$. So, suppose we have
such a punctured $3$--sphere $P$. If $X$ is the closure of a
component of the complement of $P$, then $\hat{X}$, the
$3$--manifold $X$ with its $2$--sphere boundary components
``capped-off" with $3$--cells, is a factor in a connected sum
decomposition of $M$; however,  if the boundary of $X$ is not
connected, then we must take care not to miss copies of $S^2\times
S^1$, which will also be factors in such a connected sum
decomposition. If we do have some component of the complement of
$P$ whose closure does not have connected boundary, then we can
use Lemma \ref{oneboundary} to enlarge $P$ and obtain a connected
sum of a punctured $3$--sphere and, possibly, some copies of
$S^2\times S^1$ so that each component of the boundary of our new
connected sum is a non vertex-linking normal $2$--sphere and the
closure of each component of the complement of this new connected
sum has connected boundary. We will continue to use $P$ to denote
this new punctured $3$--sphere connected sum with (possibly) some
copies of $S^2\times S^1$. Next, we consider the closure of a
component of the complement of our new $P$, say $X$. The
$3$--manifold $X$ has connected boundary and has a nice
cell-decomposition, say $\C$, induced by $\T$ (all the vertices of
$\T$ are in $P$). We wish to crush the triangulation along the
$2$--sphere in $\bdy X$. To do this we must establish the three
conditions in the hypothesis of Theorem \ref{crush}. The first
condition is that $\bbb{P}(\C)\ne X$ and each component of
$\bbb{P}(\C)$ is a product $I$--bundle. To establish this
condition for components of the complement of $P$, we may need to
enlarge $P$,  constructing a connected sum of a punctured
$3$--sphere with (possibly) some copies of $S^2\times S^1$ and
(possibly) some copies of $\rp$. However, we will retain the
property  that the closure of each component of the complement of
our new connected sum  has connected boundary and we will continue
to call this new connected sum $P$. The second condition in the
hypothesis of Theorem \ref{crush} is that for the closure of each
component of the complement of $P$, we can construct a trivial
product region. To establish this, we again may need to enlarge
$P$; again, we maintain that $P$ is a connected sum of a punctured
$3$--sphere with (possibly) some copies of $S^2\times S^1$ and
(possibly) some copies of $\rp$. The third condition in the
hypothesis of Theorem \ref{crush} is that there are no cycles of
truncated prisms in the induced cell decomposition $\C$ on $X$,
which are not in $\bbb{P}(X)$. To establish this condition, we may
need to
 enlarge $P$ again but now we may need to add some copies of $L(3,1)$ as
 connected summands of $P$.

Finally, we arrive at the situation where we have constructed in
$M$ a connected sum of a punctured $3$--sphere  with, possibly,
some copies of $S^2\times S^1$, possibly, some copies of $\rp$
and, possibly, some copies of $L(3,1)$, say $P$, so that either $M
= P$ or if $X$ is the closure of a component of the complement of
$P$ and if $\C$ is its induced cell-decomposition, then $X$ has
connected boundary, $\bbb{P}(\C)\ne X$, each component of
$\bbb{P}(\C)$ is a product $I$--bundle,  $\bbb{P}(X)$ is a trivial
product region for $X$ and there are no cycles of truncated prisms
in $\C$, which are not in $\bbb{P}(X)$. Thus, we can crush the
triangulation $\T$ along $S$. We obtain a triangulation, $\T_{X}$,
of $\hat{X}$. The $3$--manifold $\hat{X}$ is a factor (and so is
$\hat{P}$) in a connected sum decomposition of $M$; furthermore,
the total number of tetrahedra in the triangulations $\T_{X}$ over
all components $X$ in the complement of $P$ is less than the
number of tetrahedra in $\T$. Now, for each $\hat{X}$ we consider
if the triangulation $\T_{X}$ is $0$--efficient. If it is, we set
it aside as one of the desired factors in a connected sum
decomposition of $M$; if it is not, then we go through the above
routine, now with $\hat{X}$ replacing $M$. After a finite number
of steps, we have the desired connected sum decomposition of $M$.

So, this is the idea; we follow with the details.

 The $3$--manifold $M$
is given via some triangulation, say $\T$. By Proposition
\ref{0-decide}, we can decide if $\T$ is $0$--efficient;  if it is
not, then our algorithm will construct a non vertex-linking normal
$2$--sphere. If $\T$ is $0$--efficient, there is nothing to prove;
so, we may assume $\T$ is not $0$--efficient and we have a non
vertex-linking, normal $2$--sphere, $S$.

\vspace{.15 in} \noindent {\it Claim. We can construct a punctured
$3$--sphere connected sum with (possibly) some copies of
$S^2\times S^1$, say $P$, embedded in $M$, so that
 $P$ contains all
the vertices of $\T$, each $2$--sphere component in $\bdy P$ is
normal and not vertex-linking and if $X$ is the closure of a
component of the complement of $P$, then $X$ has precisely one
normal $2$--sphere in its boundary or $P = M$.}

\vspace{.15 in} \noindent {\it Proof of Claim.} By Proposition
\ref{engulf-surface} there are compression bodies $H$ and $H'$
embedded in $M$, $H\cap H' = \bdy_+H = S = \bdy_+H' $, each
component of $\bdy_-H\cup \bdy_-H'$ is a normal $2$--sphere and
$H\cup H'$ contains the vertices of $\T$. Let $P = H\cup H'$. Then
$P$ is a punctured $3$--sphere containing all the vertices of $\T$
and each $2$--sphere in $\bdy P$ is normal and not vertex-linking.

If $P$ has no boundary components, then $M$ is $S^3$. So, we may
assume that $\bdy P \ne \emptyset$. If the closure of each
component $X$ of the complement of $P$ has precisely one boundary
component, then we are done.  So, suppose $X$ is the closure of a
component of the complement of $P$ and $X$ has more than one
boundary component. Hence, by Lemma \ref{oneboundary} there is a
punctured $3$--cell, $P_X$,  embedded in $X$, each component of
$\bdy P_X$ is a normal $2$--sphere, $\bdy X\subset \bdy P_X$ and
no component of the frontier of $P_X$ in $X$ is normally isotopic
to a component of $\bdy X$. It is possible, of course, that there
are no components in the frontier of $P_X$ ($P_X = X$) and so we
have $X$ is itself a punctured $3$--sphere. In any case, we add
the punctured $3$--sphere $P_X$ to $P$ along their intersection,
$\bdy X$. We now may have introduced some copies of $S^2\times
S^1$ as factors in a connected sum decomposition. We continue to
call the resulting $3$--manifold $P$. We consider the closure of a
component of the complement of this new $P$. If such a component
does not have connected boundary, we repeat this procedure. The
only difference from the first stage is that we are adding
punctured $3$--spheres to a connected sum of a punctured
$3$--sphere with some number of copies of $S^2\times S^1$ along
$2$--sphere boundary components; however, the result is still a
punctured $3$--sphere connected sum some number of copies of
$S^2\times S^1$. By Kneser's Finiteness Theorem, Theorem
\ref{kneser}, the procedure must stop and we end up having
constructed a connected sum of a punctured $3$--sphere with some
number of copies of $S^2\times S^1$, again say $P$, satisfying the
conclusions of our claim. Namely, the closure of a component of
the complement of $P$ has precisely one normal, non vertex-linking
$2$--sphere boundary component. This completes the proof of our
Claim.

If the triangulation $\T$ has more than one vertex, there are
other ways to achieve this claim without using the algorithm of
Proposition \ref{0-decide}. For example, we could take a
vertex-linking normal $2$--sphere and then use Proposition
\ref{engulf-sphere} to obtain a punctured $3$--sphere containing
all the vertices of $\T$ with each $2$--sphere in its boundary
normal and not vertex-linking. Similarly, if $\T$ has more than
one vertex, then we can take a collection of pairwise disjoint
vertex-linking normal $2$--spheres, one for each vertex of $\T$
and apply Lemma \ref{oneboundary} where we take for $X$ of that
lemma the closure of the component of the complement of our
vertex-linking normal $2$--spheres, which does not meet any
vertex. Now, we use the punctured $3$--sphere from Lemma
\ref{oneboundary}, which has each of our vertex-linking normal
$2$--spheres in its boundary and add each $3$--cell about the
vertices to get a punctured $3$--sphere containing all the
vertices of $\T$ with each $2$--sphere in its boundary normal and
not vertex-linking.

If we split $M$ along the $2$-spheres in $\bdy P$, from the
previous claim, and cap off the boundary components with
$3$--cells, than $M$ is a connected sum of the resulting
$3$--manifolds. Furthermore,  from our construction we know
precisely how many copies of $S^2\times S^1$ we have in $P$. Let
$X$ be the closure of a component of the complement of $P$. Then
$X$ has only one $2$--sphere in its boundary; denote it by $S$.

We wish to crush the triangulation $\T$ along $S$. If we can do
this, then we will have an ideal triangulation of $\open{X}$;
however, since $S$ is a single $2$--sphere in $\bdy X$, we will
have a triangulation of $\hat{X}$ the manifold obtained by capping
off the boundary $2$--sphere of $X$ with a $3$--cell, which we
have observed is a factor of $M$ in a connected sum decomposition.

 As in the proof of Theorem \ref{0-eff-exists}, no vertices of $\T$ are
in $X$; so,  $X$ has a nice induced cell decomposition consisting
of truncated tetrahedra, truncated prisms, and product triangular
and quadrilateral pieces (See Figure \ref{f-cell-decomp}). Now,
following the above constructions where we crush a triangulation
along a normal surface, we need to establish the conditions in the
hypothesis of Theorem \ref{crush}. Let $\C$ denote the
cell-decomposition on $X$ induced from $\T$ and let $\bbb{P}(\C)$
denote the $I$--bundle consisting of all edges in $\C$, all cells
of type III and IV, along with all trapezoids in the faces of the
cells of types II, III, and IV in $\C$.

\vspace{.15 in}\noindent {\it Claim. We can construct a connected
sum of a punctured $3$--sphere  with , possibly, copies of
$S^2\times S^1$ and, possibly, copies of $\rp$, embedded in $M$,
say $P$, so that either
 $P$ contains all
the vertices of $\T$, each $2$--sphere component in $\bdy P$ is
normal and not vertex-linking, the closure of each component of
the complement of $P$ has connected boundary, and if $X$ is the
closure of a component of the complement of $P$ and $\C$ is its
induced cell-decomposition, then $\bbb{P}(\C)\ne X$ and each
component of $\bbb{P}(\C)$ is a product $I$--bundle or $P = M$.}

\vspace{.15 in}\noindent {\it Proof of Claim.} By our earlier
claim, there is a connected sum of a punctured $3$--sphere  with,
possibly, copies of $S^2\times S^1$, say $P$, embedded in $M$, so
that
 $P$ contains all
the vertices of $\T$, each $2$--sphere component in $\bdy P$ is
normal and not vertex-linking and  the closure of each component
of the complement of $P$,  has precisely one normal $2$--sphere in
its boundary. Let $X$ be the closure of a component of the
complement of $P$ and let $\C$ denote the induced
cell-decomposition of $X$.

If $\bbb{P}(\C) = X$, then since $\bdy X = S$ is connected, $X$ is
a twisted $I$--bundle over $\rpp$; so,
 we add $X$ to $P$ and we have a connected
 sum of a punctured $3$--sphere with, possibly, some copies of $S^2\times S^1$, along with a copy of $\rp$.

If $\bbb{P}(\C) \ne X$ and is not a product $I$--bundle, then
there is a M\"obius band, say $A$, embedded in $X$ with its
boundary in $S$. We may assume that $A$ is normal in the induced
cell-decomposition $\C$ on $X$. A small regular neighborhood
$N(S\cup A)$ of $S\cup A$ is a punctured $\rp$ and its frontier,
say $S'$, is a $2$--sphere, which, by Theorem \ref{barrier}, Item
2, is a barrier surface in the component of its complement not
meeting $S\cup A$. We shrink the $2$-sphere $S'$. We get a
punctured $\rp$, say $P_X$, which contains $N(S\cup A)$, and each
component of the boundary of $P_X$ is either a normal $2$--sphere
or a $2$--sphere contained entirely in the interior of a cell of
$\C$. We fill in each of the $0$-weight $2$--spheres contained
entirely in the interiors of cells of $\C$ with $3$--cells missing
$S\cup A$. We will continue to denote the resulting punctured
$\rp$ by $P_X$. If there are no normal $2$--spheres beside $S$ in
the boundary of $P_X$, then $X = P_X$ is a punctured $\rp$ and we
can again add $X$ to $P$.

So, we assume there are  normal $2$--spheres other than $S$ in the
boundary of $P_X$. In this situation, if the closure of every
component of the complement of $P_X$ in $X$ has precisely one
$2$--sphere in its boundary, then we add $P_X$ to $P$.   If there
is a component of the complement of $P_X$ in $X$ whose closure has
more than one $2$--sphere in its boundary, then we apply Lemma
\ref{oneboundary}. Specifically, if $X'$ is the closure of a
component of the complement of $P_X$ in $X$ and $X'$ has more than
one boundary component, then by, possibly, repeated applications
of Lemma \ref{oneboundary} and Kneser's Finiteness Theorem, there
is a punctured $3$--sphere, $P_{X'}$, embedded in $X'$ so that
$\bdy X'\subset \bdy P_{X'}$ and each component of the complement
of $P_{X'}$ in $X'$ has precisely one normal $2$--sphere boundary
component. We add $P_{X'}$ to $P_X$ and still have a connected sum
of a punctured $\rp$  with, possibly, some factors which are
copies of $S^2\times S^1$. We do this for each component of the
complement of $P_X$ in $X$, which does not have connected
boundary.

So, we have shown that if $\bbb{P}(\C)$ is not a product
$I$--bundle, then either we have that $\hat{X}$ is an $\rp$  or
there is a connected sum of a punctured $\rp$  with, possibly,
some copies of $S^2\times S^1$, say $P_X$, embedded in $X$ so that
$S\subset \bdy P_X$, each component of the boundary of $P_X$
distinct from $S$ is a non vertex-linking normal $2$--sphere not
normally isotopic to $S$, and the closure of a component of the
complement of $P_X$ in $X$ has precisely one  boundary component.
This latter situation adds $2$--spheres to our collection, none of
which are normally isotopic to ones we have; hence, again by
Kneser's Finiteness Theorem, this situation can occur only
finitely many times. If $P_X = X$ for every component $X$ of the
complement of $P$, then we have that $M$ is itself a connected sum
of a punctured $3$--sphere with possibly some copies of $S^2\times
S^1$ and possibly some copies of $\rp$. This completes the proof
of the claim.

\vspace{.15 in} \noindent {\it Claim. We can construct a connected
sum of a punctured $3$--sphere  with, possibly, copies of
$S^2\times S^1$ and, possibly, copies of $\rp$, embedded in $M$,
say $P$, so that either
 $P$ contains all
the vertices of $\T$, each $2$--sphere component in $\bdy P$ is
normal and not vertex-linking, the closure of each component of
the complement of $P$ has connected boundary, and if $X$ is the
closure of a component of the complement of $P$, then we can
construct a trivial induced product region $\bbb{P}(X)$ for $X$ or
$P = M$.}

\vspace{.1 in}\noindent{\it Proof of Claim.} We start with the
conclusion of the preceding Claim and let $P$  be the  connected
sum of a punctured $3$--sphere  with, possibly, copies of
$S^2\times S^1$ and, possibly, copies of $\rp$ in the conclusion
of that Claim. Then
 $P$ contains all
the vertices of $\T$, each $2$--sphere component in $\bdy P$ is
normal and not vertex-linking, the closure of each component of
the complement of $P$ has connected boundary, and if $X$ is the
closure of a component of the complement of $P$ and $\C$ is its
induced cell-decomposition, then $\bbb{P}(\C)\ne X$ and each
component of $\bbb{P}(\C)$ is a product $I$--bundle.

 Let $X$ be the closure of a component of the complement of $P$ and let $\C$ be its
 induced cell-decomposition. Let $K_1\times [0,1],\ldots,K_i\times
[0,1],\ldots,K_n\times [0,1]$ be some ordering of the components
of $\bbb{P}(\C)$. We will, just as in the proof of Theorem
\ref{0-eff-exists}, attempt to construct a trivial product region
for $X$.

Again, the issue is to show that each plug for an appropriate
$K_i\times [0,1]$ is a $3$--cell. So, let us consider an arbitrary
component $K_i\times [0,1]$ of $\bbb{P}(\C)$.

If $K_i$ is simply connected, then let $D_i = K_i$. If $K_i$ is
not simply connected, let $D_{i}^0$ and $D_{i}^1$ be as defined
above. Then, using the above notation, we let $N_i$ be a small
regular neighborhood of $K_i\times [0,1]\cup D_{i}^0\cup D_{i}^1$.
Now, we have a possibly different situation from that above. We
still have that the boundary components of $N_i$ consist of one
properly embedded annulus and some $2$--spheres (there are some
$2$--spheres as we are assuming $K_i$ is not simply connected) and
they form  barrier surfaces in the components of their complements
not meeting $K_i\times [0,1]\cup S$. However, if we let $S_{i,j},
j = 1,\ldots, n_i$ denote the $2$--sphere boundary components of
$N_i$, then a component of the complement of $N_i$ may have more
than one  $S_{i,j}$ in its boundary (an $S_{i,j}$ need not
separate $X$) and the closure of a component of the complement of
$N_i$ in $X$ may be reducible and so has no hope of being a
$3$--cell plug.

So, as in the proof of Theorem \ref{0-eff-exists}, suppose for the
components $K_1\times [0,1],\ldots, K_k\times [0,1]$ of
$\bbb{P}(\C)$, we have simply connected planar complexes
$D_1,\ldots,D_{k'}, k'\le k$ and embeddings $D_j\times [0,1], 1\le
j\le k'$ into $X$ so that $D_j\times\ve = D_{j}^\ve, \ve = 0,1$.
Furthermore, suppose for $j\ne j'$, either $(D_j\times [0,1])\cap
(D_{j'}\times [0,1]) = \emptyset$ or $D_j\times [0,1]\subset
D_{j'}\times [0,1]$ or vice versa; $\bigcup_1^k (K_i\times
[0,1])\subset \bigcup_1^{k'}(D_j\times [0,1])$ and the frontier of
$\bigcup_1^{k'}(D_j\times [0,1])$ is contained in the frontier of
$\bigcup_1^k (K_i\times [0,1])$. So, we have replaced a number of
the $K_i\times [0,1]$ with trivial products. See Figures
\ref{f-product-new} and \ref{f-product-1}.

If $k< n$, consider the component $K_{k+1}\times [0,1]$ of
$\bbb{P}(\C)$. If  $K_{k+1}\times [0,1] \subset
\bigcup_1^{k'}(D_j\times [0,1])$, then $\bigcup_1^{k+1} (K_i\times
[0,1])\subset \bigcup_1^{k'}(D_j\times [0,1])$ and the frontier of
$\bigcup_1^{k'}(D_j\times [0,1])$ is contained in the frontier of
$\bigcup_1^{k+1} (K_i\times [0,1])$ and there is nothing to do.
So, suppose $K_{k+1}\times [0,1]\not\subset
\bigcup_1^{k'}(D_j\times [0,1])$. If $K_{k+1}$ is simply
connected, we set $D_{k'+1} = K_{k+1}$; $D_{k'+1}\times [0,1]$ is
disjoint from each $D_j\times [0,1], 1\le j\le k'$. So, suppose
$K_{k+1}$ is not simply connected. We let $D_{k+1}^\ve, \ve = 0,1$
be defined as above and let $N_{k+1} = N((K_{k+1}\times [0,1])\cup
D_{k+1}^0\cup D_{k+1}^1)$ denote a small regular neighborhood of
$(K_{k+1}\times [0,1])\cup D_{k+1}^0\cup D_{k+1}^1$ and let $N =
N(\bigcup_1^{k'}(D_j\times [0,1])$ be a small regular neighborhood
of $\bigcup_1^{k'}(D_j\times [0,1])$. The frontier of $N$ consists
of annuli and is a barrier surface in the component of its
complement in $X$ not containing $\bigcup_1^{k'}(D_j\times
[0,1])$. If we denote the $2$--spheres in the frontier of
$N_{k+1}$ by $S_{k+1,j}, j = 1,\ldots,n_{k+1}$, then we can shrink
the $S_{k+1,j}$ in the complement of $(K_{k+1}\times [0,1])\cup S$
and $\bigcup_1^{k'}(D_j\times [0,1])$. After shrinking the
$S_{k+1,j}$, we have a collection of $0$--weight $2$--spheres
entirely contained in the interior of cells of $\C$ and, possibly,
some normal $2$--spheres. We can fill in any of the $0$--weight
$2$--spheres with $3$--cells missing $(K_{k+1}\times [0,1])\cup
S\cup \bigcup_1^{k'}(D_j\times [0,1])$. Suppose we have done this.
If there are no normal $2$--spheres, then it follows that each
$S_{k+1,j}$ separates $X$ and each plug for $N_{k+1}$ is a
$3$--cell. However, if there are normal $2$--spheres, then we let
$\mathcal{S}$ denote this collection of normal $2$--spheres along
with the $2$--sphere $S$.  We apply Lemma \ref{oneboundary} using
the collection $\mathcal{S}$. It follows, there is a punctured
$3$--sphere, $P_{X}$, embedded in $X$, each component of $\bdy
P_{X}$ is a normal $2$--sphere, $\bdy X\subset \bdy P_{X}$ and no
component of the frontier of $P_{X}$ in $X$ is normally isotopic
to a component of $\bdy X$. It is possible that there is a
component of the complement of $P_{X}$ in $X$ whose closure does
not have connected boundary, again, we can enlarge $P_{X}$ in a
finite number of steps (Kneser's Finiteness Theorem) so that we
have a punctured $3$--sphere connected sum with, possibly, some
copies of $S^2\times S^1$, which we  continue to call $P_X$. Thus
we have $P_X$ embedded in $X$, $S$ is in the boundary of $P_X$ and
the closure of each component of the complement of $P_X$ has
connected boundary. We add $P_X$ to $P$ and continue to use $P$
for the resulting connected sum of a punctured $3$--sphere  with,
possibly, some copies of $S^2\times S^1$ and, possibly, some
copies of $\rp$.

Now, it is possible that our new $P$ does not satisfy the
conclusions of our previous Claim; namely, we no longer know that
for $X$ the closure of a component of the complement of $P$ and
$\C$ its induced cell-decomposition that $\bbb{P}(\C)\ne X$ and
each component of $\bbb{P}(\C)$ is a product $I$--bundle. However,
we can go back and fix this and by Kneser's Finiteness Theorem, we
will only need to go back a finite number of times. Thus, we
eventually are able to construct a $P$ satisfying the conclusions
of our claim.

\vspace{.15 in} \noindent {\it Claim.  We can construct a
connected sum of a punctured $3$--sphere  with, possibly, copies
of $S^2\times S^1$, possibly, copies of $\rp$, and, possibly,
copies of the lens space $L(3,1)$ embedded in $M$, say $P$, so
that
 $P$ contains all
the vertices of $\T$, each $2$--sphere component in $\bdy P$ is
normal and not vertex-linking, the closure of each component of
the complement of $P$ has connected boundary, and if $X$ is the
closure of a component of the complement of $P$ and $\C$ is the
induced cell-decomposition on $X$, then all the conditions in the
hypothesis of Theorem \ref{crush} are satisfied to crush the
triangulation along the normal $2$--sphere in the boundary of $X$
or $P = M$.}

\vspace{.1 in} \noindent {\it Proof of Claim.} We have, from the
above Claims, constructed a connected sum of a punctured
$3$--sphere with, possibly, copies of $S^2\times S^1$ and,
possibly, copies of $\rp$, say
 $P$, which contains all
the vertices of $\T$, each $2$--sphere component in $\bdy P$ is
normal and not vertex-linking, the closure of each component of
the complement of $P$ has connected boundary, and if $X$ is the
closure of a component of the complement of $P$ and $\C$ is its
induced cell-decomposition, then we can construct a trivial
induced product region $\bbb{P}(X)$ for $X$ (in particular,
$\bbb{P}(\C)\ne X$ and each component of $\bbb{P}(\C)$ is a
trivial product $I$--bundle). We need to show that $P$ can be
constructed so that, in addition, if $X$ is the closure of a
component of the complement of $P$ and $\C$ is its induced
cell-decomposition, then there are no cycles of truncated prisms
in $\C$, which are not in $\bbb{P}(X)$.

As in the proof of Theorem \ref{0-eff-exists}, there are two types
of cycles of truncated prisms: one is a cycle about an edge $e$ of
$\C$ (see Figure \ref{f-cycle}(A)) and the other cycles about more
than one edge of $\C$ (see Figure \ref{f-cycle}(B)).

If there is a complete cycle about a single edge $e$ as in Figure
\ref{f-cycle}(A), then a surgery on $S$ at $D$ gives two {\it
normal} $2$--spheres $S_{0}$ and $S_1$, neither of which is
vertex-linking. See Figure \ref{f-surgery-cycle}. Furthermore, $S$
along with $S_0$ and $S_1$ bound a punctured $3$--cell, say $P_X$.
As before, it is possible that the closure of a component of the
complement of $P_X$ in $X$ does not have connected boundary,
neither $S_0$ nor $S_1$ separates. We apply Lemma
\ref{oneboundary} to the closure of the complement of $P_X$ and if
necessary, we repeatedly apply Lemma \ref{oneboundary} until we
have that there is a connected sum of a punctured $3$--sphere,
possibly, with some copies of $S^2\times S^1$ embedded in $X$,
which we continue to denote $P_X$, $S$ is in the boundary of $P_X$
along with possibly other normal $2$--spheres distinct from $S$
and the closure of each component of the complement of $P_X$ in
$X$ has connected boundary. We add $P_X$ to $P$. Now, we may need
to go back and construct an even larger $P$ so that the hypotheses
of this Claim are still satisfied; however, again, we will need to
do this at most a finite number of times by Theorem \ref{kneser}.

If there is a complete cycle about more than one edge (see Figure
\ref{f-cycle}(B)), then, as above, the collection of cells of type
II (truncated prisms), form a solid torus with, possibly, some
self identifications in its boundary. Using the same notation as
earlier, we use $\hat{\tau}$ to denote the cycle of truncated
prisms and $\tau$ to denote the cycle of truncated prisms minus
the bands of trapezoids. If $\tau\cap S$ is three open annuli,
then, as above, it follows that $\hat{\tau}$ is in the induced
product region for $X$, $\bbb{P}(X)$.

So, we suppose that $\tau\cap S$ is a single open annulus, meeting
a meridional disk of $\tau$ three times. See Figure
\ref{f-cycle-torus}. Our argument here is a bit different from
that above in the similar situation, as we do not have that the
manifold $M$ is irreducible. We now consider the cycle of
truncated prisms with the trapezoidal faces slightly pulled into
the truncated prisms. We get an embedded solid torus, which we
will denote by $\tau'$, which meets $S$ in a single annulus $A$
and $A$ meets the meridional disk of $\tau'$ three times. Let $A'$
denote the closure of the annulus in $\bdy\tau'$ complementary to
$A$.

Now, let $N$ be a small regular neighborhood of $\tau' \cup S$.
$N$ is the punctured lens space $L(3,1)$ and its frontier is a
$2$--sphere $S'$, which is a barrier surface in the component of
its complement not meeting $\tau'\cup S$ ($S'$ is parallel to the
$2$--sphere $(S\setminus A)\cup A'$). We shrink $S'$ and obtain a
punctured $L(3,1)$ having $S$ in its boundary along with,
possibly, some normal $2$--spheres and some $0$--weight
$2$--spheres contained entirely in the interiors of tetrahedra.
This is a familiar situation and after using Lemma
\ref{oneboundary}, Kneser's Finiteness Theorem (Theorem
\ref{kneser}), and, if necessary, returning to the previous steps,
we obtain a connected sum of a punctured $L(3,1)$ along with,
possibly, some copies of $S^2\times S^1$ and $\rp$, say $P_X$. We
add $P_X$ to $P$. This is what adds the possibility that we now
have copies of $L(3,1)$ in our connected sum.

This completes the proof of this Claim.

As we remarked earlier, if $X$ is the closure of a component of
the complement of $P$, then $\hat{X}$ and $\hat{P}$ are factors in
a connected sum decomposition of $M$. The $3$--manifold $\hat{P}$
is a connected sum of $S^3$ and (possibly) some copies  of
$S^2\times S^1$, $\rp$ and $\L(3,1)$. If $X$ is the closure of a
component of $P$, then we can crush the triangulation along the
$2$--sphere in the boundary of $X$; hence, we construct from the
truncated tetrahedra in $\C$ which are not in $\bbb{P}(X)$,  an
ideal triangulation $\T_X^*$ of $\open{X}$. However, since
boundary $X$ is a single $2$--sphere, $\abs{\T_X^*}$ is
homeomorphic to $\widehat{X}$ and $\T_X^*$ is a triangulation of
$\widehat{X}$. Since each $2$--sphere in the boundary of $P$ is
non vertex-linking, there is at least one truncated prism in $X$
and so there are strictly fewer truncated tetrahedra in $\C$ than
there are in $\T$. It follows that the triangulation $\T_X^*$ of
$\hat{X}$ has fewer tetrahedra than the triangulation $\T$. In
fact, the entire number of tetrahedra in the triangulations
obtained by crushing $\T$ along each $2$--sphere in the boundary
of $P$ is less than the number of tetrahedra in $\T$.

Now, we consider each $\hat{X}$ and triangulation $\T_X^*$ of
$\hat{X}$. If $\T_X^*$ is $0$--efficient, then we are satisfied;
if not, then we construct in $\hat{X}$ a punctured $3$--sphere
connected sum with, possibly, copies of $\S^2\times S^1$, copies
of $\rp$ and copies of $L(3,1)$ embedded in $\hat{X}$, say $P_X$,
so that
 $P_X$ contains all
the vertices of $\T_X^*$, each $2$--sphere component in $\bdy P_X$
is normal and not vertex-linking, the closure of each component of
the complement of $P_X$ has connected boundary, and if $X'$ is the
closure of a component of the complement of $P_X$ and $\C'$ is the
induced cell-decomposition on $X'$, then all the conditions in the
hypothesis of Theorem \ref{crush} are satisfied to crush the
triangulation along the normal $2$--sphere in the boundary of
$X'$. Since the number of tetrahedra in each of our triangulations
obtained after crushing is strictly decreasing, and any factors in
a connected sum decomposition of $\hat{X}$ are factors in a
connected sum decomposition of $M$, the process must terminate in
the desired connected sum decomposition of $M$.

This completes the proof of Theorem \ref{reduce}.\end{proof}

Next we show that given a triangulation of a closed, orientable
$3$--manifold $M$, there is a construction of a maximal, pairwise
disjoint collection of distinct normal $2$--spheres. As an
immediate corollary of this construction, we get an irreducible
decomposition of $M$. This  can also be used as an alternate
approach to Theorem \ref{reduce}. We are using that a pairwise
disjoint collection of distinct normal $2$--spheres $\mathcal{S}$
is {\it maximal} if and only if for $\S'$ a pairwise disjoint
collection of distinct normal $2$--sphere with $\S\subset \S'$,
then $\S' = \S$.

\begin{thm}\label{max} Suppose $\T$ is a triangulation of the closed,
orientable $3$--manifold $M$. There is an algorithm to construct a
maximal, pairwise disjoint collection of distinct normal
$2$--spheres in $M$.\end{thm}

\begin{proof} Let $\S_1$ be the (pairwise disjoint) collection of
distinct, normal, vertex-linking $2$--spheres in $M$.  Now, we
want to know if there is a normal $2$--sphere disjoint and
distinct from those in the collection $\S_1$. To do this we use a
re-writing of the normal surfaces in the complement of the
collection $\S_1$.

Let $M_1$ denote  the manifold obtained by splitting $M$ along the
$2$--spheres in $\S_1$ and let $\C_1$ be the induced
cell-decomposition on $M_1$. Suppose there is a normal $2$--sphere
in $M$, which is disjoint and distinct from the $2$--spheres in
$\S_1$. Then there is a normal $2$--sphere in $M_1$, which is
disjoint from the boundary of $M_1$ and is not normally isotopic
into a $2$--sphere in the boundary of $M_1$. Recall from Lemma
\ref{0-decide} we proved that if there is a non vertex-linking
normal $2$--sphere at all, then there is one at a vertex of the
projective solution space. We shall do the same here but we
replace the projective solution space $\P(M,\T)$ by the projective
solution space of the sub-cone of parameterizations of normal
surfaces in $\C_1$ obtained by setting all variables corresponding
to normal disks types which meet the boundary of $M_1$ equal to
zero, say $\overline{\P}(M_1,\C_1)$. To this end suppose $S_1$ is
a normal $2$-sphere in $M_1$, which is disjoint and distinct from
the $2$--spheres in $\bdy M_1$; suppose $\mathcal{C}(S_1)$ is the
carrier of $S_1$; and $S_1$ has been chosen so that among all
normal $2$--spheres in $M_1$, which are disjoint and distinct from
the $2$--spheres in $\bdy M_1$, the dimension of
$\mathcal{C}(S_1)$ is a minimum. If $\mathcal{C}(S_1)$ is not a
vertex of the projective solution space,
$\overline{\P}(M_1,\C_1)$, then there are normal surfaces $X$ and
$Y$ carried by proper faces of $\mathcal{C}(S_1)$ and nonnegative
integers $k, n,$ and $m$ with $kS_1 = nX + mY$. Since, $\chi(S_1)
> 0$, we must have either $\chi(X) > 0$ or $\chi(Y)
> 0$ and so a component of one has positive Euler characteristic,
say a component $X'$ of $X$. The carrier of $X'$ is also a proper
face of $\mathcal{C}(S_1)$. If $X'$ is a projective plane, then
its double is a normal $2$--sphere carried by a proper face of
$\mathcal{C}(S_1)$ and so, as a proper face of $\mathcal{C}(S_1)$,
has dimension less than that of $\mathcal{C}(S_1)$. The normal
surface $2X'$ is disjoint from the $2$--spheres in $\bdy M_1$ and
so, the only possibility by our choice of $S_1$ is that $2X'$ is
normally isotopic into a component of $\bdy M_1$. But then $2X'$
would have to be a component of the sum $kS_1 = nX + my$, which
contradicts our choice of $S$. We get a similar contradiction if
$X'$ is, itself, a $2$--sphere. Thus $S_1$ has to be carried by a
vertex of the projective solution space $\overline{\P}(M_1,\C_1)$.
The vertices of  this projective solution space can be constructed
and
 we can decide if a normal surface is a
$2$--sphere distinct from a $2$--sphere in $\bdy M_1$. We let
$\S_2 = \S_1\cup S_1$.

Now, suppose we have constructed a collection $\S_n$ of pairwise
disjoint, distinct normal $2$--spheres, $n\ge 2$. Let $M_{n+1}$
denote the manifold obtained by splitting $M$ along the
$2$--spheres in the collection $\S_n$. Then just as above, there
is any normal $2$--sphere in $M$, which is disjoint and distinct
from the $2$--spheres in the collection $\S_n$, if and only if in
re-writing of the normal surfaces, there is a normal $2$--sphere
which is a vertex solution.  However, by Kneser's Finiteness
Theorem, we eventually have a collection $\S$ so there are no
normal $2$--spheres which are disjoint and distinct from the
$2$--spheres in this collection.
\end{proof}

The preceding proof uses the re-writing process. Each $2$--sphere
in our collection is normal in $M$. However, we find these
$2$--spheres by re-writing normal surfaces in $M$ as normal
surfaces in the induced cell-decomposition of $M$ split along a
pairwise disjoint collection of normal $2$--spheres. In this
re-writing, we show that if there is a normal $2$--sphere disjoint
and distinct from our collection, then there is one  at a vertex
in a new projective solution space determined by re-writing the
parametrization of the normal surfaces in $M$, which are disjoint
from our collection of $2$--spheres. Hence, we can decide if there
is such a $2$--sphere and we know when we have indeed constructed
a maximal collection. In \cite{jac-ree}, it is shown that a
maximal pairwise disjoint collection of distinct normal
$2$--spheres actually exists at the vertices of $\P(M,\T)$; and in
fact, the $2$--spheres in the collection have parameterizations
which project to vertices of a face of $\P(M,\T)$, which is a
simplex. However, in Theorem \ref{max}, we do not necessarily have
that the $2$--spheres in our collection have projective
representatives which are vertices (or for that matter, even
fundamental). It is not clear which of these algorithms might be
the most efficient in constructing such a maximal collection.

There is a straight forward observation about the way a maximal
collection of normal $2$--spheres decompose a $3$--manifold, which
was used by the authors in 1988 and is well known. It also was
used by both Rubinstein and Thompson in their proofs of the
$3$--sphere recognition algorithm. This observation follows
immediately from Proposition \ref{oneboundary}.

\begin{remark}\label{r-form} Suppose $\T$ is a triangulation of the closed, orientable $3$--manifold
$M$ and $\S$ is a maximal, pairwise disjoint collection of
distinct normal $2$--spheres in $M$. Then if $X$ is a component of
$M$ split along $\S$, $X$ satisfies one of the following:
\begin{enumerate}
\item $X$ has connected boundary, which is a vertex-linking normal
$2$--sphere, and $X$ is $3$--cell about a vertex of $\T$, \item
$X$ has disconnected boundary and is a punctured $3$--sphere, or
\item $X$ has connected boundary but is not of the
 first type. \end{enumerate}\end{remark}

 It seems possible that if $X$ satisfies Item 2 (does not have
 connected boundary), then there can be at most three $2$--spheres
 in boundary $X$. This does not seem to be important but is
 a curiosity. Below we also leave open a (seemingly)
 related question about $0$--efficient triangulations of the
 $3$--cell.

We remarked earlier that Theorem \ref{max} gives an alternate
proof of the decomposition obtained in Theorem \ref{reduce}. To
see this, suppose $\T$ is a triangulation of the closed,
orientable $3$--manifold $M$ and $\S$ is a maximal, pairwise
disjoint collection of distinct normal $2$--spheres in $M$. Now,
if $X$ is a component of $M$ split along the $2$--spheres in the
collection $\S$ and $X$ is a punctured $3$--cell, then $X$
determines a factor in a connected sum decomposition of $M$ which
is either $S^3$ or a connected sum of some number of copies of
$S^2\times S^1$. On the other hand, if $X$ is as in Item 3 (Remark
\ref{r-form}), then $\hat{X}$, the manifold obtained from $X$ by
filling in its boundary $2$--sphere with a $3$--cell, is a factor
in a connected sum decomposition of $M$. We attempt to crush the
triangulation $\T$ along the $2$--sphere in $\bdy X$. If there are
obstructions to crushing the triangulation, then since the
collection $\S$ is maximal, we have that $\hat{X}$ is either $S^3,
\rp$ or $L(3,1)$. Otherwise, we can crush the triangulation $\T$
along $\bdy X$. In this way we get a one-vertex triangulation
$\T_X$ of $\hat{X}$. We observe that the number of tetrahedra in
$\T_X$ is strictly smaller than the number of tetrahedra in $\T$
if $\S$ has at least one $2$--sphere which is not vertex-linking.
So, we can use our argument on $\hat{X}$, constructing a maximal
collection of pairwise disjoint, distinct normal $2$--spheres in
$\hat{X}$. Since the number of tetrahedra is strictly decreasing,
the process must stop with a decomposition of $M$ into factors
which are either $S^3, S^2\times S^1, \rp, L(3,1)$ or have
triangulations, which do not have any non vertex-linking normal
$2$--spheres; hence, these latter factors have $0$--efficient
triangulations.

We suspect that if we have a maximal collection of pairwise
disjoint, distinct, normal $2$--spheres embedded in $M$, say $\S$,
and $X$ is a component of $M$ split along $\S$ which is as in Item
3 of Remark \ref{r-form}, then if the triangulation $\T$ can be
crushed along $\bdy X$, then the resulting triangulation is
immediately $0$--efficient. However, we have not been able to
establish this. Certainly if $\T_X$ is the resulting one-vertex
triangulation of $\hat{X}$ and if $S$ is a non vertex-linking
normal $2$--sphere in $\hat{X}$, it seems likely $S$ would give
rise to a normal $2$--sphere embedded in $X$ which is distinct
from $\bdy X$ and thereby contradict the maximality of the
collection $\S$. Such a $2$--sphere is the image of a $2$--sphere
in $X$ (crushing is a cellular map) but it may not be normal and
in shrinking, it may become equivalent to $\bdy X$. We have not
spent much time considering this. Recall that our basic philosophy
in this work is to crush a triangulation along a normal
$2$--sphere as soon and as often as we can. The construction of a
maximal collection of pairwise disjoint, distinct, normal
$2$--spheres  seems to be computationally quite expensive, whether
using the method above or the method in \cite {jac-ree}. On the
other hand, crushing a triangulation along an embedded normal
surface seems to quite sharply reduce the complexity of any
further computations.

Recall the solution to the $3$-sphere recognition problem in both
\cite {rub, tho} uses almost normal $2$--spheres. In one
direction, it is argued that in any triangulation of $S^3$ there
is an embedded almost normal $2$--sphere. Conversely, under
special circumstances, namely, in each component of the complement
of a maximal collection of embedded, normal $2$--spheres, it is
argued that one can decide if there is an almost normal
$2$--sphere and if there is one, then that component is a
$3$--cell. One reason for investigating $0$-efficient
triangulations is that they provide a constructible and natural
environment for this latter direction in the $3$--sphere
recognition algorithm. This is discussed in detail in
\cite{jac-let-rub1, jac-let-rub2}. Here we give only this latter
direction in the $3$--sphere recognition algorithm. We note,
however, the subtleness for getting an almost normal, {\it
octagonal} $2$--sphere in the $3$--sphere recognition algorithm.
There are triangulations of $S^3$ where the only almost normal
$2$--spheres are tubed; for example, the triangulation of the
$3$--sphere coming from the boundary of the $4$--simplex. We do,
however, observe in the proof of Lemma \ref{anoctagonal} below,
that in a one-vertex, $0$--efficient triangulation an almost
normal $2$--sphere must be octagonal.

\begin{thm} [$3$--Sphere Recognition Problem \cite{rub,tho}]  Given a $3$--manifold $M$,
it can be decided if
$M$ is homeomorphic with the $3$--sphere.\end{thm}

Suppose we are given a $3$--manifold $M$ via a triangulation $\T$.
By  Theorem \ref{reduce}, we can construct  a connected sum
decomposition of $M$, where we know  each factor is either $S^3$,
$S^2\times S^1,  \rp, L(3,1)$ or has a $0$--efficient
triangulation. Of course, if any of the factors are $S^2\times
S^1$, $\rp$ or $L(3,1)$, then $M$ is not $S^3$; in practice, we
compute homology invariants directly from $\T$, which indicate if
any
 factors having nontrivial homology will show up in a connected sum decomposition. So, we may
assume that in our connected sum decomposition, we have only
factors which are either known to be $S^3$ or ones with
$0$--efficient triangulations and about which we know nothing else
(we know they are homology $3$--spheres). We now employ the next
proposition, which exhibits a basic feature of $0$--efficient
triangulations and is one of our motivating reasons for
constructing them.

\begin{prop}\label{anoctagonal} Suppose the $3$--manifold $M$ has a $0$--efficient
triangulation. Then it can be decided if $M$ has an almost normal
$2$--sphere. Furthermore, if $M$ has a $0$--efficient
triangulation and an almost normal $2$-sphere, then $M = S^3$.
\end{prop}

\begin{proof} Suppose $\T$ is a $0$--efficient triangulation of $M$. If $\T$ has more than
one vertex, then by Proposition \ref{p-0-eff}, $M$ is
$S^3$. Furthermore, $M$ has a tubed almost normal $2$--sphere determined by taking a copy of
each vertex-linking normal $2$--sphere and a tube along an edge joining them. So, if $\T$ has
two vertices, it is easy to find an almost normal (tubed) $2$--sphere and $M$ is the
$3$--sphere. Hence, we may assume
$\T$ has only one vertex. Now, with
$\T$ having only one vertex and being $0$--efficient, an almost normal $2$--sphere must be
octagonal. For if there is an almost normal
tubed, $2$--sphere, then a compression of the tube gives two normal $2$--spheres. Since
$\T$ is $0$--efficient both are vertex-linking and so must be the same $2$--sphere.
But in this case, the only possibility is that they are tubed through the normal product
region. A tube through the product region between two copies of a normal surface does not
give an almost normal surface.

Hence, if there is an almost normal $2$--sphere, there is an
octagonal one. Let $\Sigma$ be an octagonal almost normal
$2$--sphere so that $wt(\Sigma)$, the weight of $\Sigma$, is a
minimum among all such almost normal $2$--spheres. We claim,
$\Sigma$ is fundamental. For, if not, then there are normal and
almost normal surfaces $X$ and $Y$ so that $\Sigma = X + Y$ is a
non trivial Schubert sum; that is, over all possible ways to
express $\Sigma = X + Y$ as a geometric sum, we have chosen one
with $X\cap Y$ having a minimal number of components. But
$\chi(\Sigma)
> 0$; so, since, $\chi(\Sigma) = \chi(X) +\chi(Y)$, we have
$\chi(X)
> 0$, say. $X$ can not be a normal $\rpp$, since $\T$ is
$0$--efficient; if $X$ were an almost normal $\rpp$, then $Y$
would need to be a normal $\rpp$. Again, this leads to a
contradiction. (Note: There can not be a normal or almost normal
$\rpp$ in a $0$--efficient triangulation. An almost normal one
shrinks to a normal one and a normal $\rpp$ has a double which is
a non vertex-linking normal $2$--sphere.) So, the only possibility
is for $X$ to be a normal or an almost normal $2$--sphere. $X$ is
not normal, since the only normal $2$--spheres are vertex-linking
and hence are a component in any geometric sum. And by our choice
of $\Sigma$, $X$ is not almost normal, since $wt(X) < wt(\Sigma)$.
We conclude, $\Sigma$ must be fundamental.

Now, if we have an almost normal $2$--sphere $\Sigma$ in a
$0$--efficient triangulation of $M$, then since $\Sigma$ is
unstable to each side, we have that $\Sigma$ shrinks in each
component of its complement to a collection of $2$--spheres, which
are either normal or are $0$--weight and embedded entirely in the
interior of a tetrahedron. But in a $0$--efficient triangulation
any normal $2$--sphere is vertex-linking. Hence, in either case,
it follows that each component of the complement of an almost
normal $2$--sphere is a $3$--cell. It follows in this situation
that $M = S^3$.
\end{proof}

Our argument in Proposition \ref{anoctagonal}  does not conclude
that if in a $0$--efficient triangulation  there is an almost
normal $2$--sphere, there is an almost normal $2$--sphere which is
a vertex solution in the projective solution space. However, the
following argument, while very similar, provides a clever variant
to the standard argument for showing a solution is a vertex
solution and is from A.~Casson. Recall if there is an embedded
normal or almost normal surface $\Sigma$, then every surface
carried in a face of the carrier of $\Sigma$ is embedded and all
the surfaces in the cone over the carrier of $\Sigma$ are
compatible; i.e., their geometric sum is defined. Also, as we
observed in the proof of the previous lemma, in a one-vertex,
$0$--efficient triangulation, there are no almost normal tubed
$2$--spheres. Hence, we need only consider octagonal almost normal
surfaces. Instead of the weight functional on normal and almost
normal surfaces, a different functional is used. In the cone, over
the carrier of an almost normal surface a solution may have
multiple copies of an almost normal octagon or almost normal tube;
in our situation we are only concerned with octagonal surfaces. If
$F$ is a solution over the carrier of $\Sigma$, we let $O(F)$
denote the number of octagons in $F$. If $\chi(F)$ is the Euler
characteristic of $F$, we use the linear functional $L(F) =
\chi(F) - O(F)$. Note that in a closed $3$--manifold whenever we
have $\Sigma$ a connected normal or almost normal surface and
$L(\Sigma) > 0$, then $\Sigma$ is a normal or an almost normal
$2$--sphere or a normal $\rpp$; and as we have noted, in a
$0$--efficient triangulation, there is no normal $\rpp$.

Now, suppose we have a $0$--efficient triangulation and there is
an embedded almost normal $2$--sphere. Let $\Sigma$ be one that
has the lowest dimensional carrier; we will show that it is a
vertex of the projective solution space. For, if not, then there
are nonnegative integer solutions $X$ and $Y$ to the normal
equations and nonnegative integers $k, n$ and $m$ so that $k\Sigma
= nX + mY$. It follows that either $nL(X) > 0$ or $mL(Y) > 0$; say
$nL(X)
> 0$. We then have $L(X) > 0$ and hence a component of $X$, say $X'$, is a
 normal $2$--sphere or an almost normal $2$--sphere.
However, in the former case, since we have a $0$--efficient
triangulation, we have that $X'$ is a vertex-linking normal
$2$--sphere and so, $\Sigma$ would have a component a
vertex-linking normal $2$--sphere, which is a contradiction. In
the latter case, it follows there is an almost normal $2$--sphere
carried by a proper face of the carrier of $\Sigma$. However, this
situation contradicts the choice of $\Sigma$. Hence, $\Sigma$ must
be itself a vertex solution. This completes this argument.

So, a special feature of a $0$--efficient triangulation is that in
such a  triangulation, it can be decided if there is an embedded
almost normal $2$--sphere. Furthermore, if there is one, then the
$3$--manifold must be $S^3$.

We comment that in seeking almost normal surfaces, one can
consider the subspace of solutions obtained by setting all normal
octagon types equal to zero except for one in some tetrahedron;
similarly for that tetrahedron, we can also set all normal quad
types in it equal to zero. In this way we have at most $3^t$ such
almost normal solution spaces to consider when searching for
almost normal surfaces.

To complete the $3$--sphere recognition algorithm, after arriving
at a connected sum decomposition of $M$ into factors, which are
known to be either $S^3$ or to have $0$--efficient triangulations,
we can check each factor with a $0$--efficient triangulation to
see if it has an almost normal $2$--sphere. By our above argument,
this can be decided in a $0$--efficient triangulation. If the
answer is yes, then, again by the preceding arguments, that factor
is $S^3$. Of course, now comes a hard part of the $3$--sphere
recognition algorithm. We use the work of \cite{rub,tho} at this
point for the case when the answer is no. As pointed out above,
from \cite{rub,tho} it is known that if a factor is $S^3$, then
for any triangulation, it must have an almost normal $2$--sphere;
so, this would be true for our $0$--efficient triangulations.
Hence, if the answer is no for a factor, then that factor is not
$S^3$; and so, $M \neq S^3$.

An algorithm can be deduced from \cite{rub} to  decide if there is
an almost normal $2$--sphere in an arbitrary triangulation of a
$3$--manifold. However, notice that in an arbitrary triangulation,
the shrinking of an almost normal $2$-sphere may get hung up on a
normal $2$--sphere which we may know nothing about; hence, the
manifold may not be $S^3$. There are easily constructed one-vertex
triangulations of lens spaces, for example, which have almost
normal octagonal $2$--spheres. Also, any triangulation of a
$3$--manifold with more than one vertex has a tubed almost normal
$2$--sphere.

We have as an easy corollary to Theorem \ref{reduce} and the
$3$--sphere recognition algorithm, recognition algorithms for
$\rp$ and $S^2\times S^1$. This was observed earlier by the second
author in \cite{rub} and by M. Stocking in \cite{sto}.

\begin{cor} Given a closed, orientable $3$--manifold $M$, it can
be decided if $M$ is homeomorphic to $\rp$ or to $S^2\times S^1$.
In particular, it can be decided if $M$ is a connected sum of some
number of copies of $\rp$ along with some number of copies of
$S^2\times S^1$.\end{cor}

\begin{proof} By Theorem \ref{0-decide}, we can construct a
connected sum decomposition of $M$ into copies of $\rp, S^2\times
S^1, L(3,1)$ along with some number of closed $3$--manifolds each
having a $0$--efficient triangulation. If there are no copies of
$\rp$ or $S^2\times S^1$ or if there is a copy of $L(3,1)$, then
we know $M$ is not a connected sum of some number of copies of
$\rp$ and $S^2\times S^1$. So, we may assume there are some copies
of $\rp$ or $S^2\times S^1$ and no copies of $L(3,1)$. Each of the
other factors in our connected sum decomposition of $M$ has a
$0$--efficient triangulation; hence, we can decide, using the
$3$--sphere recognition algorithm in these $0$--efficient
triangulations, if such factors are $S^3$. If some factor with a
$0$--efficient triangulation is not $S^3$, then $M$ is not a
connected sum of some number of copies of $\rp$ or $\S^2\times
S^1$. On the other hand, if all factors with a $0$--efficient
triangulation are $S^3$, then $M$ has been decomposed into a
connected sum of some number of copies of $\rp$ and $S^2\times
S^1$ and we have the desired result.\end{proof}

Note that we are very close in the previous Corollary to being
able to decide if $M$ is $L(3,1)$. However, there are
$0$--efficient triangulations of $L(3,1)$; so, one of the factors
with $0$--efficient triangulation may be $L(3,1)$. In
\cite{jac-rub1}, we give an algorithm to decide if a $3$--manifold
is a lens space and in particular, an algorithm which decides
precisely which lens space. Similar results have been obtained by
\cite{sto}.

 Our next result is a constructive method to alter a given
triangulation to a $0$--efficient one. It uses the $3$--sphere
recognition algorithm and so it is not the most preferred in
practice. Also, it follows from Theorem \ref{reduce} (and the
$3$--sphere recognition algorithm) as we will point out in its
proof. However, we give here a proof, which implements our basic
philosophy ``crush first and ask questions later."

\begin{thm} \label{0-eff-construct}Suppose $M$ is a closed, orientable, irreducible $3$--manifold.
Then any triangulation of $M$ can be modified to a $0$--efficient
triangulation or it can be shown that $M$ is one of $S^3, \rp$ or
$L(3,1)$.
\end{thm}

\begin{proof} First, we point out how this Theorem follows from Theorem
\ref{0-decide}. We can construct $3$--manifolds $M_1,\ldots, M_n$
such that $M = M_1\#,\ldots,\#M_n$ and each $M_i$ is either $S^3,
\rp, S^2\times S^1, L(3,1)$ or has a $0$--efficient triangulation.
However, since $M$ is irreducible, we can not have $S^2\times S^1$
and if we have $\rp$ or $L(3,1)$ as a factor, then $M$ is itself
$\rp$ or $L(3,1)$. Thus, the only possibility is that we have $M$
is $S^3$ or there are only factors which are known to be $S^3$ and
some factors which have $0$--efficient triangulations. We employ
the $3$--sphere recognition algorithm on each of the
$0$--efficient factors. If any one is not $S^3$, then that factor
must be our manifold $M$ and it has a $0$--efficient
triangulation. If all factors are $S^3$, then $M$ is $S^3$. So,
this completes the argument using Theorem \ref{reduce}.

So, more directly, and the method employed by Letscher and Burton
in their computer program, REGINA.

Suppose $\T$ is the given triangulation of the $3$--manifold $M$.
By the Proposition \ref{0-decide}, we can decide if $\T$ is
$0$--efficient. If it is, then there is nothing to prove.  So, we
may suppose that $\T$ is not $0$--efficient and by Proposition
\ref{0-decide} our algorithm has constructed a  non
vertex-linking, normal $2$--sphere, say $S$. Again, if $\T$ has
more than one vertex, we can construct directly a
non-vertex-linking normal $2$--sphere or we know that $M = S^3$.
Since $M$ is irreducible, $S$ separates and, in fact, bounds a
$3$--cell; however, we do not know, {\it a priori}, which side of
$S$ is a $3$--cell.

Now, suppose we have constructed a non vertex-linking normal
$2$--sphere $S$. By Proposition \ref{engulf-surface}, we can
construct a punctured $3$--sphere $P$, which contains $S$ and all
the vertices of $\T$ or $M$ is $S^3$. Furthermore, since $M$ is
irreducible, the closure of each component of the complement of
$P$ has connected boundary (a single $2$--sphere).

So, if $X$ is a component of the complement of $P$, then there are
no vertices of $\T$ in $X$ and as above, we can crush the
triangulation along the $2$--sphere in $\bdy X$ or we have that
$M$ is $\rp$, or $M$ is $L(3,1)$ or there is a punctured $3$--cell
$P_X$ embedded in $X$, each component of $\bdy P_X$ is a normal
$2$--sphere and the closure of each component of the complement of
$P_X$ in $X$ has connected boundary. In the last case, it is
possible that $P_X = X$ and so $X$ is a $3$--cell. If $P_X\ne X$,
we add $P_X$ to $P$ getting a larger punctured $3$--sphere, which
contains all vertices of $\T$ and the closure of each component of
its complement has connected boundary or we have that $M$ is
$S^3$. We use this larger punctured $3$--sphere and observe that
the number of times we encounter this possibility is limited by
Kneser's Finiteness Theorem. Hence, we have that $M$ is $S^3, \rp,
L(3,1)$ or we can crush the triangulation $\T$ along each
$2$--sphere in $\bdy P$. Having crushed the triangulation along
each of the $2$--spheres in $\bdy P$,  we have a number of
$3$--manifolds each has a one-vertex triangulation with fewer
tetrahedra than $\T$ (in fact, the number of tetrahedra in all the
triangulations sums to some number less than the number of
tetrahedra in $\T$). Furthermore, one of these manifolds is
homeomorphic to $M$ and each of the others is a $3$--sphere.

By Proposition \ref{0-decide}, we can decide if any of these has a
$0$--efficient triangulation. Those that have $0$--efficient
triangulations, set aside. Those that do not, then we can
construct for each a non vertex-linking normal $2$--sphere. We go
through the previous argument, we used for $M$. Since the number
of tetrahedra in each new triangulation we consider is
monotonically decreasing, having begun with the number of
tetrahedra in $\T$, we eventually have (and in practice quite
quickly) either that $M$ is one of $S^3, \rp$ or $L(3,1)$ or we
have a finite collection of $3$--manifolds, each with a
$0$--efficient triangulation and one is $M$ and the others are
$S^3$.

Now, we run the $3$--sphere recognition algorithm on each of these
$3$--manifolds with $0$--efficient triangulation. \end{proof}

\subsection{$\mathbf 0$--efficient triangulations for bounded $3$--manifolds} It
is well known \cite{bin} that a compact $3$--manifold with
boundary can be triangulated with all the vertices in the
boundary. Furthermore, a triangulation with all the vertices in
the boundary can be modified by ``closing-the-book" to a
triangulation having precisely one vertex in each boundary
component, which is not a $2$--sphere, and precisely three
vertices in each boundary component, which  is a $2$--sphere (see
\cite{jac-let-rub1}). So, one can achieve a triangulation having
all vertices in the boundary and having the property that the
triangulation restricted to a boundary component is a minimal
triangulation of that boundary component. In \cite{jac-let-rub1},
it also is shown that if $M$ is a compact, orientable
$3$--manifold with boundary, then any triangulation on $\bdy M$
can be extended to a triangulation of $M$ without adding vertices.
Our methods here are different and give minimal triangulations of
boundary components while restricting normal spheres and disks in
the manifold.

 A triangulation of a  compact
$3$--manifold with non-empty boundary is {\it $0$--efficient} if
and only if all normal disks are vertex-linking (recall that a
normal surface with  boundary must be properly embedded). It turns
out that certain technical nuisances appear if our manifolds have
boundary components which are $2$--spheres. So, we will restrict
our investigations to $3$--manifolds with boundary, where no
boundary component is a $2$--sphere. In this subsection we also
need some terminology for specially embedded punctured $3$--cells.
We say that a punctured $3$--cell embedded in a $3$--manifold $M$
so that it meets $\bdy M$ in a connected (planar) subset of its
boundary is a {\it relative} punctured $3$--cell in $M$. The
frontier of a relative punctured $3$--cell consists of properly
embedded disks and $2$--spheres.

Our first result is analogous to Proposition \ref{p-0-eff}.

\begin{prop}\label{p-b-0-eff} \cite{jac-let-rub1} Suppose $M$ is
a  compact, orientable $3$--manifold with nonempty boundary, no
component  of which is a $2$--sphere. If $M$ has a $0$--efficient
triangulation, then there are no normal $2$--spheres and $M$ is
irreducible and $\bdy$--irreducible. Furthermore, the
triangulation has all its vertices in $\bdy M$ and has precisely
one vertex in each boundary component.
\end{prop}

\begin{proof} Suppose $\T$ is a $0$--efficient triangulation of $M$.
First, we show there are no normal $2$--spheres in $M$. To this
end, suppose $S$ is a normal $2$--sphere. Let $X$ be the closure
of a component of the complement of $S$ in $M$, which also
contains a component of $\bdy M$. Let $\C$ be the induced
cell-decomposition on $X$. There is an arc, say $\Lambda$, in the
$1$--skeleton of $\C$ so that $\Lambda$ meets $\bdy M$ in a vertex
of $\C$ (also a vertex of $\T$) and meets $S$ in precisely one
point (a vertex in the induced cell structure on $S$). The
frontier of a small regular neighborhood of $S\cup\Lambda$
consists of a copy $S'$ of $S$ and a properly embedded disk $E$,
both of which are barrier surfaces in the components of their
complements not containing $S\cup\Lambda$. Notice that $E\cup S'$
is the frontier of a relative punctured $3$--cell containing
$S\cup\Lambda$ and meeting $\bdy M$ in a disk $E'$ where $\bdy E'
= \bdy E$. We proceed by shrinking $E$ in the component of its
complement not meeting $S\cup\Lambda$. We have a new phenomenon
that we did not see above; namely, in the shrinking of $E$, we may
need to make $\bdy$--compressions, as well as compressions and
isotopies. A $\bdy$--compression may be to the inside of our
relative punctured $3$--cell, in which case it is clear that we
split the relative punctured $3$--cell into two relative punctured
$3$--cells with one containing $S\cup\Lambda$. On the other hand,
we may have a $\bdy$--compression to the exterior of a relative
punctured $3$--cell; however, the result of this is still a
relative punctured $3$--cell (see Figure \ref {f-rel-punct-cell}).
As a result of shrinking $E$, we have a number of relative
punctured $3$--cells, one containing $S\cup \Lambda$, and the
components of their frontiers are either normal $2$--spheres,
normal disks, or $0$--weight $2$--spheres or disks properly
embedded in some of the tetrahedra of $\T$. We may ignore all of
the relative punctured $3$--cells except that one containing
$S\cup\Lambda$, say $P$. Furthermore, we may fill in any
$0$--weight components of the frontier with $3$--cells, missing
$S\cup\Lambda$. Since $S\cup\Lambda\subset P$ and $\Lambda$ meets
$\bdy M$, any normal disk (a vertex--linking normal disk) in the
frontier of $P$ is parallel into $\bdy M$ away from
$\Lambda\cap\bdy M$. However, this is only possible if a component
of $\bdy M$ is a $2$--sphere, a contradiction. So, this is where
we use that there are no $2$--spheres in $\bdy M$ and conclude
that there are no normal $2$--spheres in $M$.

If $M$ were not irreducible, then by Theorem \ref{normalsphere}
there would necessarily be an essential, embedded $2$--sphere and
thus an essential embedded, normal $2$--sphere; so, $M$ is
irreducible. Similarly, if there were an essential, properly
embedded disk, then by Theorem \ref{normaldisk}, there would be an
essential, embedded  normal disk. But such a disk is not
vertex-linking; so $M$ is $\bdy$--irreducible.

\begin{figure}[htbp]

            \psfrag{d}{$\subset\bdy M$}
            \psfrag{a}{\footnotesize $\alpha$}
            \psfrag{b}{\footnotesize$\beta$}
            \psfrag{s}{\begin{tabular}{c}
            split\\
    ($\bdy$--compression)\\
            \end{tabular}}
             \psfrag{t}{\begin{tabular}{c}
            add a relative\\
      $2$--handle\\
            \end{tabular}}
        \vspace{0 in}
        \begin{center}
        \epsfxsize=3.5 in
        \epsfbox{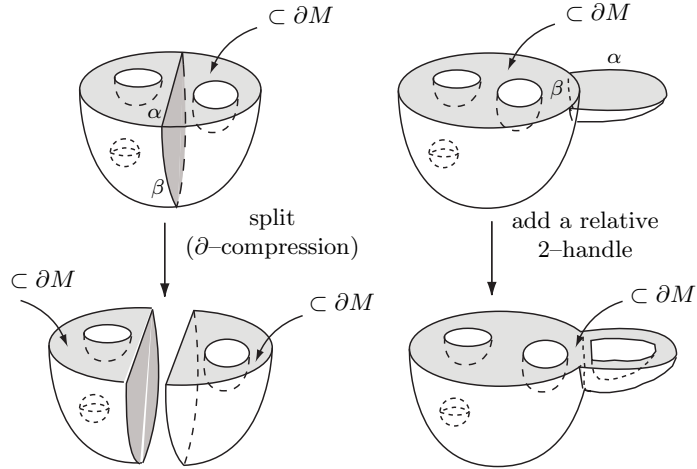}
        \caption{A relative punctured $3$--cell.}
        \label{f-rel-punct-cell}
        \end{center}

\end{figure}

Now, a vertex in $\stackrel{\circ}{M}$, the interior of $M$, would
have a normal vertex-linking $2$--sphere; so, again it follows
that all vertices are in $\bdy M$.

Suppose there is more than one vertex in a component of $\bdy M$.
Then there is an edge, $e$, in a component of $\bdy M$, which is
also an arc (has distinct vertices). The frontier of a small
regular neighborhood of $e$ is a properly embedded disk and forms
a barrier surface in the component of its complement not meeting
$e$; hence, it can be shrunk in the complement of $e$. As above,
we arrive at a relative punctured $3$--cell containing $e$, where
each component of its frontier is either a normal $2$--sphere, a
normal disk, or a $0$--weight $2$--sphere or disk properly
embedded in  a tetrahedron of $\T$. But then we conclude that not
only would $\bdy M$ be a $2$--sphere but $M$ would itself be a
$3$--cell. This completes the proof.\end{proof}

In the preceding, we considered $3$--manifolds having no
$2$--sphere boundary components. A $2$--sphere in the boundary of
a $3$--manifold is a nuisance that just doesn't fit well with
$0$--efficient triangulations.  On the other hand, one might try
to avoid the nuisance of $2$--spheres in the boundary by
considering triangulations of $3$--manifolds with nonempty
boundary, where the only normal disks are vertex-linking and there
are no normal $2$--spheres.  It seems likely that if a $3$--cell
has a $0$--efficient triangulation and has no normal $2$--spheres,
then the situation is extremely limited, as the only possibility
is the single tetrahedron triangulation of the $3$--cell with
three vertices (see Figure \ref{f-tetra}(2)). However, there is a
curiosity associated with triangulations of the $3$--cell, which
we have not resolved and seems to exhibit features similar to the
problem posed above about the number of $2$--sphere boundary
components of the closure of the complement of a maximal
collection of normal $2$--spheres. Namely, if we assume the only
normal disks are vertex-linking, then all vertices must be in the
boundary and there must be at least three vertices in the
boundary. Furthermore, there can not be an edge in the boundary
which is disjoint from any other edge. It follows that every edge
in the boundary must have a vertex in common; and, if there are
more than three vertices in the boundary, then every edge in the
triangulation must have a vertex in common. We have examples,
constructed by Ben Burton, of triangulations of a $3$--cell so
that every edge has a vertex in common and there are more than
three vertices, all in the boundary; however,  every one we know
is not $0$--efficient but has no normal $2$--spheres. An answer to
this problem is not important to this work, but it is a curiosity
we would like to resolve.

In this section, as in the preceding work, in order to obtain
$0$--efficient triangulations, we ``crush" a given triangulation
along an appropriate normal surface. Here, in the bounded case, we
need a relative version of our concept of crushing used above; for
this, the normal surface along which we crush has each component
an embedded normal disk.

Suppose $\T$ is a triangulation of a compact, orientable,
irreducible $3$--manifold with nonempty boundary. Suppose
$\mathcal{E}$ is pairwise disjoint collection of embedded normal
disks in $M$ and $X$ is the closure of a component of the
complement of the disks in $\mathcal{E}$; furthermore, suppose $X$
does not contain any vertices of $\T$. In this situation, we give
sufficient conditions for constructing a triangulation $\T'$ of
$X$, having, in particular, the property that the number of
tetrahedra in $\T'$ is strictly less than the number of tetrahedra
in $\T$. The ideas are modelled on those used in Section
\ref{sect-crush}. We will assume the reader is now familiar with
our techniques and so, we will be brief in our descriptions at
this point. However, we will organize our construction into a
theorem at the end of our discussion.

Just as above, the manifold $X$ has a nicely described
cell-decomposition, say $\mathcal{C}$, induced from the
triangulation $\T$. The induced cells are truncated tetrahedra,
truncated prisms,  triangular parallel regions, and quadrilateral
parallel regions. See Figure \ref{f-cell-decomp}. We will now show
that under the right conditions, one can construct a triangulation
of $X$ by essentially  replacing the truncated tetrahedra in $\C$
by tetrahedra.

We follow the notation used above. Let $\bbb{P}(\C) = \{$edges of
$\C$ not in the disks in $\mathcal{E}\}\cup\{$cells of type III
and type IV in $\C\} \cup\{$all trapezoidal faces of $\C\}$. Each
component of $\bbb{P}(\C)$ is an $I$-bundle. As above, we suppose
$\bbb{P}(\C) \ne X$ and all of the components of $\bbb{P}(\C)$ are
product $I$-bundles; of course, just as above, we will need  to
confirm this in practice. We continue to use $K_i\times[0,1], 1\le
i\le k,$ to denote the components of $\bbb{P}(\C)$, where
$K_{i}^\ve = K_i\times \ve, \ve = 0,1$ are isomorphic subcomplexes
of the induced normal cell structure on the disks in
$\mathcal{E}$. Notice that the vertices of all edges of $\C$ not
in $\mathcal{E}$ are themselves in $\mathcal{E}$; similarly, the
boundaries of the trapezoidal regions are in $\mathcal{E}$. While
the techniques above (in the case of crushing a triangulation
along a closed normal surface) work perfectly well for when the
surface is not connected, we applied the method in cases where the
normal surface was connected. We remark that now it is most likely
that the collection of disks $\mathcal{E}$ has a number of
components; in practice, there is one disk for each boundary
component of the $3$--manifold $M$. We assume there is a pairwise
disjoint collection of contractible planar complexes $D_j, 1\leq
j\leq k'\le k$, along with pairwise disjoint embeddings $D_j\times
[0,1]$ into $X$ so that $D_j\times [0,1]$ is a subcomplex of $X$
and for $D_{i}^\ve = D_i\times \ve, \ve =0,1$, we have $D_{i}^\ve$
is embedded as a subcomplex of the induced cell structure on the
disks in $\mathcal{E}$;  and, $\bigcup (K_i\times [0,1]) \subset
\bigcup (D_j\times [0,1])$ and the frontier of $\bigcup (D_j\times
[0,1])$ is contained in the frontier of $\bigcup (K_i\times
[0,1])$. As above, we have $K_{i}^0$ is isomorphic to $K_{i}^1$
for every $i$, but it is not necessarily the case that $D_{j}^0$
is isomorphic to $D_{j}^1$ for every $j$. See Figure
\ref{f-product-new} keeping in mind now that instead of a single
surface $S$, we might have a number of surfaces $E_1,\ldots, E_n$.


We let $\bbb{P}(X)$ be the collection of components  of $\bigcup
(D_j\times [0,1])$ and again call  $\bbb{P}(X)$ the induced
product region for $X$. Here we have that the induced  product
region $\bbb{P}(X)$ for $X$ is a trivial product region for $X$.

Essentially, the induced product region is just as in the closed
case but since the components of $\mathcal{E}$ are all disks, we
have that $\bbb{P}(X)$ is a trivial product region. On the other
hand, we do have a new twist for chains of truncated prisms.
Recall that a cell of $X$ of type $II$, a truncated prism, has two
hexagonal faces. For closed $3$--manifolds, above, these hexagonal
faces in the cell decomposition of $X$ can be identified to either
a hexagonal face of a cell of type $I$ (truncated tetrahedron) or
a hexagonal face of another cell of type $II$. This is also true
here but there is another possibility; namely, a hexagonal face
may have no identifications because it is in $\bdy M$.

If we follow a sequence of such identifications through cells of
type $II$,  we again trace out a well-defined arc which now can
terminate either at the identification of a hexagonal face of a
cell of type $II$  with one of type $I$ or in the boundary of $M$.
As above, it is also possible that these identifications form a
simple closed curve through cells of type $II$, forming a complete
cycle. We continue to call the collection of truncated prisms
identified in this way a {\it chain}. If at some point the
identification of the hexagonal faces end in a truncated
tetrahedron or in $\bdy M$, we say the chain {\it terminates}.
Note, a chain may terminate with truncated tetrahedrons at both
ends or a truncated tetrahedron only at one end and the other in
$\bdy M$ or with both ends in $\bdy M$.  If the chain begins and
ends in $\bdy M$, we call it a {\it relative cycle} and, as above,
we call a complete chain a {\it cycle} (see Figure \ref{f-cycle}
above and Figure \ref{f-bchain} below).

\begin{figure}[htbp]

            \psfrag{A}{$\sigma\subset\bdy M$}
            \psfrag{B}{$\sigma'\subset\bdy M$}
            \psfrag{e}{$e$}
            \psfrag{f}{$e'$}
            \psfrag{g}{$e''$}
        \vspace{0 in}
        \begin{center}
        \epsfxsize=3 in
        \epsfbox{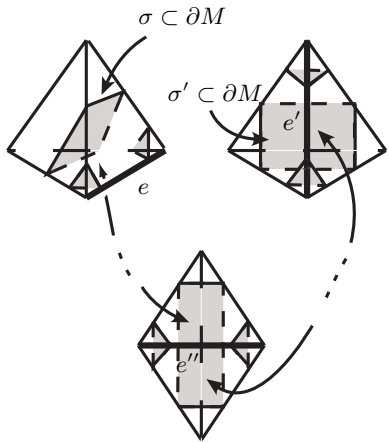}
        \caption{Relative cycles of truncated prisms - boundary case.}
        \label{f-bchain}
        \end{center}

\end{figure}

We will assume there are no cycles or relative cycles of truncated
prisms in the induced cell structure on $X$, which are not in
$\bbb{P}(X)$. Later, of course, we will have the burden of proof
to establish this in applications of crushing along a collection
of embedded normal disks.

Again, if $\{\overline{\Delta}_1,\ldots,\overline{\Delta}_n\}$ is
the collection of truncated tetrahedra in $X$, which are not in
$\bbb{P}(X)$, and $\overline{\sigma}_i$ is a face of
$\overline{\Delta}_i$, then there are four possibilities for
$\overline\sigma_i$. We have  $\overline{\sigma}_i$ is identified
with the face $\overline{\sigma}_j$ of $\overline{\Delta}_j$; we
have $\overline{\sigma}_i$ is in $\bdy M$ and is not identified
with any other face; we have that $\overline{\sigma}_i$ is at one
end of a chain of truncated prisms and the face
$\overline{\sigma}_j$ of $\overline{\Delta}_j$ is at the other end
and so there is an induced identification of $\overline{\sigma}_i$
with $\overline{\sigma}_j$ through this chain of truncated prisms;
and, lastly, we have $\overline{\sigma}_i$ is at one end of a
chain of truncated prisms and the face $\overline{\sigma}_j$ is at
the other end of this chain and also in $\bdy M$ and so there is
no identification of $\overline{\sigma}_i$. So, as before, the
faces of the truncated tetrahedra in
$\{\overline{\Delta}_1,\ldots,\overline{\Delta}_n\}$ have an
induced pairing but now, not only is it the case that faces which
are in $\bdy M$ don't have any pairings, faces which are at the
end of a chain of truncated prisms which start in $\bdy M$ do not
have a pairing. Each truncated tetrahedron in $X$ has its
triangular faces in $S$. We can identify each such triangular face
to a point (distinct points for each triangular face) and we get
tetrahedra. We will now use the notation $\td{\Delta}_i'$ for the
tetrahedron coming from the truncated tetrahedron
$\overline{\Delta}_i$ by identifying each of the triangular faces
of $\overline{\Delta}_i$ to a point (distinct points for each
triangular face) and use $\td\sigma_i$ for the triangle coming
from the hexagonal face $\overline{\sigma}_i$. Then
$\boldsymbol{\Delta}'(X) =\{\td{\Delta}'_1,\ldots,
\td{\Delta}'_n\}$ is a collection of tetrahedra with orientation
induced by that on $\T$ and the induced pairings described above
is a family $\boldsymbol{\Phi}'$ of orientation reversing affine
isomorphisms. Hence, we get a $3$--complex
$\boldsymbol{\Delta}'(X)/\boldsymbol{\Phi}'$, which has its
underlying point set a $3$--manifold, which is homeomorphic to
$X$. We will denote the associated triangulation by $\T'$.  We
call $\T'$ the triangulation obtained by {\it crushing the
triangulation $\T$ along $\mathcal{E}$}. We denote the image of a
tetrahedron $\td{\Delta}'_i$ by $\Delta'_i$ and, as above, call
$\td{\Delta}'_i$ the lift of $\Delta'_i$; we denote the image of
the face $\td\sigma_i$ of $\td{\Delta}'_i$ by $\sigma_i$. (See
Figure \ref{f-ident-prisms}, replacing in that Figure the
$\td{\Delta}_k^*,  k = i,j$ by $\td{\Delta}_k',  k = i,j$,
respectively.)

We summarize in the following theorem.

\begin{thm}\label{b-crush} Suppose $\T$ is a triangulation of  the
compact, orientable, irreducible $3$--manifold $M$ with nonempty
boundary. Suppose $\mathcal{E}$ is a pairwise disjoint collection
of embedded normal disks in $M$, $X$ is the closure of a component
of the complement of $\mathcal{E}$ and $X$
 does not contain any vertices of $\T$. Suppose there is an induced
 product region, $\bbb{P}(X)$, for $X$. If
\begin{enumerate}
\item[i)] $X\ne \bbb{P}(X)$, \item[ii)] $\bbb{P}(X)$ is a trivial
product region for $X$, and \item[iii)] there are no cycles or
relative cycles of truncated prisms in $X$, which are not in
$\bbb{P}(X)$,
\end{enumerate}
then the triangulation $\T$ can be crushed along $\mathcal{E}$ and
$\T'$ is a triangulation of $X$.
\end{thm}

If the disks in $\mathcal{E}$ are not all vertex-linking and $t$,
$t'$ are the number of tetrahedra of $\T$ and $\T'$, respectively,
then $t > t'$; in fact, $\T' = \T$ if an only if all the disks in
$\mathcal{E}$ are vertex-linking.

The following is an existence theorem and is the version of
Theorem \ref{0-eff-exists} for bounded $3$--manifolds. We have
organized this study of $0$--efficient triangulations for bounded
$3$--manifolds analogously to that above for closed
$3$--manifolds; however, we point out here and below that if one
is given a compact, orientable $3$--manifold with boundary, which
is known to be irreducible and $\bdy$--irreducible, then any
triangulation can be modified to a $0$--efficient triangulation or
the $3$--manifold is the $3$--cell and (here is the major point)
we do {\it not} need to employ the $3$--sphere recognition
algorithm. See Theorem \ref{0-eff-construct} above and Theorem
\ref{b-0-eff-construct} below.

\begin{thm}\label{b-0-eff-exists} A compact, orientable,
irreducible, $\bdy$--irreducible
$3$--manifold with nonempty boundary has  a $0$--efficient triangulation.
\end{thm}

\begin{proof} Suppose $M$ is a compact, orientable, irreducible
and $\bdy$--irreducible $3$--manifold
with nonempty boundary. Let $\T$ be a triangulation of $M$. We
wish to mimic, with appropriate modification, the proof of Theorem
\ref{0-eff-exists}.

Let $S_1,...,S_n$ denote the boundary components of $M$. Let
$\T_i$ denote the triangulation of $S_i$ induced by $\T$. If a
component of the boundary is a $2$--sphere, then since $M$ is
irreducible, $M$ is a $3$--cell. If $M$ is the $3$--cell, we can
give $M$ a one tetrahedron, three vertex triangulation, which is a
$0$--efficient triangulation of a $3$--cell.  Hence, we will
assume $M$ is not a $3$--cell and, so, no component of $\bdy M$ is
a $2$--sphere.

If there is a normal $2$--sphere, there is a normal disk which is not
vertex--linking or, following the arguments from above, $\bdy M$ would be a
$2$--sphere. So, in our situation, to show $\T$ is $0$--efficient is completely
dependent on whether every normal disk is vertex-linking.

If every normal disk is vertex-linking, we are done. So, we may
assume there is a normal disk, which is not vertex-linking. Notice
that if $E$ is any properly embedded disk, then there is a disk
$E'\subset \bdy M$, $E'\cap E = \bdy E = \bdy E'$ and $E\cup E'$
bounds a $3$--cell in $M$. Furthermore, if $E$ is a normal non
vertex-linking disk and $F$  is a properly embedded normal disk so
that $E$ is contained in the $3$--cell that $F$ co-bounds with a
disk $F'\subset \bdy M$, then $F$ is not vertex-linking ($M$ is
not a $3$--cell). So, under the assumption that there is a non
vertex-linking disk, we consider a collection $\mathcal{E}$ of
pairwise disjoint, properly embedded, normal disks,
$E_1,\ldots,E_n$ where $\bdy E_i\subset S_i, i = 1,\ldots,n$ and
$E_1$, say, is not vertex-linking; furthermore, we may take the
disks in this collection {\it maximal} in the sense that  if $F$
is a properly embedded, normal disk in $M$, $F$ is disjoint from
all the disks in $\mathcal{E}$, $F$ co-bounds a $3$--cell $B$ with
a disk in $\bdy M$ and some $E_i\subset B$, then $F = E_i$.

For each disk $E_i\in\mathcal{E}$, let $E_i'$ denote the disk in
$S_i$ so  $E_i\cap E_i' = \bdy E_i = \bdy E_i'$ and $E_i\cup E_i'$
bounds a $3$--cell. Denote these $3$--cells bounded by $E_i\cup
E_i'$, $B_i, i = 1,\ldots,n$.

We claim that all the vertices of $\T$ are in $\bigcup_{i=1}^n
B_i$. For suppose there is a vertex $v$ of $\T$ not in
$\bigcup_{i=1}^n B_i$. Then there is an arc $\Lambda$ in the
$1$--skeleton of $\T$, $\Lambda$ has one end at $v$ and meets
$\bigcup_{i=1}^n B_i$ in a single point in, say $E_j$, for some
$j, 1\leq j \leq n$. Furthermore, if $v$ is in $\T_j$ we can take
$\Lambda$ in the $1$-skeleton of $\T_j$.  We consider a small
regular neighborhood, $N$, of $ B_j \cup\Lambda$. The frontier of
$N$ consists of a properly embedded disk $F_j$, which may not be
normal. However, the frontier of $N$ along with the collection of
normal disks $E_i, i\ne j$ is a barrier surface in the component
of their complement not meeting $\bigcup_{i=1}^n B_i \cup\Lambda$.
Furthermore, there is a disk $F_j'\subset S_j, \bdy F_j' = \bdy
F_j$ and $F_j'\cup F_j$ bounds a relative $3$--cell in $M$. Note
that $F_j' = E_j'$ (normally isotopic with the same normal curve
as a boundary) if $\Lambda$ is not in $S_j$ and, in general,
$F_j'\supset E_j'$. We shrink $F_j$, using Theorem \ref{barrier},
obtaining a relative punctured $3$--cell, say $P$, embedded in
$M$, which contains $B_j\cup\Lambda$, $P\cap \bdy M\subset S_j$,
$P\cap E_i = \emptyset$ and each component of the frontier of $P$
in $M$ is either a normal disk, a normal $2$--sphere or a
$0$--weight $2$--sphere or disk properly embedded in a tetrahedron
of $\T$. Each $0$--weight $2$--sphere in the frontier of $P$
bounds a $3$--cell, whose interior misses $P$; similarly, for each
$0$--weight disk in the frontier of $P$ there is a disk in $\bdy
M$, which together with the disk in the frontier of $P$ bounds a
$3$--cell whose interior misses $P$. We fill in these $0$-weight
frontiers with these $3$--cells. We still have a relative
punctured $3$--cell, which we continue to call $P$,
$B_j\cup\Lambda\subset P$, $P\cap\bdy M\subset S_j$, and $P\cap
E_i = \emptyset, i\ne j$. If there is a normal $2$--sphere in
$\bdy P$, then, since $M$ is irreducible, such a $2$--sphere
bounds a $3$--cell in $M$, which can not meet $\bdy M$. We add all
such $3$--cells to $P$. We get a relative punctured $3$--cell, we
continue to call it $P$,  $B_j\cup\Lambda\subset P$, $P\cap\bdy
M\subset S_j$, and each component of the frontier of $P$ is a
normal disk. Now, suppose $F$ is a normal disk in the frontier of
$P$. Then there is a disk $F'\subset S_j, \bdy F = \bdy F'$ and
$F\cup F'$ bounds a $3$--cell in $M$. If this $3$--cell contains
$P$ and thus $E_j$, then we contradict the maximality of the disks
in the collection $\mathcal{E}$; so, we may assume that the
interior of the $3$--cell is disjoint from $P$. However, if this
happens for every normal disk in the frontier of $P$, then $S_j$
would be a $2$--sphere, which leads to a contradiction.  Thus
assuming that a vertex of $\T$ is not in $\bigcup_{i=1}^n B_i$
leads to a contradiction.

Let $X$ be the closure of the component of $M\setminus \bigcup_{i
= 1}^n B_i$ that does not contain the vertices of $\T$. First we
notice that $X$ is homeomorphic with $M$. As in the proof of
Theorem \ref{0-eff-exists}, the manifold $X$ has a nicely
described cell decomposition induced from the triangulation $\T$
and the way in which the normal disks $E_1,\ldots,E_n$ sit in $\T$
(all the vertices are in the $3$--cells complementary to $X$).

The cells of $X$ are the same four types as above:  truncated
tetrahedra, truncated prisms,  triangular parallel regions, and
quadrilateral parallel regions (see Figure \ref{f-cell-decomp}),
except that some of the truncated  faces (those faces also in
faces of tetrahedra of $\T$) may now be in $\bdy M$, including
some trapezoidal faces.

We need to show the conditions of Theorem \ref{b-crush} are
satisfied and hence, we can crush the triangulation $\T$ along the
collection of normal disks $\mathcal{E}$.

We let $\C$ denote the induced cell-decomposition of $X$ and, as
above, let $\bbb{P}(\C) = \{$edges of $\C$ not in the disks in
$\mathcal{E}\}\cup\{$cells of type III and type IV in $\C\}
\cup\{$all trapezoidal faces of $\C\}$. Each component of
$\bbb{P}(\C)$ is an $I$-bundle.

\vspace{.15 in} \noindent {\it Claim. $\bbb{P}(\C)\ne X$ and each
component of $\bbb{P}(\C)$ is a product $I$--bundle.}

\vspace{.1 in} \noindent {\it Proof of Claim}. If $\bbb{P}(\C) =
X$, then $X$ is an I-bundle with corresponding $\bdy I$--bundle a
collection  of disks. Hence, $X$ is the product of a disk and an
interval; but then $M$ would be the $3$--cell, a contradiction.

If some component of $\bbb{P}(\C)$ is not a product $I$--bundle,
then there is a M\"obius band properly embedded in $X$ with its
boundary in one of the disks $D_i$. However, then we would have a
projective plane embedded in $M$; this is impossible since we have
assumed $\bdy M\ne \emptyset$ and $M$ is irreducible. This
completes the proof of the Claim.

\vspace{.15 in} \noindent {\it Claim. There is a trivial induced
product region $\bbb{P}(X)$ for $X$. }

\vspace{.1 in} \noindent {\it Proof of Claim}. Essentially, we
follow the argument as in the analogous Claim in the proof of
Theorem \ref{0-eff-exists}. However, here in the bounded case, we
have two additional and new phenomena  from what we had above, in
the case of $M$ closed.

We have established that each component of $\bbb{P}(\C)$ is a
product $I$--bundle. As above, we write each I-bundle component of
$\bbb{P}(\C)$ as $K_i\times I, i = 1,\ldots,k$, where $k$ is the
number of components of $\bbb{P}(\C)$ and again, we set $K_{i}^\ve
= K_i\times \ve, \ve = 0, 1$. The subcomplexes $K_i^0$ and   $K_i^
1$ are isomorphic subcomplexes in the induced normal cell
structures on the disks in our collection $\mathcal{E}$. Here, we
may have $K_{i}^0$ and $K_{i}^1$ in the same or in distinct disks
in $\mathcal{E}$.

We would like to have for each of the  I-bundles $K_i\times [0,1]$
that $K_i$ is simply connected; but just as above, this is not
necessarily the case. So, let's consider such a $K_i\times [0,1]$.
We may as well assume that $K_i$ is not simply connected. Thus if
$K_{i}^\ve$ is in the disk $E_j$ in the collection $\mathcal{E}$,
then $K_i^\ve$ separates $E_j$; and if $K_i^\ve$ meets $\bdy E_j$,
then every component of the complement of $K_i^\ve$ is simply
connected. This is just as in the closed case; however, if
$K_i^\ve$ does not meet $\bdy E_j$, then the component of the
complement of $K_i^\ve$ in $E_j$, which contains $\bdy E_j$, is
not simply connected. Of course, all the other components of the
complement of $K_i^\ve$ in $E_j$ are simply connected. Because of
the boundary of $E_j$, we have a new consideration.

Suppose $K_i$ is not simply connected, $K_i^\ve\subset E_j, \ve =
0,1$ and $K_i\times [0,1]$ does not meet $\bdy E_j$. We claim that
by having chosen the collection $\mathcal{E}$ maximal, we must
have that $K_i^0$ is in the same component of the complement of
$K_i^1$ as $\bdy E_j$ and $K_i^1$ is in the same component of the
complement of $K_i^0$ as $\bdy E_j$. For suppose $K_i^1$ is not in
the component of the complement of $K_i^0$ that meets $\bdy E_j$.
Define $D_i^0$ to be $K_i^0$ along with all the components of its
complement not meeting $\bdy E_j$. Then we have that $D_i^0$ is
simply connected and also we have $K_i^1\subset D_i^0$. See Figure
\ref{f-ideal-pushout}. Let $N_i = N(D_i^0\cup (K_i\times [0,1]))$
be a small regular neighborhood of $D_i^0\cup (K_i\times [0,1])$.
Since $K_i^1\subset D_i^0$, the frontier of $N_i$ consists of an
annulus, possibly some $2$--spheres and a disk $F_i$, properly
embedded in $X$ and having its boundary in $E_j$.  There is a disk
$F_i'$ in $E_j$ so that $\bdy F_i' = \bdy F_i$ and $F_i\cup F_i'$
bounds a $3$--cell, say $B$, in $X$; actually, we want to think of
$B$ as a relative (punctured) $3$--cell. Furthermore, $B$ contains
$N_i$ and $F_i$ along with $\cup E_i, i\ne j$, form a barrier
surface in the component of the complement of $F_i$ not meeting
$N_i$. We shrink $F_i$ in the component of the complement of $F_i$
not meeting $N_i$. Since $M$ (and hence, $X$) is not a $3$--cell,
our standard arguments give that in this situation we would have a
 {\it normal} disk $E$ in $X$ and a disk $E'$ in
$\bdy X$ so that $\bdy E' = \bdy E$ and $E\cup E'$ bounds a
$3$--cell containing $B$. But this is a contradiction to our
choice of  the collection $\mathcal{E}$ being maximal.

\begin{figure}[htbp]

            \psfrag{a}{$K_i\times [0,1]$}
            \psfrag{b}{$B$}
            \psfrag{c}{$F_i$}
            \psfrag{d}{$F_i'$}
            \psfrag{e}{\Large $E_j$}
            \psfrag{0}{$D_i^0$}
            \psfrag{1}{$D_i^1$}
            \psfrag{A}{$D_i^0\subset D_i^1\subset E_j$}
        \vspace{0 in}
        \begin{center}
        \epsfxsize=3 in
        \epsfbox{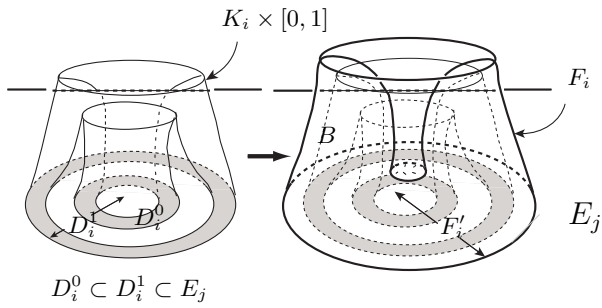}
        \caption{$K_i^1$ is not in the same component of the complement
of $K_i^0$ as $\bdy E_j$.}
        \label{f-ideal-pushout}
        \end{center}

\end{figure}

So, having made this observation, let's consider the possibilities
for $K_i\times [0,1]$. One possibility is that $K_i^0$ and $K_i^1$
are in distinct disks of our collection $\mathcal{E}$, say
$K_i^0\subset E_j$ and $K_i^1\subset E_{j'}, j\ne j'$. In this
case, $K_i^0$ does not meet $\bdy E_j$ and $K_i^1$ does not meet
$\bdy E_{j'}$. So, we let $D_i^0$ be the union of $K_i^0$ along
with the components of the complement of $K_i^0$ in $E_j$, which
do not meet $\bdy E_j$; similarly, we let $D_i^1$ be the union of
$K_i^1$ along with the components of the complement of $K_i^1$ in
$E_{j'}$, which do not meet $\bdy E_{j'}$. Then $D_i^\ve, \ve =
0,1$ is simply connected. (If $K_i$ is simply connected, then
$K_i^\ve, \ve = 0,1$ does not separate and $D_i^\ve = K_i^\ve$.)
We let $N_i = N(D_i^0\cup (K_i\times [0,1])\cup D_i^1)$ be a small
regular neighborhood of $D_i^0\cup (K_i\times [0,1])\cup D_i^1$.
Then the frontier of $N_i$ consists of a single annulus and
possibly some $2$--spheres. Furthermore, since $M$ is irreducible
and $\bdy M\ne\emptyset$, we have immediately that any such
$2$--sphere bounds a $3$--cell not meeting $D_i^0\cup (K_i\times
[0,1])\cup D_i^1$ and so each plug for $N_i$ is a $3$--cell. It
follows that there is a simply connected planar complex $D_i$ and
an embedding of $D_i\times [0,1]$ into $X$ so that $D_i\times 0 =
D_i^0, D_i\times 1 = D_i^1, K_i\times [0,1]\subset D_i\times
[0,1]$, and the frontier of $D_i\times [0,1]$ is contained in the
frontier of $K_i\times [0,1]$. See Figure \ref{f-product-new}.

Another possibility is that  both $K_i^0$ and $K_i^1$ are in the
same disk $E_j$ and  $K_i^0$ (and hence, $K_i^1$) does not meet
$\bdy E_j$. By our observation above, we have that $K_i^0$
($K_i^1$, respectively) is in the component of the complement of
$K_i^1$ ($K_i^0$, respectively) that meets $\bdy E_j$. In this
situation, we let $D_i^0$ (respectively $D_i^1$) be the union of
$K_i^0$ (respectively $K_i^1$) along with all components of the
complement of $K_i^0$ (respectively $K_i^1$) in $E_j$ not meeting
$K_i^1$ (respectively $K_i^0$). Then $D_i^\ve, \ve = 0, 1$ is
simply connected. Here again, we get straight away that for $N_i =
N(D_i^0\cup (K_i\times [0,1])\cup D_i^1)$, a small regular
neighborhood of $D_i^0\cup (K_i\times [0,1])\cup D_i^1$, then the
plugs for $N_i$ are $3$--cells. Hence, in this case as in the
previous case,  it follows that there is a simply connected planar
complex $D_i$ and an embedding of $D_i\times [0,1]$ into $X$ so
that $D_i\times 0 = D_i^0, D_i\times 1 = D_i^1, K_i\times
[0,1]\subset D_i\times [0,1]$, and the frontier of $D_i\times
[0,1]$ is contained in the frontier of $K_i\times [0,1]$. Again,
see Figure \ref{f-product-new}.

The only remaining possibility is that both $K_i^0$ and $K_i^1$
are in the same disk $E_j$ and now  $K_i^0$ (and hence, $K_i^1$)
does meet $\bdy E_j$. As we observed above, if we let $D_i^0$
(respectively $D_i^1$) be the union of $K_i^0$ (respectively
$K_i^1$) along with all components of the complement of $K_i^0$
(respectively $K_i^1$) in $E_j$ not meeting $K_i^1$ (respectively
$K_i^0$), then $D_i^\ve, \ve = 0, 1$, is simply connected. Just as
in all the similar situations above, we let $N_i = N(D_i^0\cup
(K_i\times [0,1])\cup D_i^1)$ be a small regular neighborhood of
$D_i^0\cup (K_i\times [0,1])\cup D_i^1$. However, now we have a
new situation. The components of the frontier of $N_i$ still might
contain some $2$--spheres but now must contain properly embedded
disks, whose boundaries are in $\bdy M$ (see Figure
\ref{f-b-product}). Just as before, any $2$--sphere frontier
bounds a $3$--cell in $M$ which misses $N_i$ and we have a
$3$--cell plug. So, let's consider what happens for a component of
the frontier of $N_i$ which is a disk, say $E$. Then there is a
disk $E'$ contained in $\bdy X$, $\bdy E' = \bdy E$, and $E\cup
E'$ bounds a $3$--cell $B$ in $X$. If $B$ does not contain $N_i$
($E'$ does not contain $E_j'$), then we have a $3$--cell plug and
are satisfied. However, {\it a priori}, it might be possible that
the $3$--cell $B$ contains $N_i$. We shall show that this is
impossible, having chosen the  collection $\mathcal{E}$ to be
maximal. First, note that $E$ along with the normal disks in
$\mathcal{E}$, distinct from $E_j$, form a barrier surface in the
component of the complement of $E$ not meeting $N_i$; so, we can
shrink $E$. However, just as above, since $M$ is irreducible,
$\bdy$--irreducible and not itself a $3$--cell, this situation
leads to a normal disk, say $F$, properly embedded in $M$ and a
disk $F'\subset S_j$ so that $\bdy F' = \bdy F$ and $F\cup F'$
bounds a $3$--cell containing $E_j$. This contradicts our choice
of $\mathcal{E}$ maximal. Hence, all plugs for $N_i$ must be
$3$--cells and again we have that there is a simply connected
planar complex $D_i$ and an embedding of $D_i\times [0,1]$ into
$X$ so that $D_i\times 0 = D_i^0, D_i\times 1 = D_i^1, K_i\times
[0,1]\subset D_i\times [0,1]$, and the frontier of $D_i\times
[0,1]$ is contained in the frontier of $K_i\times [0,1]$ (see
Figure \ref{f-b-product}).

\begin{figure}[htbp]

            \psfrag{0}{$K_i^0\subset D_i^0$}
            \psfrag{1}{$K_i^1\subset D_i^1$}
            \psfrag{k}{$D_k\times [0,1]$}
            \psfrag{j}{$E_j$}
            \psfrag{J}{$E_j'\subset S_j\subset\bdy M$}
            \psfrag{F}{$F$}
            \psfrag{G}{$F'$}
            \psfrag{i}{\begin{tabular}{c}
            $K_i\times [0,1]\subset$\\
    $D_i\times [0,1]$\\
            \end{tabular}}
            \psfrag{B}{(B)}
             \psfrag{A}{(A)}
        \vspace{0 in}
        \begin{center}
        \epsfxsize=5 in
        \epsfbox{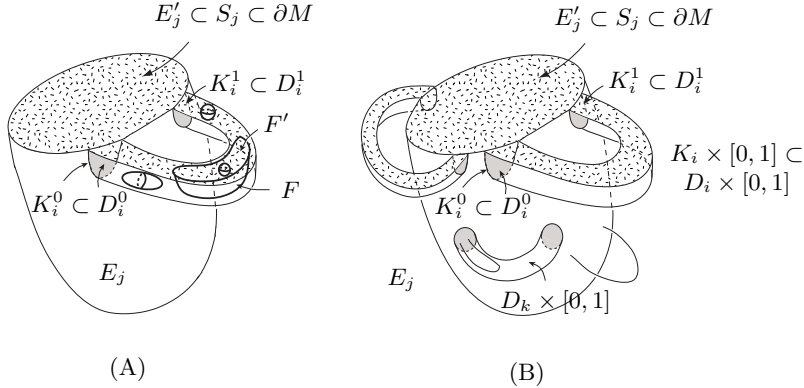}
        \caption{(A) $K_i^1$ and $K_i^0$ both meet $\bdy E_j$, leading to a
$3$--cell plug or a new normal disk $F$; (B) induced product
region, $\bbb{P}(X)$.}
        \label{f-b-product}
        \end{center}

\end{figure}

Now, having made these observations, we can choose some order for
the components of $\bbb{P}(\C)$, say $K_1\times
[0,1],\ldots,K_k\times [0,1]$, and from this, just as in the proof
of Theorem \ref{0-eff-exists}, derive that there is a trivial
induced product region $\bbb{P}(X)$ for $X$. Again, see Figure
\ref{f-b-product}.

This completes the proof of the Claim.

 \vspace{.15
in} \noindent {\it Claim. There is no  cycle (complete or
relative) of cells of type II in $X$, which are not in
$\bbb{P}(X)$.}

\vspace{.1 in} \noindent {\it Proof of Claim}. In this situation
we have four possibilities: a complete cycle about a single edge,
a complete cycle about more than one edge, a relative cycle about
a single edge and a relative cycle about more than one edge.

Suppose there were a complete cycle about a single edge $e$ as in
Figure \ref{f-cycle}, Part A. As in the similar situation in the
proof of Theorem \ref{0-eff-exists}, there would be a disk $D$
meeting $e$ in precisely one point and meeting a non
vertex-linking disk, say $E_j$, of the collection $\mathcal{E}$ in
$\bdy D$. A surgery on $E_j$ at $D$ (see Figure
\ref{f-surgery-cycle}), which is the same as adding a $2$--handle
to $E_j$ along the curve $\bdy D$, gives a {\it normal}
$2$--sphere, $S_j$, and a new {\it normal} disk, $F_j$.
Furthermore, $\bdy F_j = \bdy E_j$ and $F_j\cup E_j'$ bounds a
relative punctured $3$--cell with frontier $F_j\cup S_j$; this
relative punctured $3$--cell contains the $3$--cell bounded by
$E_j\cup E_j'$ (the disk $D$ is not in the $3$--cell bounded by
$E_j\cup E_j'$). Since $M$ is irreducible, $S_j$ bounds a
$3$--cell whose interior does not meet the $3$--cell bounded by
$E_j\cup E_j'$. But then $F_j$ co-bounds a $3$--cell with $E_j'$
and this contradicts maximality of the disks in the collection
$\mathcal{E}$. Hence, if the collection $\mathcal{E}$ is maximal,
there is no complete cycle of truncated prisms about a single
edge.

Now, suppose there is a complete cycle of truncated prisms about
more than one edge (see Figure \ref{f-cycle}(B)), then as in the
consideration for closed manifolds, the collection of cells of
type II (truncated prisms), form a solid torus with, possibly,
some self identifications in its boundary (possibly, some of the
trapezoidal faces in the boundary are identified); furthermore,
each hexagonal face of a truncated prism in the cycle of truncated
prisms is a meridional disk for the solid torus. As above, we
distinguish between the cycle of truncated prisms, which we denote
$\hat{\tau}$, and the cycle of truncated prisms minus the bands of
trapezoids, which we denote by $\tau$. We have that $\tau$ meets
the collection $\mathcal{E}$ in either three open annuli, each
meeting a meridional disk of $\tau$ precisely once, or a single
open annulus, meeting a meridional disk of $\tau$ three times. See
Figure \ref{f-cycle-torus}.

First, we consider the case we have three annuli in $\tau$, say
$A_1, A_2$ and $A_3$, each meeting a hexagonal face of a truncated
prism in $\hat{\tau}$ precisely once. The frontier of an $A_i, i =
1,2,3$ in a disk in the collection $\mathcal{E}$ is in the
collection of trapezoids in the faces of the cycle of truncated
prisms $\hat{\tau}$; thus in the induced $I$--bundle region
$\bbb{P}(\C)$.  Now, just as above, there are two components of
the frontier of each $A_i$ and each component of the frontier of
$A_i$ separates the disk of $\mathcal{E}$, which it is in.
Consider $A_1$ and suppose $A_1$ is in $E_j$ and denote the two
components of the frontier of $A_1$ by $a_1$ and $a_1'$. Then both
$a_1$ and $a_1'$ are in trapezoids in $\hat{\tau}$ and so in
$\bbb{P}(\C)$. By our above construction of the induced product
region for $X$, $\bbb{P}(X)$, it follows that $a_1$ and $a_1'$ are
each in a simply connected region of $E_j$, which does not meet
$\bdy E_j$. But this is possible only if $A_1$ is also in such a
simply connected region of $E_j$. Therefore, $A_1$ and hence
$\hat{\tau}$ is in $\bbb{P}(X)$, the induced product region for
$X$. See Figure \ref{f-annuli}.

While it can be shown directly, we have, in particular, from the
conclusion of this argument that the torus $\tau$ can meet at most
two distinct components of $\mathcal{E}$.

Suppose now we have just one annulus, say $A_1$, common to $\tau$
and the collection $\mathcal{E}$. We have that $A_1$ meets each
hexagonal face in the chain of truncated prisms, $\hat{\tau}$,
three times. We can wipe this situation away quite easily by
noting that this would imply that there is a lens space $L(3,1)$
embedded in $M$ as a connected summand but this is impossible as
$M$ is irreducible. On the other hand, the reader may wonder why
this differs from the previous case (where there were three annuli
common to $\tau$ and $\mathcal{E}$), since we only used that there
was a single annulus in the disk $E_j\in\mathcal{E}$). The same
argument could be used to conclude that $\hat{\tau}$ is in the
induced product region $\bbb{P}(X)$; however, this skirts the fact
that above we used quite strongly that $M$ is irreducible to get
the $3$--cell plugs in order to get a trivial induced product
region. This completes the proof of the Claim and  there is no
complete cycle of truncated prisms, which is not in $\bbb{P}(X)$.

Finally, we have the completely new possibility in the case when
$M$ has boundary; this is the possibility of a relative cycle of
truncated prisms. In this case, we also have two possible
situations: a relative cycle about a single edge or a relative
cycle about different edges.

In the first situation, a relative cycle about a single edge $e$,
then necessarily the edge $e$ is in $\bdy M$. There is a disk $D$
which meets the edge $e$ in a single point, meets some $E_i$ in a
spanning arc $\beta$ and meets $S_i\subset\bdy M$ in a spanning
arc $\alpha$.  See Figure \ref{f-rel-cycle}. Furthermore, $\alpha$
meets $E_i'$ only in its end points, which are in $\bdy E_i'$. The
disk $D$, which lies outside the $3$--cell $B_i$, gives a
``$\bdy$--compression" (adds a relative $2$--handle, see Figure
\ref{f-rel-punct-cell}) to the $3$--cell $B_i$ along the disk
$E_i$, resulting in a relative punctured $3$--cell with two (new)
{\it normal}, non vertex-linking disks in its frontier. Consider
the disks in the frontier of this relative punctured $3$--cell.
Since $M$ is irreducible and $\bdy$--irreducible, each co-bounds a
$3$--cell in $M$ with a disk in $S_i$ and since $S_i$ is not a
$2$--sphere at least one must co-bound a $3$--cell with a disk in
$S_i$, which also contains $B_i$. However, this contradicts the
maximality of the disks in the collection $\mathcal{E}$.

\begin{figure}[htbp]

            \psfrag{A}{$\sigma\subset\bdy M$}
            \psfrag{B}{$\sigma'\subset\bdy M$}
            \psfrag{e}{$e$}
            \psfrag{a}{$\alpha\subset\bdy M$}
            \psfrag{b}{$\beta$}
            \psfrag{d}{$e\subset\bdy M$}
            \psfrag{s}{\begin{tabular}{c}
            surgery\\
    ($\bdy$--compression)\\
            \end{tabular}}
        \vspace{0 in}
        \begin{center}
        \epsfxsize=4 in
        \epsfbox{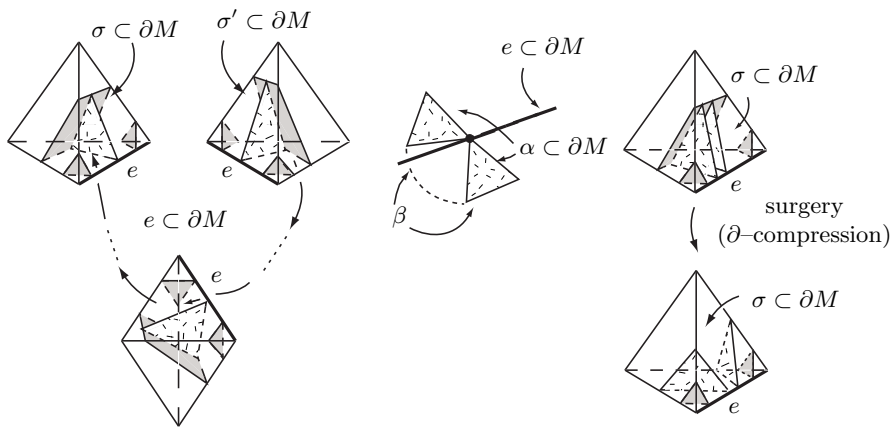}
        \caption{Relative cycle about a single edge gives a new normal
disk.}
        \label{f-rel-cycle}
        \end{center}

\end{figure}

Now, suppose we have a relative cycle about more than one edge
(see Figure \ref{f-bchain}), then the collection of cells of type
II (truncated prisms), is a $3$--cell having the form of a product
$I$-bundle, $D^2\times I$, with, possibly, some self
identifications in its vertical boundary $\bdy D^2\times I$
(possibly, some of the trapezoidal faces in the boundary are
identified). Furthermore, each hexagonal face of a truncated prism
in the relative cycle of truncated prisms is a section in the
$I$-bundle structure for the $3$--cell and there is a component of
$\bdy M$, say $S_j$, which contains $D^2\times\ve, \ve = 0,1$. As
above in the case of a complete cycle of truncated prisms, and
because of the possible singularities, we distinguish between the
relative cycle of truncated prisms, which we denote
$\hat{\delta}$, and the cycle of truncated prisms minus the bands
of trapezoids, which we denote by $\delta$. We have that $\delta$
meets the collection $\mathcal{E}$ only in $E_j$ ($\hat{\delta}$
meets only $S_j$) and then in three vertical bands of the form
$A_i = \alpha_i\times I$ for the pairwise disjoint open intervals
$\alpha_i, i = 1,2,3$ in $\bdy D^2$. See Figure
\ref{f-3strips-new}.

Now, as in the complete cycle case, the frontier of each $A_i, i =
1,2,3$ in $E_j$ is in the collection of trapezoids in the faces of
the relative cycle of truncated prisms $\hat{\delta}$; thus, the
frontier of each $A_i$ is in the induced $I$--bundle region
$\bbb{P}(\C)$. Furthermore, there are two components of the
frontier of each $A_i$ and each component of the frontier of $A_i$
separates $E_j$. Consider $A_1$ and denote the two components of
the frontier of $A_1$ by $a_1$ and $a_1'$. Then both $a_1$ and
$a_1'$ are in trapezoids in $\hat{\delta}$ and so in
$\bbb{P}(\C)$. By our construction of the induced product region
for $X$, $\bbb{P}(X)$, in the case that a component of
$\bbb{P}(\C)$ meets the boundary of a disk $E_j$, we have that
$a_1$ and $a_1'$ are each in a simply connected region of $E_j$,
which is also in $\bbb{P}(X)$. Hence, we either have $A_1$, and
therefore $\hat{\delta}$ in $\bbb{P}(X)$ or $\bbb{P}(X)$ meets
$E_j$ in the complement of $A_1$. But the latter is impossible,
since $A_2$ and $A_3$ are in regions of $E_j$ complementary to
$A_1$. So, the only possibility is that $A_1$ is in $\bbb{P}(X)$;
that is, such a relative cycle of truncated prisms, $\hat{\tau}$,
is in the induced product region for $X$, $\bbb{P}(X)$. See Figure
\ref{f-3strips-new}.

\begin{figure}[htbp]

            \psfrag{a}{$a_1$}
            \psfrag{1}{$A_1$}
            \psfrag{2}{$A_2$}
            \psfrag{3}{$A_3$}
            \psfrag{B}{$B_j$}
            \psfrag{E}{\Large{$E_j$}}
            \psfrag{b}{$a_1'$}
            \psfrag{F}{\Large{$E_j'\subset\bdy M$}}
            \psfrag{d}{$\delta$}
            \psfrag{M}{\Large $\open{M}$}
            \psfrag{D}{$D = \mathbb{D}^2$}
        \vspace{0 in}
        \begin{center}
        \epsfxsize=2.5 in
        \epsfbox{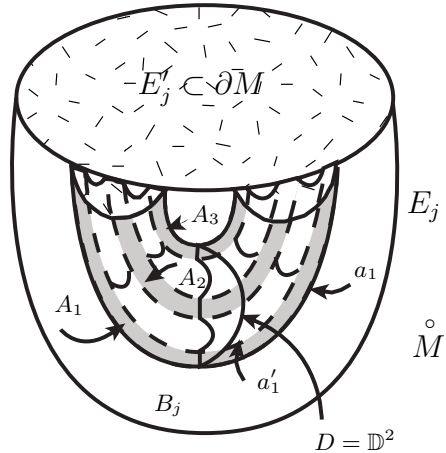}
        \caption{Relative cycle about more than one edge.}
        \label{f-3strips-new}
        \end{center}

\end{figure}

This completes the proof of the Claim.

Hence, we have established that all the conditions of Theorem \ref
{b-crush} (namely, $\bbb{P}(X)$ exists and $\bbb{P}(X)\ne X$,
$\bbb{P}(X)$ is a trivial product region for $X$, and there are no
cycles or relative cycles of truncated prisms in $X$, which are
not in $\bbb{P}(X)$). Thus by Theorem \ref{b-crush}, there is a
triangulation $\T'$ of $X$, where the tetrahedra of $\T'$ are in
one-one correspondence with the truncated tetrahedra of $\C$ not
in $\bbb{P}(X)$. Since, we had that at least one of the disks in
the collection $\mathcal{E}$ is not vertex-linking, then the
number of tetrahedra of $\T'$ is smaller than the number of
tetrahedra in the original triangulation $\T$. By repeated
applications of the above techniques, it follows that we must
eventually get the desired triangulation of $M$.\end{proof}

As we pointed out at the beginning of this section, there are
various  ways to prove that a compact, orientable $3$--manifold
with nonempty boundary, no component of which is a $2$--sphere,
admits a triangulation having all vertices in the boundary and
precisely one vertex in each boundary component. The methods of
the previous theorem give a ``natural method of identifying
maximal trees to points" to achieve such triangulations for
bounded manifolds, which are also irreducible and
$\bdy$--irreducible. We have the following corollary to Theorem
\ref{b-0-eff-exists} and Proposition \ref{p-b-0-eff}. After the
Corollary, we give an algorithm for constructing a $0$--efficient
triangulation, which also provides a construction of
triangulations having all vertices in the boundary and just one
vertex in each boundary component.

\begin{cor}\label{b-one-vertex} Suppose $M$ is a compact, orientable,
irreducible, $\bdy$--irreducible $3$--manifold with nonempty
boundary, no component of which is a $2$--sphere. Then $M$ admits
a triangulation having all vertices in the boundary and precisely
one vertex in each boundary component.
\end{cor}

We need the following proposition. It is analogous to Proposition
\ref{0-decide}. An alternate proof appears in \cite{jac-ree}.

\begin{prop} \label{0-b-decide} Given a triangulation $\T$ of a compact,
orientable $3$--manifold with nonempty boundary, no component of
which is a $2$--sphere, there is an algorithm to decide if $\T$ is
$0$--efficient; furthermore, the algorithm will construct a non
vertex-linking  normal disk, if one exists.
\end{prop}

\begin{proof} First, we observe that if we are given a normal $2$--sphere, we can
construct a non vertex-linking normal disk. This uses the
techniques from above; namely if $S$ is a normal $2$--sphere, then
there is an arc in the $1$--skeleton of $\T$, meeting $S$ in a
single point and meeting the boundary of $M$ in a vertex of the
triangulation. A small regular neighborhood of this arc and $S$
has boundary a normal $2$--sphere and a properly embedded disk,
possibly not normal; however, they are barrier surfaces in the
components of their complement not containing $S$ (and the arc).
We shrink the disk. We get  a  non vertex-linking, normal disk.
This is where we use that no component of the boundary of $M$ is a
$2$--sphere.

So, we may assume there are no interior vertices and hence, no
vertex-linking normal $2$--spheres. Among all (non vertex-linking)
normal $2$--spheres and all non vertex-linking normal disks,
suppose $\Sigma$ is one for which the dimension of its carrier in
the projective solution space of $\T$ is a minimum. We claim
$\Sigma$ is carried by a vertex.

If $\Sigma$ is not carried by a vertex of the projective solution
space, then there are non  negative integers $k, n$, and $m$ and
normal surfaces $X$ and $Y$ so that $k\Sigma = nX + mY$, where $X$
and $Y$ are carried by {\it proper} faces of the carrier of
$\Sigma$. It follows as above that either $\chi (X) \ge 0$ or $
\chi (Y) \ge 0$. Suppose $\chi(X)\ge 0$. Then a component of $X$,
say $X'$, is either a normal $2$--sphere, projective plane or disk
and is carried by a proper face of the carrier of $\Sigma$. Now,
by our observations, $X'$ is not a vertex-linking normal
$2$--sphere (this could only come from an interior vertex). We may
also assume that $X'$ is not a vertex-linking normal disk; for
then, it would be, itself, a factor of $k\Sigma$. Hence, $X'$ is
either a non vertex-linking disk, a non vertex-linking $2$--sphere
or it is a projective plane and therefore, its double is a non
vertex-linking $2$--sphere. But this contradicts our choice of
$\Sigma$.
\end{proof}

We could have begun this section with the following theorem  which
has as an immediate corollary that a compact, orientable
$3$--manifold with nonempty boundary, which is known to be
irreducible and $\bdy$--irreducible, admits a $0$--efficient
triangulation. It does {\it not} use the $3$--sphere recognition
algorithm as the analogous theorem, Theorem \ref{0-eff-construct},
does in the closed case. This is because, for a bounded
$3$--manifold, anytime we have a $2$--sphere and we know it bounds
a $3$--cell, then we know for which side of the $2$--sphere we
have the $3$--cell.

\begin{thm} \label{b-0-eff-construct} Suppose $M$ is a compact, orientable, irreducible and
$\bdy$--irreducible $3$--manifold  with non-empty boundary. Then
any triangulation of $M$ can be modified to a $0$--efficient
triangulation or $M$ is a $3$--cell.
\end{thm}
\begin{proof} Suppose $M$ is given as in the hypothesis and $\T$ is
 a triangulation of $M$. We can decide
if $\bdy M$ is a $2$--sphere. If it is, then $M$ is a $3$--cell
and there is nothing to prove; so, we may assume no component of
$\bdy M$ is a $2$--sphere.

Let $S_1,\ldots,S_n$ denote the boundary components of $M$ and let
$\T_i$ denote the triangulation $\T$ restricted to $\S_i$. We want
to use the techniques of the proof of Theorem
\ref{b-0-eff-exists}; there we had a pairwise disjoint collection
of normal disks, one for each component of $\bdy M$ to guide a
crushing of the triangulation. Here we can proceed in several
ways. For example, if there is more than one vertex in a component
of $\bdy M$ or if there is an interior vertex, then we can begin
straight away to construct the desired collection of normal disks.
On the other hand, if there is only one vertex in each boundary
component of $M$ and no vertices in $\stackrel{\circ}{M}$, then we
need to use the algorithm in Proposition \ref{0-b-decide} to
decide if $\T$ is $0$--efficient and if it is not, then to find a
non vertex-linking disk. It turns out that this latter approach
works in both situations but is a bit over the top in the case it
is obvious that $\T$ is not $0$--efficient.

So, we apply the algorithm of Proposition \ref{0-b-decide}. If
$\T$ is $0$--efficient there is nothing to prove. If it is not
then there is a non vertex-linking  normal disk and the algorithm
will construct one for us. Say $E_1$ is such a disk and notation
has been chosen so that $\bdy E_1\subset S_1$. Since $M$ is
irreducible and $\bdy$--irreducible, there is a disk $E_1'\subset
S_1$, $\bdy E_1 = \bdy E_1'$ and $E_1\cup E_1'$ bounds a $3$--cell
$B_1$. Since no component of $\bdy M$ is a $2$--sphere, we know
precisely which component of the complement of $E_1\cup E_1'$
meets $B_1$. Let $E_i$ be a vertex-linking disk in $S_i, i\ne 1$;
then for each $i$ we have a disk $E_i'\subset S_i, \bdy E_i' =
\bdy E_i$ and $E\i\cup E_i'$ bound a $3$--cell, say $B_i$. Thus we
have a pairwise disjoint collection of normal disks $\mathcal{E} =
\{E_1,\ldots,E_n\}$, so that $\bdy E_i\subset S_i$, for each $i$
there is a disk $E_i'\subset S_i, \bdy E_i = \bdy E_i'$, $E\i\cup
E_i'$ bounds a $3$--cell $B_i$ and $E_1$ is not vertex-linking.
However, we want to have all the vertices of $\T$ contained in
$\bigcup B_i$.

So, suppose not all the vertices of $\T$ are in $\bigcup B_i$.
First, we want to make sure all the vertices of $\T$ in $\bdy M$
are in some $B_i$. If this is not the case, then there is a vertex
$v$ of $\T$, $v$ is not in $\bigcup B_i$ and $v$ is in some $\T_j$
($v\in S_j$). In this case, there is an arc $\alpha_j$ in the
$1$--skeleton of $\T_j$ having one end point $v$ and the other in
$\bdy E_j'$ and otherwise, missing $E_j'$.  Let $N_j$ be a small
regular neighborhood of $B_j\cup\alpha_j$. Then since
$\alpha_j\subset S_j$, we have $N_j$ is a $3$--cell (a relative
punctured $3$--cell); and if $F_j$ is the frontier of $N_j$, then
$F_j$ is a properly embedded disk and $F_j$ along with the $E_i,
i\ne j$, form a barrier surface in the component of the complement
of $F_j$ not meeting $B_j\cup\alpha_j$. We shrink $F_j$. Then just
as so many times above, in this way, we get a normal disk
$\hat{E}_j$, $\bdy \hat{E}_j\subset S_j$ and there is a disk
$\hat{E'}_j\subset S_j$ so that $\bdy\hat{E}_j = \bdy\hat{E'}_j$
and $\hat{E}_j\cup\hat{E'}_j$ bounds a $3$--cell $\hat{B}_j$ which
contains $B_j\cup\alpha_j$. Hence, $\hat{B}_j$ contains $v$.
Furthermore, $\hat{E}_j$ (as well as $E_1$) is not vertex-linking.
It is possible that $j = 1$, but still $\hat{E}_j$ is not
vertex-linking.

Thus, if we stick with our original notation in order to keep
things simpler, we now may assume that we have constructed a
pairwise disjoint collection of normal disks $\{E_1,\ldots,E_n\}$,
so that $\bdy E_i\subset S_i$, for each $i$ there is a disk
$E_i'\subset S_i, \bdy E_i = \bdy E_i'$, $E\i\cup E_i'$ bounds a
$3$--cell $B_i$, all the vertices of $\T$ in $\bdy M$ are in
$\bigcup B_i$ and $E_1$, at least, is not vertex-linking.

If all the vertices of $\T$ are not in $\bigcup B_i$, then there
is a vertex $v\in\stackrel{\circ}{M}$ not in some $B_i$. Thus
there is an arc $\alpha$ in the $1$--skeleton of $\T$ so that $v$
is at one end of $\alpha$ and the other end of $\alpha$ is in some
$E_j$ and otherwise $\alpha$ misses $\bigcup B_i$; furthermore,
since all the vertices in $\bdy M$ are in some $B_i$, the entire
arc $\alpha\subset\stackrel{\circ}{M}$. Let $N_j$ be a small
regular neighborhood of $B_j\cup\alpha$. Then  $N_j$ is a
$3$--cell (a relative punctured $3$--cell); and if $F_j$ is the
frontier of $N_j$, then $F_j$ is a properly embedded disk and
$F_j$ along with the $E_i, i\ne j$, form a barrier surface in the
component of the complement of $F_j$ not meeting $B_j\cup\alpha$.
From here the argument is just as above.

So, we can construct a pairwise disjoint collection of normal
disks $\mathcal{E} = \{E_1,\ldots,E_n\}$ so that $\bdy E_i\subset
S_i$, for each $i$ there is a disk $E_i'\subset S_i, \bdy E_i =
\bdy E_i'$, $E\i\cup E_i'$ bounds a $3$--cell $B_i$, all the
vertices of $\T$ are contained in $\bigcup B_i$ and $E_1$, at
least, is not vertex-linking.

We are now set up to follow the steps in the proof of Theorem
\ref{b-0-eff-exists}; however, the difference is that we do not
have here that the normal disks in the collection $\mathcal{E}$
are maximal. (Recall we are using maximal here in the sense that
if $F$ is a normal disk in $M$, $F$ is disjoint from the disks in
the collection $\mathcal{E}$ and there is a disk $F'\subset S_j$,
$\bdy F' = \bdy F$, $F'\cup F$ bounds a $3$--cell $B$, $B_j\subset
B$, then we must have $F = E_j$.) We will follow the steps of the
proof of Theorem \ref{b-0-eff-exists}; while we will be very brief
in doing this, we do need some terminology for the situations
where we are not able to use that the disks in the collection
$\mathcal{E}$ are maximal.

So, suppose $\mathcal{E}$ is as above and there is a pairwise
disjoint collection of normal disks $\mathcal{F} =
\{F_1,\ldots,F_n\}$ so that $\bdy F_i\subset S_i$, for each $i$
there is a disk $F_i'\subset S_i, \bdy F_i = \bdy F_i'$, $F\i\cup
F_i'$ bounds a $3$--cell $B_i'$, and for all $i, B_i\subseteq
B_i'$. In this case we say that $\mathcal{F}$ is {\it larger} than
$\mathcal{E}$ and if for some $j$, $F_j$ is not equivalent to
$E_j$, we say $\mathcal{F}$ is {\it strictly larger} than
$\mathcal{E}$.

Let $X$ be the closure of the component of the complement of
$\bigcup B_i$ in $M$. Then $X$ is homeomorphic to $M$ and if $\C$
is the induced cell decomposition on $X$, the cells of $\C$ are of
Type I, II, III, and IV (all the vertices of $\T$ are in $\bigcup
B_i$).

\vspace{.15 in} \noindent {\it Claim. $\bbb{P}(\C)\ne X$ and each
component of $\bbb{P}(\C)$ is a product $I$--bundle.}

\vspace{.1 in}\noindent {\it Proof of Claim.} This follows just as
in the proof of Theorem \ref{b-0-eff-exists}.

\vspace{.15 in} \noindent {\it Claim. There is a trivial induced
product region $\bbb{P}(X)$ for $X$ or we can construct a strictly
larger collection of normal disks than $\mathcal{E}$. }

\vspace{.1 in} \noindent {\it Proof of Claim}. We have established
that each component of $\bbb{P}(\C)$ is a product $I$--bundle. As
above, we write each I-bundle component of $\bbb{P}(\C)$ as
$K_i\times I, i = 1,\ldots,k$, where $k$ is the number of
components of $\bbb{P}(\C)$ and again, we set $K_{i}^\ve =
K_i\times ve, ve = 0, 1$. The subcomplexes $K_i^0$ and   $K_i^1$
are isomorphic subcomplexes in the induced normal cell structures
on the disks in our collection $\mathcal{E}$. Here, we may have
$K_{i}^0$ and $K_{i}^1$ in the same or in distinct disks in
$\mathcal{E}$.

We may as well assume that $K_i$ is not simply connected. Recall
that if $K_{i}^\ve$ is in the disk $E_j$ in the collection
$\mathcal{E}$, then $K_i^\ve$ separates $E_j$. If $K_i^\ve$ meets
$\bdy E_j$, then every component of the complement of $K_i^\ve$ is
simply connected; however, if $K_i^\ve$ does not meet $\bdy E_j$,
then the component of the complement of $K_i^\ve$ in $E_j$, which
contains $\bdy E_j$ is not simply connected.

As above, we consider the possibilities for $K_i\times [0,1]$.

If $K_i^0$ and $K_i^1$ are in distinct disks of our collection
$\mathcal{E}$, say $K_i^0\subset E_j$ and $K_i^1\subset E_{j'},
j\ne j'$. Then we let $D_i^0$ be the union of $K_i^0$ along with
the components of the complement of $K_i^0$ in $E_j$, which do not
meet $\bdy E_j$; similarly, we let $D_i^1$ be the union of $K_i^1$
along with the components of the complement of $K_i^1$ in
$E_{j'}$, which do not meet $\bdy E_{j'}$. Then $D_i^\ve, \ve =
0,1$ is simply connected. So, as in the proof of Theorem
\ref{b-0-eff-exists}, it follows that there is a simply connected
planar complex $D_i$ and an embedding of $D_i\times [0,1]$ into
$X$ so that $D_i\times 0 = D_i^0, D_i\times 1 = D_i^1, K_i\times
[0,1]\subset D_i\times [0,1]$, and the frontier of $D_i\times
[0,1]$ is contained in the frontier of $K_i\times [0,1]$.

If both $K_i^0$ and $K_i^1$ are in the same disk $E_j$ and $K_i^0$
(and hence, $K_i^1$) does not meet $\bdy E_j$, then there are two
possibilities.

One possibility is that $K_i^0$ is in the same component of the
complement of $K_i^1$ as $\bdy E_j$ and $K_i^1$ is in the same
component of the complement of $K_i^0$ as $\bdy E_j$. Then we let
$D_i^0$  be the union of $K_i^0$ along with the components of the
complement of $K_i^0$ in $E_j$, which do not meet $\bdy E_j$;
similarly, we let $D_i^1$ be the union of $K_i^1$ along with the
components of the complement of $K_i^1$ in $E_{j}$, which do not
meet $\bdy E_{j}$. Then $D_i^\ve, \ve = 0,1$ is simply connected.
It follows, just as above and as in the proof of Theorem
\ref{b-0-eff-exists}, that there is a simply connected planar
complex $D_i$ and an embedding of $D_i\times [0,1]$ into $X$ so
that $D_i\times 0 = D_i^0, D_i\times 1 = D_i^1, K_i\times
[0,1]\subset D_i\times [0,1]$, and the frontier of $D_i\times
[0,1]$ is contained in the frontier of $K_i\times [0,1]$.

So, suppose $K_i^1$ is not in the component of the complement of
$K_i^0$ that meets $\bdy E_j$. Define $D_i^0$ to be $K_i^0$ along
with all the components of its complement not meeting $\bdy E_j$.
Then we have that $D_i^0$ is simply connected and also we have
$K_i^1\subset D_i^0$. See Figure \ref{f-b-product}. Let $N_i =
N(D_i^0\cup (K_i\times [0,1]))$ be a small regular neighborhood of
$D_i^0\cup (K_i\times [0,1])$. Since $K_i^1\subset D_i^0$, the
frontier of $N_i$ consists of an annulus, possibly some
$2$--spheres and a disk $F_i$, properly embedded in $X$ and having
its boundary in $E_j$.  There is a disk $F_i'$ in $E_j$ so that
$\bdy F_i' = \bdy F_i$ and $F_i\cup F_i'$ bounds a $3$--cell, say
$B$, in $X$; actually, we want to think of $B$ as a relative
(punctured) $3$--cell. Furthermore, $B$ contains $N_i$ and $F_i$
along with $E_k, k\ne j$ form a barrier surface in the component
of the complement of $F_i$ not meeting $N_i$. We shrink $F_i$ in
the component of the complement of $B$ not meeting $N_i$.  Since
$M$ (and hence, $X$) is not a $3$--cell, our standard arguments
give that in this situation we have a  {\it normal} disk $E$ in
$X$ and a disk $E'$ in $\bdy X$ so that $\bdy E' = \bdy E$ and
$E\cup E'$ bounds a $3$--cell containing $P$. We replace $E_j$ in
the collection $\mathcal {E}$ by the disk $E$ and thus construct a
strictly larger collection of normal disks than $\mathcal{E}$. We
now begin our considerations over again, using this larger
collection of normal disks. This can only happen a finite number
of times by Theorem \ref{kneser}.

The only remaining possibility is that both $K_i^0$ and $K_i^1$
are in the same disk $E_j$ and now  $K_i^0$ (and hence, $K_i^1$)
does meet $\bdy E_j$. We can go through the steps in this case
just as in the proof of Theorem \ref{b-0-eff-exists} and we find
that we can either construct an appropriate product $D_i\times
[0,1]$ or we construct a strictly larger collection of normal disk
than $\mathcal{E}$. In the latter situation, as before, we go back
to the earlier steps and use this new collection of normal disks;
again, this phenomenon can happen at most a finite number of
times.

Now, having made these observations, we can choose some order for
the components of $\bbb{P}(\C)$, say $K_1\times
[0,1],\ldots,K_k\times [0,1]$, and from this, just as in the proof
of Theorem \ref{0-eff-exists} derive that there is a trivial
induced product region $\bbb{P}(X)$ for $X$.

This completes the proof of the Claim.

 \vspace{.15
in} \noindent {\it Claim. There is no  cycle (complete or
relative) of cells of type II, which are not in $\bbb{P}(X)$ .}

\vspace{.1 in} \noindent {\it Proof of Claim}. In this situation
we have four possibilities: a complete cycle about a single edge,
a complete cycle about more than one edge, a relative cycle about
a single edge and a relative cycle about more than one edge.

Suppose there were a complete cycle about a single edge $e$ as in
Figure \ref{f-cycle}, Part A. Just as in the proof of Theorem
\ref{b-0-eff-exists}, there would be a disk $D$ meeting $e$ in
precisely one point and meeting a non vertex-linking disk, say
$E_j$, of the collection $\mathcal{E}$ in $\bdy D$. A surgery on
$E_j$ at $D$ (see Figure \ref{f-surgery-cycle}), which is the same
as adding a $2$--handle to $E_j$ along the curve $\bdy D$, gives a
{\it normal} $2$--sphere, $S_j$, and a new {\it normal} disk,
$F_j$.  But here the existence of the disk $F_j$, which co-bounds
a $3$--cell with $E_j'$ gives us a strictly larger collection of
disks than the collection $\mathcal{E}$. We use this new
collection and go back to the beginning.

Now, suppose there is a complete cycle of truncated prisms about
more than one edge (see Figure \ref{f-cycle}(B)), then as in the
consideration in the proof of Theorem \ref{b-0-eff-exists}, the
collection of cells of type II (truncated prisms), form a solid
torus with, possibly, some self identifications in its boundary.
As above, we distinguish between the cycle of truncated prisms,
which we denote $\hat{\tau}$, and the cycle of truncated prisms
minus the bands of trapezoids, which we denote by $\tau$. We have
that $\tau$ meets the collection $\mathcal{E}$ in either three
open annuli, each meeting a meridional disk of $\tau$ precisely
once, or a single open annulus, meeting a meridional disk of
$\tau$ three times. See Figure \ref{f-cycle-torus}.

In the case we have three annuli in $\tau$ meeting the collection
$\mathcal{E}$, we have just as above that the cycle of truncated
prisms $\hat{\tau}$ is in $\bbb{P}(X)$, the induced product region
for $X$.

Suppose we have just one annulus common to $\tau$ and the
collection $\mathcal{E}$. Just as in the proof of Theorem
\ref{b-0-eff-exists}, we have that  there is a lens space $L(3,1)$
embedded in $M$ as a connected summand but this is impossible as
$M$ is irreducible.

Now, we consider the possibility of a relative cycle of truncated
prisms. If there is a relative cycle about a single edge, then we
can find a strictly larger collection of normal disk than our
collection $\mathcal{E}$. If there is a relative cycle about more
than one edge, then the argument in the proof of Theorem
\ref{b-0-eff-exists} can be used to show that such a relative
cycle of truncated prisms is in the induced product region for
$X$, $\bbb{P}(X)$.

This completes the proof of the Claim.

Hence, as above, we have established all the conditions of Theorem
\ref {b-crush} and so there is a triangulation $\T'$ of $X$, where
the tetrahedra of $\T'$ are in one-one correspondence with the
truncated tetrahedra of $\C$ not in $\bbb{P}(X)$. Since, we had
that at least one of the disks in the collection $\mathcal{E}$ is
not vertex-linking, then the number of tetrahedra of $\T'$ is
smaller than the number of tetrahedra in the original
triangulation $\T$. By repeated applications of the above
techniques, it follows that we start with $\T$ and through a
sequence of such modifications get a $0$--efficient triangulation
of $M$.\end{proof}

In the previous subsection our methods could be used to give a
connected sum decomposition of any $3$--manifold into known
factors and factors with $0$--efficient triangulations. A similar
result in the bounded case for disk sum decompositions is possible
but is not too interesting. On the other hand, there are a couple
of related and curious questions, which our methods do not answer
and for which we do not know the answer. Recall that we have
defined inessential $2$--spheres and inessential properly embedded
disks. We will say a properly embedded disk $D$ in the
$3$--manifold $M$ ($\bdy M\ne\emptyset$) is {\it trivial} if there
is a disk $D'\subset \bdy M, D\cap D' = \bdy D = \bdy D'$ and
$D\cup D'$ bounds a $3$--cell ($D$ is {\it parallel into $\bdy
M$}). An inessential disk is trivial in an irreducible
$3$--manifold. Now, does any closed, orientable $3$--manifold have
a triangulation in which every inessential, normal $2$--sphere is
vertex-linking? Does any compact, orientable $3$--manifold with
nonempty boundary have a triangulation in which there are no
inessential, normal $2$--spheres and every properly embedded,
trivial, normal disk is vertex-linking? We know from
\cite{jac-let-rub1} that there are one-vertex triangulations of
handlebodies, which have no normal $2$--spheres and every trivial,
normal disk is vertex-linking.

\section{$0$--efficient and minimal triangulations.}\label{min-triang}

 A triangulation
$\T$ of a $3$--manifold $M$ is said to be a {\it minimal} triangulation if for
any other triangulation $\T'$ of $M$, the number of tetrahedra of $\T$ is no
greater than the number for $\T'$. Similarly one defines a {\it minimal}
triangulation for a surface in terms of number of triangles.

Recall for a closed surface $S$, we have $f = 2(\nu - \chi (S))$,
where $f$ is the number of triangles and $\nu$ the number of
vertices in a triangulation of $S$. Hence, a minimal triangulation
occurs when the triangulation has precisely one vertex, except for
$S^2$ and $\rpp$, which have minimal triangulations with three and
two vertices, respectively. There does not seem to be any simple
way, such as this, to determine much about a minimal triangulation
of a $3$--manifold. However, directly from the proof of Theorem
\ref{0-eff-exists}, we have the following theorem.

\begin{thm} A minimal triangulation of a closed, orientable,
irreducible $3$--manifold is $0$--efficient or the manifold is
homeomorphic with either $\rp$ or $L(3,1)$; hence, a minimal
triangulation of a closed, orientable, irreducible $3$--manifold
has one-vertex unless the $3$--manifold is one of $S^3, \rp$ or
$L(3,1)$.
\end{thm}

There are two one-tetrahedron (minimal) triangulations of $S^3$;
one has one vertex and one has two vertices (see Figure
\ref{f-tetra} (4) and (5)). Both are $0$--efficient. There are two
distinct minimal triangulations of $\rp$, each having precisely
two tetrahedra. One of these triangulations has one vertex and the
other has two vertices. For the example with two vertices, see
Figure \ref{f-not-0-eff}; and for the one with one vertex, see
Figure \ref{f-onevertex-RP3-S2xS1}(A). Of course, as observed
above, neither is $0$--efficient. There are four distinct minimal
triangulations of $L(3,1)$, each having two tetrahedra. Two of
these triangulations have one vertex (see Figure \ref{f-two-L3_1})
and two have two vertices (see Figure \ref{f-not-0-eff}). Of the
two triangulations of $L(3,1)$ with one vertex, only one is
$0$--efficient. Neither of the two tetrahedra, two-vertex
triangulations of $L(3,1)$ are $0$--efficient (see Proposition
\ref{p-0-eff}).

In the case of manifolds with boundary, we have the following,
again, directly from the proof of Theorem \ref{b-0-eff-exists}.

\begin{thm} A minimal triangulation of a compact, orientable,
irreducible and $\bdy$--irreducible $3$--manifold with nonempty
boundary is $0$--efficient, or it is the underlying point set of a
tetrahedron (the $3$--cell triangulated by a single tetrahedron
with no identifications). Hence, a minimal triangulation of such a
$3$--manifold has all vertices in the boundary and just one vertex
in each boundary component or is a $3$--cell and is triangulated
by just one tetrahedron and has either three or four vertices, all
in the boundary.
\end{thm}

As pointed out in the previous Theorem, there are two (minimal)
one-tetrahedron triangulations of the $3$--cell. One is the
tetrahedron with no identifications, see Figure \ref{f-tetra} (1),
and is not $0$--efficient. The other is a single  tetrahedron with
two faces identified,  see Figure \ref{f-tetra}(2); it has three
vertices and is $0$--efficient.

\begin{figure}[htbp]

            \psfrag{a}{$a$}
            \psfrag{b}{$b$}
            \psfrag{c}{$c$}
            \psfrag{d}{$d$}
            \psfrag{e}{$e$}
            \psfrag{f}{$f$}
            \psfrag{D}{$\partial$}
            \psfrag{B}{genus 2 handlbody - layered triangulation}
            \psfrag{A}{\begin{tabular}{c}
            once punctured\\
   Klein bottle\\
            \end{tabular}}
        \vspace{0 in}
        \begin{center}
        \epsfxsize=3.5 in
        \epsfbox{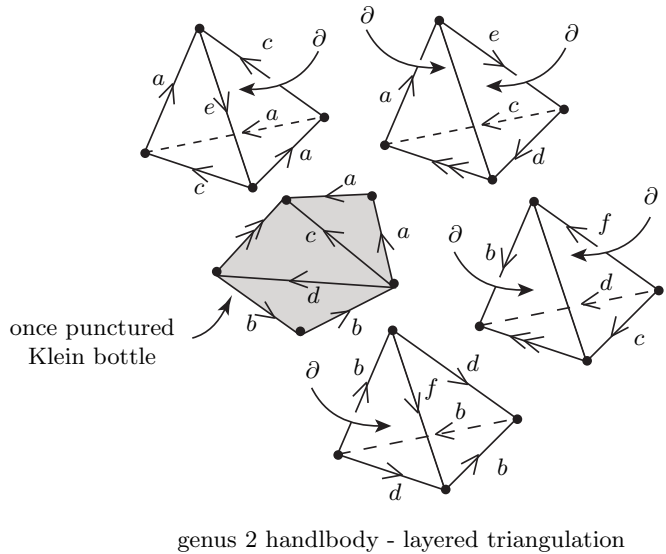}
        \caption{A minimal triangulation of the genus 2 handlebody, having
$t = 3g - 2 = 4$ tetrahedra and one vertex. Four tetrahedra
``layered" on a once-punctured Klein bottle.}
        \label{f-genus2}
        \end{center}

\end{figure}

We do not know if a minimal triangulation of a closed,  reducible
and orientable $3$--manifold  needs to be a one-vertex
triangulation. However, we can use Euler characteristic to show
that a minimal triangulation of a compact, orientable
$3$--manifold with nonempty boundary, no component of which is a
$2$--sphere, has precisely one vertex in each boundary component
and the number of tetrahedra can be expressed as $$t =
\sum_{j=1}^{b}(3g_j - 2) + e_{int} - v_{int},$$ where $b$ is the
number of boundary components, $g_j$ is the genus of the $j^{th}$
boundary component, and $e_{int}$ and $v_{int}$ are the number of
interior edges and interior vertices of the triangulation,
respectively. Since $e_{int} - v_{int} \geq 0$, if there are no
interior edges, then we would have the number of tetrahedra in a
minimal possible triangulation. Of course, if there are no
interior edges, then if the manifold is also irreducible, it is a
handlebody. There are triangulations of a genus $g$ handlebody,
which have just one vertex (of course, in the boundary) and
realize this minimum number of tetrahedra, $t = 3g - 2$. See
Figure \ref{f-tetra}(3) for $g = 1$ and Figure \ref{f-genus2} for
$g = 2$. A study of ``layered triangulations" of handlebodies is
given in \cite{jac-let-rub1}. The minimal  triangulation of
$S^2\times S^1$, also,  is one-vertex (see Figure
\ref{f-onevertex-RP3-S2xS1}).

We complete this section on  $0$--efficient and minimal
triangulations with some observations about the order of edges in
these triangulations. We observed above, Corollary \ref{nodisk},
that a $0$--efficient triangulation has no edge of order one and
no faces which are cones unless the manifold is $S^3$.

\begin{prop} A minimal and $0$--efficient triangulation $\T$ of
the closed, orientable and
irreducible
$3$--manifold $M$ has
\begin{itemize}
\item no edge of order one unless $M = S^3$, \item no edge of
order two unless $M = L(3,1)$ or $L(4,1)$, and \item no edge of
order three unless either $M = L(5,2)$ or $\T$ contains, as a
subcomplex, the two-tetrahedron, geodesic layered triangulation
$\{4,3,1\}$ of the solid torus.
\end{itemize}
\end{prop}

\begin{figure}[htbp]

            \psfrag{1}{$1$}
            \psfrag{2}{$2$}
            \psfrag{3}{$3$}
            \psfrag{4}{$4$}
            \psfrag{5}{$2\leftrightarrow 4$}
             \psfrag{C}{$\{3,4,1\}\leftrightarrow\{3,2,1\}$}
              \psfrag{D}{\Large{$L(6,1)$}}
            \psfrag{A}{\begin{tabular}{c}
           $\{4,3,1\}$ layered\\
   solid torus\\
            \end{tabular}}
             \psfrag{B}{\begin{tabular}{c}
           $\{3,2,1\}$ layered\\
   solid torus\\
            \end{tabular}}
        \vspace{0 in}
        \begin{center}
        \epsfxsize=3 in
        \epsfbox{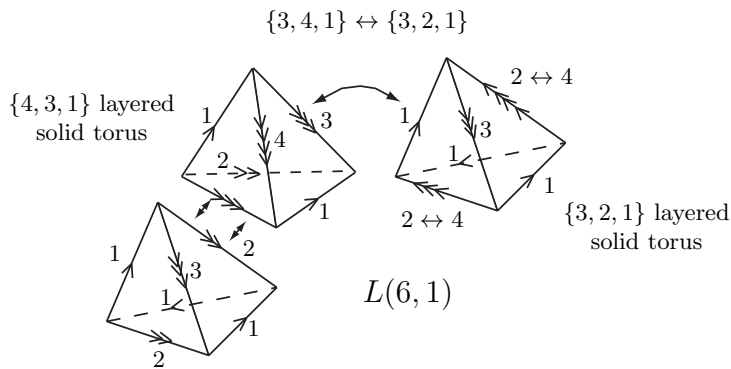}
        \caption{Layered triangulation of $L(6,1)$, which is a minimal
layered triangulation, is $0$--efficient,  has two edges of order
three and contains the layered solid torus $\{4,3,1\}$ as a
subcomplex.}
        \label{f-layered}
        \end{center}

\end{figure}

Before giving the proof, we recall that, with just few exceptions,
a minimal triangulation is $0$--efficient; however, we have
combined the two notions to reduce the number of exceptions in the
statement of the theorem.  Also, a $\{4,3,1\}$ geodesic layered
triangulation of the solid torus is given as a subcomplex in
Figure \ref{f-layered}. This is just one in a family of
triangulations of the solid torus having a central role in our
study of triangulations of $3$--manifolds, see \cite{jac-rub1,
jac-let-rub1, jac-let-rub2, jac-sed1, jac-let-rub-sed}. The
geodesic, layered triangulation of the lens space $L(6,1)$, shown
in Figure \ref{f-layered}, is minimal and $0$--efficient and has
two edges, those labeled $``2"$ and $``4"$, each of order three.
There are many examples of $0$--efficient triangulations of lens
spaces, which have edges of order three; however, while we suspect
the geodesic, layered triangulations of lens spaces are minimal,
we do not have a proof of this.

\begin{proof} Let $\T$ be the given minimal and $0$--efficient triangulation.
If there is an edge of order one, then we saw above that $M =
S^3$. So, suppose we have an edge, $e$, of order two. If it is in
just one tetrahedron, then the triangulation $\T$ has just one
tetrahedron. In this case, we see in Figure \ref{f-tetra}(6) that
we must have $L(4,1)$. So, we suppose there are two tetrahedra
$\td{\Delta}'$ and $\td{\Delta}''$ in our triangulation $\T$,
which are identified along two faces from each and the faces have
the edge $e$ in common. Let $e'\subset\td{\Delta}'$ and
$e''\subset\td{\Delta}''$ be edges disjoint from $e$ in
$\td{\Delta}'$ and $\td{\Delta}''$, respectively, so that
$\td{\Delta}' = e * e'$ and $\td{\Delta}'' = e * e''$. See Figure
\ref{f-order-two}. Now, $e'\cup e''$ is a loop in the
$1$--skeleton of our triangulation. There are three possibilities
for identifications of $e'$ and $e''$ we must consider: $e'$ and
$e''$ are not identified, $e'$ and $e''$ are identified with
``opposite" orientations, and $e'$ and $e''$ are identified with
the ``same" orientation.  In the latter situation, there would be
an embedded $\rpp$ in the manifold $M$; but this contradicts the
triangulation being $0$--efficient. So, we may assume the latter
possibility does not occur.

Let $A'B'C'$ and $D'B'C'$ denote the two faces of $\td{\Delta}'$
containing $e' = \overline{B'C'}$ and let $A''B''C''$ and
$D''B''C''$ denote the two faces of $\td{\Delta}''$ containing
$e'' = \overline{B''C''}$. Again, see Figure \ref{f-order-two}. In
general, we should find that we can collapse the two tetrahedra
$\td{\Delta}'$ and $\td{\Delta}''$ identifying the faces $A'B'C'$
with $A''B''C''$ to get a single face $ABC$ and similarly for
$D'B'C'$ and $D''B''C''$. However, if we could collapse, then the
triangulation would not be minimal; so, in a minimal
triangulation, if we have an edge, $e$, of order two, there is an
obstruction to such a collapse. The obstruction is that we have
the face $ABC$  already identified to the face $A'B'C'$  in some
way  (or symmetrically, $DEF$ is already identified with $D'E'F'$
in some way) or we have a two tetrahedron triangulation.

Suppose $ABC$ is identified with $A'B'C'$ in some way. The
possibilities are $ABC \leftrightarrow A'B'C'$, $ABC
\leftrightarrow B'C'A'$ and $ABC \leftrightarrow C'A'B'$. Either
of the latter two identifications show that $M = L(3,1)$; hence,
if $\T$ is minimal and $0$--efficient, we have $\T$ the
triangulation given in Figure \ref{f-two-L3_1}(A). If we have $ABC
\leftrightarrow A'B'C'$ and we are not in a two-tetrahedron
triangulation, then we have a cone and we could identify $DEF$
with $D'E'F'$ and reduce the number of tetrahedra. So, we must
have a two tetrahedron triangulation with an edge of order two.
Since, it is $0$--efficient, it is $L(3,1)$. Both $\rp$ and
$\S^2\times S^1$ have minimal two-tetrahedra triangulations with
edges of order two (see Figure \ref{f-not-0-eff} and
\ref{f-onevertex-RP3-S2xS1}(B)); however, these are not
$0$--efficient.

\begin{figure}[htbp]

            \psfrag{a}{$A'$}
            \psfrag{b}{$B'$}
            \psfrag{c}{$C'$}
            \psfrag{d}{$D'$}
            \psfrag{e}{$\Large\mathbf{e}$}
            \psfrag{f}{$e'$}
            \psfrag{g}{$e''$}
            \psfrag{x}{$A''$}
            \psfrag{y}{$B''$}
            \psfrag{z}{$C''$}
            \psfrag{w}{$D''$}
            \psfrag{A}{$A$}
            \psfrag{B}{$B$}
            \psfrag{C}{$C$}
            \psfrag{D}{$D$}
            \psfrag{X}{\LARGE{$\td{\Delta}'$}}
            \psfrag{Y}{\LARGE{$\td{\Delta}''$}}
            \psfrag{Z}{\LARGE{$\Delta'$}}
            \psfrag{W}{\LARGE{$\Delta''$}}
 \psfrag{i}{\begin{tabular}{c}
           faces\\
   identified\\
            \end{tabular}}
            \psfrag{Q}{\begin{tabular}{c}
           $A'B'D'\leftrightarrow A''B''D''$\\
   and\\ $A'C'D'\leftrightarrow A''C''D''$\\
            \end{tabular}}
        \vspace{0 in}
        \begin{center}
        \epsfxsize=2.5 in
        \epsfbox{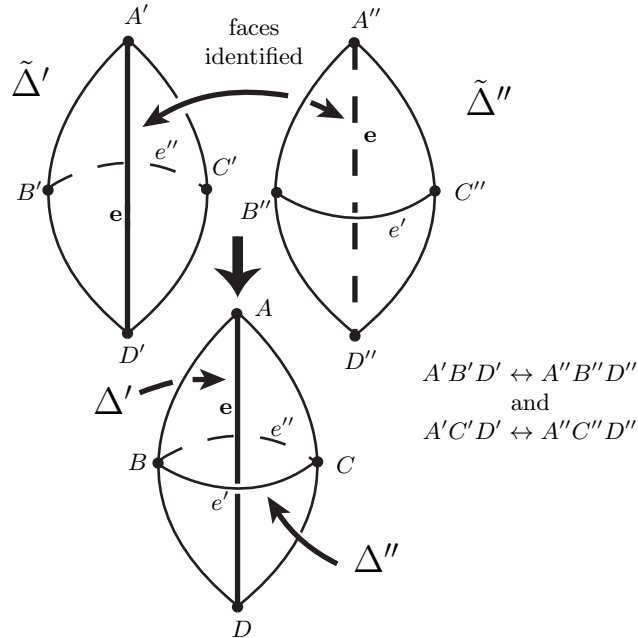}
        \caption{An edge $e$ of order two.}
        \label{f-order-two}
        \end{center}

\end{figure}

So, suppose we have an edge of order three. If there are three
distinct tetrahedra around the edge, then we could make a
$3\rightarrow 2$ move and reduce the number of tetrahedra. So, if
we assume $\T$ is minimal, we must have that the edge of order
three is in just one or in just two tetrahedra.

If we have an edge $e$ of order three in just one tetrahedron,
then we must have two faces, having $e$ as a common edge
identified. Since the edge $e$ is order three, this is only
possible if the triangulation has just one tetrahedron. Hence, we
must have $L(5,2)$; see Figure \ref{f-tetra}(7). So, we assume we
have the edge of order three is in two tetrahedra. So, again, we
have a tetrahedron with two faces identified but the edge is not
that edge common to both faces. See Figure \ref{f-tetra}, Parts(2)
and (3). However, we can not have an edge of order one; so, we
have a one-tetrahedron solid torus, Figure \ref{f-tetra}(3). Since
the edge is of order two in this single tetrahedron solid torus,
and we are assuming it is of order three, we have a tetrahedron
layered on the faces adjacent to these edges. It follows we have
the ``layered," two tetrahedron, $\{4,3,1\}$ solid torus as a
subcomplex.\end{proof}

\section{$0$--efficient and minimal ideal triangulations}

Suppose $M$ is a compact $3$--manifold with boundary, no component
of which is a $2$--sphere, and $\T$ is  an ideal triangulation of
the interior of $M$. Recall that in this work,  we do not allow
any ``regular" vertices in an ideal triangulation; i.e., vertices
with vertex-linking surface a $2$--sphere. Hence, all vertices are
ideal and the index of a vertex is $\geq 1$.  We say the ideal
triangulation $\T$ of $\stackrel{\circ}{M}$ is a {\it
$0$--efficient ideal triangulation} if and only if there are no
normal $2$--spheres.  We say an ideal triangulation of
$\stackrel{\circ}{M}$ is
 a {\it
minimal ideal triangulation} if and only if for any ideal triangulation $\T'$ of
$\stackrel{\circ}{M}$, $\T$ has no more tetrahedra than $\T'$.
Clearly, if $M$ has a
$0$--efficient ideal triangulation, then $M$ is irreducible.

Our first result in this section is an existence theorem for ideal
triangulations. Typically, one incurs ideal triangulations in
methods to construct a complete hyperbolic metric on a link
complement in the $3$--sphere. There are several algorithms
starting with a link projection to construct an ideal
triangulation of the link complement. The ones we know all give at
some point in the construction an ideal triangulation of the link
complement as a combinatorial object, just as we have defined.
However, beyond ideal triangulations of link complements in the
$3$--sphere, we have not seen any results on existence of ideal
triangulations. Our methods provide useful generalizations as in
\cite{jac-rub1}, where we construct ideal triangulations of link
complements in $3$--manifolds other than $S^3$,
 and in the next theorem, where we obtain ideal triangulations of the interior of
a compact $3$--manifold, allowing genus $> 1$ boundary components.
These methods also allow us to modify any ideal triangulation of
these $3$--manifolds to a $0$--efficient ideal triangulation and
show that a minimal ideal triangulation is $0$--efficient.

We will say the $3$--manifold $M$ is {\it anannular} if every
properly embedded annulus in $M$ is parallel into $\bdy M$.

\begin{thm} \label{i-exist}Suppose $M$ is a compact, irreducible, $\bdy$--irreducible,
anannular $3$--manifold. Then $\stackrel{\circ}{M}$ admits an
ideal triangulation.
\end{thm}

 \begin{proof} Let $\T$ be a triangulation of the $3$--manifold $M$ so
  that each vertex of
$\T$ is in $\bdy M$. Recall by \cite{bin3}, the manifold $M$ has a
triangulation with all vertices in the boundary; of course, under
our hypothesis and  Theorem \ref{b-one-vertex}, there is a
triangulation of $M$ that not only has all vertices in $\bdy M$
 but has precisely one vertex in each component of $\bdy M$. Let
$S_1,\ldots,S_n$ denote the components of $\bdy M$; let $N$ denote
a collared neighborhood of $\bdy M$; and let $E_i$ denote the
component of the frontier of $N$ that is isotopic to $S_i, i =
1,\ldots,n$. The frontier of $N$  is a barrier surface in the
component of its complement not meeting $\bdy M$. Hence, we can
shrink each $E_i$. It may be the case that some of the $E_i$ are
normal and thus already stable and the shrinking does not change
these surfaces. By shrinking, we arrive at a collection of normal
surfaces and possibly some $2$--spheres embedded entirely in the
interior of tetrahedra. However, since $M$ is irreducible and
$\bdy$--irreducible, each $E_i$ is incompressible; so, we have
precisely one normal surface for each $E_i$, which is isotopic to
$E_i$ and thus isotopic to the component $S_i$ of $\bdy M$. We may
also have some number of normal $2$--spheres. We discard all the
$2$--spheres. We continue to denote the (now) normal surface,
which is parallel to the boundary component, $S_i$ by $E_i, 1\leq
i\leq n$; and we will call the product region determined by the
isotopy between $S_i$ and $E_i$, $P_i$. Having shown existence of
the normal surfaces $E_i$, we want to take a maximal such
collection in the following sense. We choose a collection of
normal surfaces $E_1,\ldots,E_n$ so that $E_i$ is isotopic to
$S_i, 1\le i\le n$, and having the property that if  $E_i'$ is a
normal surface isotopic to $S_i$ and $E_i'\cap E_j = \emptyset,
j\ne i$, then $E_i'$ is normally isotopic into $P_i$, the product
region between $E_i$ and $S_i$. This is possible by Kneser's
Finiteness Theorem, Theorem \ref{kneser}.

Let $X$ denote the closure of the complement of $\bigcup_{i =1}^n
P_i$. The compact, bounded  $3$--manifold $X$ is homeomorphic to
$M$ and has a nice cell structure $\C$ induced by the
triangulation $\T$ and the fact that each vertex of $\T$ is in
$\bdy M$ (see Figure \ref{f-cell-decomp}). Now, we want to crush
the triangulation $\T$ along the surfaces in the collection
$E_1,\ldots, E_n$, distinct points for each $E_i$.

Just as above, we need to show the conditions of Theorem
\ref{crush} are satisfied; however, here we have that after
crushing we do not have a manifold at the vertex points; however,
we arrive at an ideal triangulation of $\stackrel{\circ}{X}$,
which is homeomorphic to $\stackrel{\circ}{M}$.

We proceed as in the proof of Theorem \ref{0-eff-exists} but,
since the boundary components of $X$ have genus $\ge 1$, we have
some new considerations.

\vspace{.15 in} \noindent {\it Claim. $\bbb{P}(\C)$ is a product
$I$--bundle and $\bbb{P}(\C)\ne X$.} \vspace{.1 in}

\noindent {\it Proof of Claim}. If all cells in $X$ were of type
III and IV, $\bbb{P}(\C) = X$, then $X$ is an $I$--bundle over a
closed surface, having non positive Euler characteristic (no
component of $\bdy M$ is a $2$--sphere). But then $M$ also would
be such an $I$--bundle and this would contradict that $M$ is
anannular. Similarly, if $\bbb{P}(\C)$ is not a product
$I$--bundle, then there is a M\"obius band properly embedded in
$X$ with its boundary in some $E_i$. If this were the case, then
the frontier of a small regular neighborhood of this M\"obius band
would give an essential annulus in $M$.

This completes the proof of the {\it Claim}.

\vspace{.15 in} \noindent {\it Claim. There is a trivial induced
product region $\bbb{P}(X)$ for $X$. } \vspace{.15 in}

\vspace{.1 in} \noindent {\it Proof of Claim.} We will show that
any obstruction to this claim leads either to a contradiction that
$M$ is anannular or to a contradiction of the collection
$E_1,\ldots,E_n$ being maximal.

We use the same notation as in earlier sections. From our first
Claim, we have that each component of $\bbb{P}(\C)$ is a product
$I$-bundle, which we write $K_i\times [0,1]$. Just as above, we
have that $K_i^\ve = K_i\times \ve, \ve = 0,1$, has an induced
cell-decomposition from $K_i^0\subset  E_j$ and $K_i^1\subset
E_{j'}$ with $K_i^0$ isomorphic to $K_i^1$.

First, we observe that $K_i^\ve, \ve = 0,1$ are contained in
simply connected regions of the surfaces $E_j$ and $E_{j'}$. For
otherwise, we would have, say, $K_i^0$ meeting $E_j$ in a region
which is not simply connected. It follows that there is a properly
embedded annulus $A$ in $K_i\times [0,1]$, $A$ has one of its
boundaries in $K_i^0$ and the other in $K_i^1$ and the boundary of
$A$ in $K_i^0$ is not contractible in $E_j$. See Figure
\ref{f-annulus}. Since $M$ is $\bdy$--irreducible it follows that
the boundary of the annulus $A$ in $K_i^1$ is also not
contractible in $E_{j'}$. But since $M$ is anannular, the only
possibility is that $j' = j$ and the annulus $A$ is parallel into
$E_j$. Let $T$ denote the solid torus determined by this region of
parallelism. Then a small regular neighborhood of $P_j\cup T$ has
frontier, $E_j'$, isotopic to $E_j$; furthermore, $E_j'$, along
with all the $E_k, k\ne j$, form a barrier surface in the
component of the complement of $E_j'$ not meeting $P_j\cup T$. We
can shrink $E_j'$ but then we must get a normal surface isotopic
to $E_j$ but not normally isotopic to it. This is a contradiction
to the collection $E_1,\ldots,E_n$ being maximal. Thus, we have
that any component of $\bbb{P}(\C)$ has a product structure
$K_i\times [0,1]$ and its ``top," $K_i^1$, and ``bottom," $K_i^0$,
meet the surfaces $E_j$ in a union of triangles and quadrilaterals
in the induced cell structure on $E_j$ forming disjoint isomorphic
subcomplexes each of which are contained in simply connected
regions of the surfaces $E_1,\ldots, E_n$.

\begin{figure}[htbp]

            \psfrag{K}{$\Large{K_i\times[0,1]}$}
            \psfrag{0}{$K_i^0$}
            \psfrag{1}{$K_i^1$}
            \psfrag{B}{$A =$ annulus}
            \psfrag{J}{$K_i^0\subset E_j$}
            \psfrag{L}{$K_i^1\subset E_j$}
            \psfrag{D}{Cells of Type III and IV}
            \psfrag{C}{\begin{tabular}{c}
           $K_i\times [0,1]$\\
  $I$--bundle\\
            \end{tabular}}
            \psfrag{A}{$A$}
        \vspace{0 in}
        \begin{center}
        \epsfxsize=2.5 in
        \epsfbox{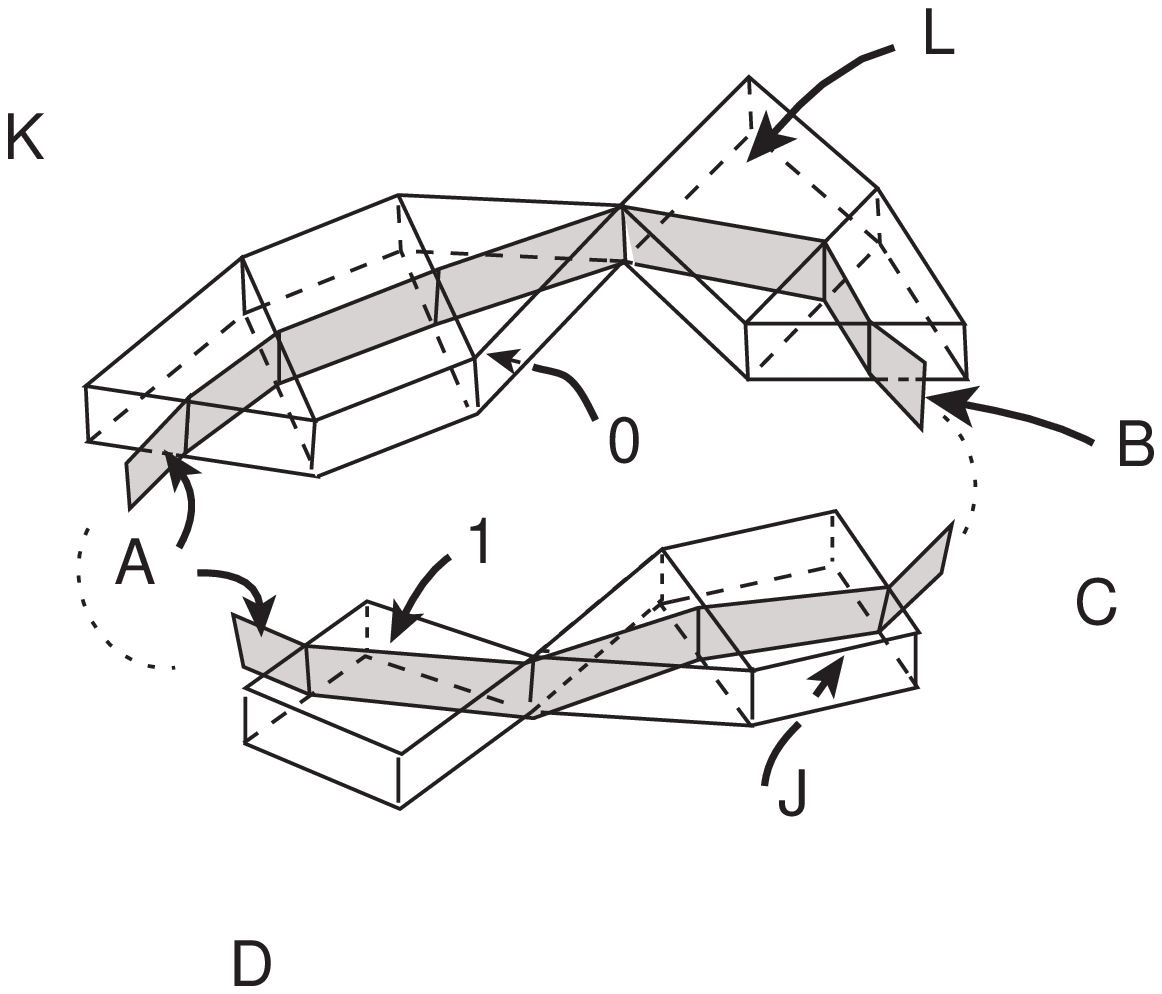}
        \caption{A $0$-weight annulus, $A$, in cells of type III and IV.}
        \label{f-annulus}
        \end{center}

\end{figure}

This seems good; however, again we have a complication here that
did not occur in the proof of Theorem \ref{0-eff-exists}; however,
we did see a similar situation in the proof of Theorem
\ref{b-0-eff-exists}. In a $2$--sphere, if we have two disjoint,
connected subcomplexes, then we can find disjoint simply connected
subcomplexes containing them. In a surface of higher genus, even
if  both subcomplexes are connected and contained in simply
connected subcomplexes, this may not be the case; one subcomplex
may be contained in the simply connected subcomplex determined by
the other.

So, suppose $K_i$ is not simply connected, $K_i^\ve\subset E_j,
\ve = 0,1$. We have shown that $K_i^\ve, \ve =0,1$ are contained
in simply connected regions of $E_j$. Thus there are simply
connected regions of the complement of $K_i^\ve$ in $S_j'$ so that
we can add them to $K_i^\ve$ and get simply connected subcomplexes
$D_i^\ve\supset K_i^\ve, \ve = 0,1$. We claim that by having
chosen the collection $E_1,\ldots,E_n$ maximal, we must have that
$K_i^0$ is not in $D_i^1$ and $K_i^1$ is not in $D_i^0$. For
suppose $K_i^1$, say, is in $D_i^0$. Just as in the proof of
Theorem \ref{b-0-eff-exists}, let $N_i = N(D_i^0\cup (K_i\times
[0,1]))$ be a small regular neighborhood of $D_i^0\cup (K_i\times
[0,1])$. Since $K_i^1\subset D_i^0$, the frontier of $N_i$
consists of an annulus, possibly some $2$--spheres and a disk
$F_i$, properly embedded in $X$ and having its boundary in $E_j$.
Since $M$ is $\bdy$--irreducible ($E_j$ is incompressible), there
is a disk $F_i'$ in $E_j$ so that $\bdy F_i' = \bdy F_i$ and
$F_i\cup F_i'$ bounds a $3$--cell, say $B$, in $X$. Furthermore,
$B$ contains $N_i$. The boundary of a small regular neighborhood
of $P_j\cup B$, say $E_j'$, along with each of the $E_k, k\ne j$
form a barrier surface. Now, if we shrink $E_j'$, we arrive at a
contradiction to the collection $E_1,\ldots,E_n$ being maximal
just as in the proof of Theorem \ref{b-0-eff-exists}. Again, see
Figure \ref{f-ideal-pushout}.

To meet the conditions of Theorem \ref{crush}, we now need to show
that there are product regions $\bbb{P}(X)$ and they are trivial.
We have that each component of $\bbb{P}(\C)$ is a product
$I$--bundle, $K_i\times [0,1]$, and for each $i, 1\le i\le k$,
where $k$ is the number of components of $\bbb{P}(\C)$.
Furthermore, we have shown, there are simply connected
subcomplexes $D_i^\ve, \ve =o,1$ in the collection $E_1,\ldots,
E_n$ so that $K_i^\ve\subset D_i^\ve$ and $D_i^0\cap D_i^1 =
\emptyset$. So, just as so many times above, to show that the
product regions $\bbb{P}(X)$ exist and are trivial, we need to
show that for $N_i = N(D_i^0\cup(K_i\times [0,1])\cup D_i^1)$ a
small regular neighborhood of $D_i^0\cup(K_i\times [0,1])\cup
D_i^1$, then $N_i$ has $3$--cell plugs. However, this is straight
forward as $M$ is irreducible and each of the surfaces $E_j$ is
parallel into a component of $\bdy M$.

So, by choosing some order, say $K_1\times [0,1],\ldots, K_k\times
[0,1]$ for the components of $\bbb{P}(\C)$, we can construct the
trivial product regions $D_j\times [0,1], 1\le j\le k'\le k$ so
that $\bigcup (K_i\times [0,1])\subset \bigcup (D_j\times [0,1])$
and the frontier of $\bigcup (D_j\times [0,1])$ is contained in
the frontier of $\bigcup (K_i\times [0,1])$. We let $\bbb{P}(X)$
be the components of $\bigcup (D_j\times [0,1])$.

This completes the proof of the Claim.

\vspace{.15 in} \noindent {\it Claim. There is no cycle of
truncated prisms in $X$, which is not in $\bbb{P}(X)$.} \vspace{.1
in}

\noindent {\it Proof of Claim}. Just as above, there is the
possibility of two types of cycles of truncated prisms: one is a
cycle about an edge $e$ of $\T$ (see Figure \ref{f-cycle}(A)) and
the other cycles about more than one edge of $\T$ (see Figure
\ref{f-cycle}(B)).

If there is a complete cycle about an edge $e$ as in Figure
\ref{f-cycle}(A), then one of the surfaces $E_j$ contains a thin
tube of elementary quads about the edge $e$. In this case, there
is a properly embedded disk $D$ in $X$ meeting $e$ in precisely
one point and meeting $E_j$ in $\bdy D$. A surgery on $E_j$ at $D$
gives two {\it normal} surfaces. One is a normal $2$--sphere,
since $\bdy M$ in $\bdy$--irreducible ($E_j$ is incompressible),
and the other is isotopic with $E_j$. However, this contradicts
the choice of the collection $\mathcal{E} = \{E_1,\ldots, E_n\}$
being maximal. See Figure \ref{f-surgery-cycle}.  It follows that,
having chosen $\mathcal{E}$  maximal, there is no complete cycle
of truncated prisms in the induced cell structure on $X$ about a
single edge.

If there is a complete cycle about more than one edge (see Figure
\ref{f-cycle}(B)), then, as above, the collection of cells of type
II (truncated prisms), form a solid torus with, possibly, some
self identifications in its boundary. Again, because of the
possible singularities, we distinguish between the cycle of
truncated prisms, which we denote $\hat{\tau}$, and the cycle of
truncated prisms minus the bands of trapezoids, which we denote by
$\tau$. We have that $\tau$ meets the surfaces in the collection
$\mathcal{E}$ either in three open annuli, each meeting a
meridional disk of $\tau$ precisely once, or a single open
annulus, meeting a meridional disk of $\tau$ three times. See
Figure \ref{f-cycle-torus}.

In the case we have three annuli common to $\tau$ and the surfaces
in $\mathcal{E}$, we conclude as in the argument for Theorem
\ref{b-0-eff-exists} that the cycle of truncated prisms must be in
the induced product region $\bbb{P}(X)$.

So, consider the case when we have just one annulus, say $A_1$, in
$\tau\cap E_j$. Then $A_1$ meets each hexagonal face in the chain
of truncated prisms, $\hat{\tau}$, three times. If the annulus
$A_1$ is trivial in $E_j$ (i.e., the core curve in $A_1$ is
contractible in $E_j$), then there would be an $L(3,1)$ as a
connected summand of $M$. But this is impossible as $M$ is
irreducible and $\bdy M\ne \emptyset$. Thus, we must have that the
annulus $A_1$ is not trivial in $E_j$. Again, by possibly
shrinking the torus $\hat\tau$, we have an embedded torus,
$\overline\tau$, $\overline\tau$ meets $E_j$ in an annulus
$A\subset\bdy\overline\tau$ ($A$ is just the, possibly singular,
annulus $A_1$ slightly shrunk and is embedded), which meets the
meridional disk of $\overline\tau$ exactly three times. Let $A'$
denote the closure of the complement of $A$ in $\bdy\overline\tau$
($A'$ is the band of trapezoids, $A_1'$, in $\hat{\tau}$ pulled
slightly into the truncated prisms and is embedded). Then since
$M$ can not be a solid torus, we have that $A'$ is an essential
annulus in $X$, which contradicts that $M$ is annular. This
completes the proof of the Claim.

Thus we have shown that the conditions in the hypothesis of
Theorem \ref{crush} are satisfied and  $\T^*$ is an ideal
triangulation of $\open{X}$. The ideal triangulation $\T^*$ has a
distinct ideal vertex for each component $S_i$ of $\bdy M$; its
index is the genus of $S_i$. Also, notice that every tetrahedron
of $\T$ which meets a $P_i$ in a truncated tetrahedron or meets an
$E_i$ in a quadrilateral piece, is  crushed and does not occur in
the ideal triangulation.

This completes the proof of the Theorem. \end{proof}

The next theorem shows that under the same hypothesis, we can get
$0$--efficient ideal triangulations for the interiors of these
manifolds.

\begin{thm} Suppose $M$ is a compact, irreducible, $\bdy$--irreducible, anannular
$3$--manifold. Then any ideal triangulation of
$\stackrel{\circ}{M}$ can be modified to a $0$--efficient ideal
triangulation of $\stackrel{\circ}{M}$.
\end{thm}

\begin{proof} Assume $\T$ is an ideal triangulation of
$\stackrel{\circ}{M}$. If $\T$ is
$0$--efficient, there is nothing to prove; so, we assume $\T$ is
not $0$--efficient. Hence, there is a normal $2$--sphere $\Sigma$
embedded in $\stackrel{\circ}{M}$. Now, for some normal
vertex-linking surface $S$ we have an arc $\Lambda$ in the
$1$--skeleton of $\T$, which meets $\Sigma$ in a single point,
meets $S$ in a single point and does not meet any of the other
vertex-linking surfaces. Since $M$ is irreducible, $\Sigma$ bounds
a $3$--cell $B$ in $M$.

Let $N$ be a small regular neighborhood of $S\cup\Lambda\cup B$.
Then one component of the boundary of $N$ is a copy of $S$ but the
other is a surface $S'$ isotopic to $S$ and, along with other
vertex-linking surfaces, a barrier surface in the component of its
complement not meeting $S\cup\Lambda\cup B$. We can shrink $S'$.
Again, as above, since $M$ is irreducible and $\bdy$--irreducible,
we arrive at a normal surface $S''$, which is isotopic with $S'$
and hence, isotopic with $S$, and possibly some normal
$2$--spheres along with some $0$--weight $2$--spheres embedded
entirely  in the interior of tetrahedra. We discard any such
$2$--spheres. Since $\Sigma$ must contain a quadrilateral and $M$
is not an $I$--bundle ($M$ is anannular), the normal surface $S''$
is not vertex-linking. It follows that $S''$ contains a
quadrilateral. The surface $S''$ is isotopic with the
vertex-linking surface $S$ and so $S''$ bounds a product $P =
S\times [0,1)$ containing $S$ and with $S'' = S\times 0$.

Using the same techniques as in the proof of the previous theorem,
we crush the ideal triangulation along the normal surface $S''$.
Actually, we can think of this as a crushing along the normal
surface $S''$ along with all the vertex-linking surfaces distinct
from $S$. We arrive at a new ideal triangulation $\T'$ of the
interior of a manifold $M'$ homeomorphic with $M$. However, we
have that the triangulation $\T'$ has strictly fewer tetrahedra
than the triangulation $\T$. It follows that after possibly some
further, but finite number of steps, we have an ideal
triangulation having no normal $2$--spheres and so is
$0$--efficient.\end{proof}

\begin{cor}  Suppose $M$ is a compact, irreducible, $\bdy$--irreducible, anannular
$3$--manifold. Then a minimal
ideal triangulation of $\stackrel{\circ}{M}$ is  $0$--efficient.
\end{cor}

\begin{proof} If an ideal triangulation of  $\stackrel{\circ}{M}$
is not $0$--efficient, then by
the proof of the previous theorem we can find another ideal
triangulation having strictly fewer tetrahedra. However, if our
ideal triangulation is minimal, this is impossible. It follows
that a minimal ideal triangulation must be
$0$--efficient.\end{proof}

 We observed above that there are several known algorithms for constructing ideal
triangulations of the complements of links in the $3$--sphere. We
also use different techniques in \cite{jac-let-rub1} and get
results similar to Theorem \ref{i-exist} in more general
circumstances. The methods used here lead to an algorithm for
constructing ideal triangulations of irreducible,
$\bdy$--irreducible, anannular $3$--manifolds. Furthermore, there
is an algorithm to decide if an ideal triangulation is or is not
$0$--efficient and, if it is not $0$--efficient, the algorithm
will construct a normal $2$--sphere. It follows there is an
algorithm to construct $0$--efficient ideal triangulations for
irreducible, $\bdy$--irreducible and anannular $3$--manifolds.
Finally, as an aside, we note
 it can be decided if a
$3$--manifold is irreducible, or  if it is $\bdy$--irreducible, or
if it is anannular, see \cite{haken1, sch, jac-tol, jac-ree, rub,
tho}.

We end this section with the following observation about the ideal
triangulations we have constructed.

\begin{thm} Suppose $M$ is a compact, irreducible, $\bdy$--irreducible, anannular
$3$--manifold. Then  $\stackrel{\circ}{M}$ has an ideal
triangulation which is $0$--efficient and the only closed normal
surface in a collared neighborhood of a vertex is vertex-linking.
\end{thm}

\begin{proof} We know such a $3$--manifold $M$ has an ideal
triangulation, say $\T$ and we may as well assume this
triangulation is $0$--efficient.

Now, suppose $S$ is a vertex-linking surface and in a collared
neighborhood of the vertex, $v_S$ at $S$, there is a closed normal
surface $F$. Then there is a collared neighborhood $N$ of $v_S$ so
that $F\subset \stackrel{\circ}{N}$ and the boundary $S'$ of $N$
is isotopic to $S$; furthermore, $S'$ is a barrier surface in the
component of its complement not meeting $N$. Now, we shrink $S'$.
If $F$ is not a vertex-linking surface, then we arrive at a normal
surface $S''$ which is isotopic to $S$, is not vertex-linking
(contains a quadrilateral) and bounds a product $P = S\times
[0,1), S'' = S\times \{0\}$ and $F\subset P$. We wish to crush the
triangulation along the surface $S''$.  By following the steps,
exactly as above, we are able to show that the conditions of
Theorem \ref{crush} are satisfied and we arrive at a new ideal
triangulation $\T'$ of the interior of a $3$--manifold $M'$ which
is homeomorphic with $M$. Furthermore, since $S''$ contains a
quadrilateral, the ideal triangulation $\T'$ has strictly fewer
tetrahedra than $\T$.

Now, it might be the case we have introduced some normal
$2$--spheres. But if we do then by using the same methods, we can
get still another ideal triangulation with even fewer tetrahedra.
Since there were only a finite number of tetrahedra in $\T$ to
start with, the process must terminate in the desired
triangulation. \end{proof}

\section{$0$--efficient triangulations and irreducible knots.}

In \cite{bin3} it is shown that for any closed, orientable
$3$--manifold $M$ there is a knot $K$ embedded in $M$ so that the
complement of $K$, $M\setminus K$, is irreducible. We shall say a
knot in a closed $3$--manifold $M$ is an {\it irreducible knot} if
its complement in $M$ is an irreducible $3$--manifold. Theorem
\ref{0-eff-construct} provides a constructive method for finding
such knots in irreducible $3$--manifolds.

\begin{thm} Given the closed, irreducible $3$--manifold $M$, $M\neq\rp$, a
one-vertex triangulation $\T$ of
$M$ may be constructed so that every edge of $\T$ (a knot in $M$) is an
irreducible knot.
\end{thm}
\begin{proof} By Theorem \ref{0-eff-construct}, given any triangulation of $M$,
we can modify that triangulation to a $0$--efficient triangulation
or we can show the manifold is $S^3, \rp$ or $L(3,1)$. For $S^3$
and $L(3,1)$, we also have one-vertex, $0$--efficient
triangulations. See Figures \ref{f-tetra} (5) and Figure
\ref{f-two-L3_1} A, respectively. Hence, except for $\rp$ we can
construct a $0$--efficient triangulation, with just one vertex.

 Now, if
$\T$ is a one-vertex triangulation, then any edge, $e$, in $\T$ is
a knot in $M$. Furthermore, if $\T$ is also $0$--efficient, then
for any edge $e$, $M\setminus\{e\}$ is irreducible. For otherwise,
there would be an essential $2$--sphere in $M\setminus\{e\}$;
however, the boundary of a small neighborhood of $e$ acts as a
barrier surface and so such a $2$--sphere would lead to a normal
$2$--sphere in $M$, rel $\T$, which is not vertex-linking (it
misses the edge $e$). This contradicts $\T$ being $0$--efficient.
\end{proof}

We do not know if $\rp$ has a one vertex triangulation in which
every edge is an irreducible knot; however, in the two-tetrahedra,
one-vertex triangulation of $\rp$ in Figure~
\ref{f-onevertex-RP3-S2xS1}, there are three edges, two are
irreducible knots and the third is not irreducible; it misses an
embedded $\rpp$. On the other hand, this third edge bounds a disk,
is a trivial knot and so, could not be irreducible in any manifold
except $S^3$. Anyhow, for all irreducible $3$--manifolds, there is
an algorithm to construct irreducible knots; our algorithm gives
them as edges in nice triangulations. This is somewhat a basic
theme of these methods and triangulations; namely, conditions on
the triangulations give edges of the triangulations, which provide
special properties for the associated knot complements.

\bibliographystyle{plain}
\bibliography{0-efficient-new}

\end{document}